\documentclass[10pt]{article}
\usepackage{amssymb,amsmath,amsthm}
\usepackage{color}
\usepackage{graphicx}
\usepackage{float}
\usepackage{caption}
\newcommand{\R}{\mathbb{R}}
\UseRawInputEncoding

\begin{document}
\title{\Large\bf{ Multiplicity and uniqueness of positive solutions for a superlinear-singular $(p,q)$-Laplacian equation on locally finite graphs}}
\date{}
\author {  Xuechen Zhang,$^{1}$  Xingyong Zhang$^{1,2}$\footnote{Corresponding author, E-mail address: zhangxingyong1@163.com}\\
       {\footnotesize $^1$Faculty of Science, Kunming University of Science and Technology,}\\
       {\footnotesize Kunming, Yunnan, 650500, P.R. China.}\\
       {\footnotesize $^2$Research Center for Mathematics and Interdisciplinary Sciences, Kunming University of Science and Technology,} \\
       {\footnotesize Kunming, Yunnan, 650500, P.R. China.}}

 \date{}
 \maketitle

 \begin{center}
 \begin{minipage}{15cm}
 \par
\small  {\bf Abstract:} We investigate the multiplicity and uniqueness of positive solutions for the superlinear singular $(p,q)$-Laplacian equation
\begin{eqnarray*}
 \begin{cases}
-\Delta_p u-\Delta_q u+a(x)u^{p-1}+b(x)u^{q-1}=f(x)u^{-\gamma}+\lambda g(x)u^{\alpha}, \;\;\;\;\hfill \mbox{in}\;\; V,\\
u>0,\;\;u\in W_a^{1,p}(V) \cap W_b^{1,q}(V),
 \end{cases}
\end{eqnarray*}
on a weighted locally finite graph $G=(V,E)$, where $0<\gamma<1<q\leq p<\alpha+1$, $\lambda$ is a parameter, the potential functions $a(x)$ and $b(x)$ satisfy some suitable conditions, $f>0, g \geq 0$, $f\in L^1(V)\cap L^{\frac{p}{p-1+\gamma}}(V) \cap L^{\frac{q}{q-1+\gamma}}(V)$ and $g\in L^1(V)\cap L^\infty(V)$. By making use of the method of Nehari manifold and the Ekeland's variational principle, we prove that there exist two positive solutions for $\lambda$ belonging to some precise interval. Besides, we also investigate the existence and uniqueness of positive solution for $\lambda<0$.
We overcome some difficulties which are caused by: $(i)$ the singular term; $(ii)$ the definition of gradient $|\nabla u|$ on graph which is different from that on $\R^N$; $(iii)$ the lack of compactness of Sobolev embedding.
\par
{\bf Keywords:}  $(p,q)$-Laplacian equation, singular term, locally finite graph, Ekeland's variational principle, Nehari manifold.
 \par
{\bf 2020 Mathematics Subject Classification.} 35A02; 35B38; 35J75; 35R02
 \end{minipage}
 \end{center}
  \allowdisplaybreaks
 \vskip2mm
 {\section{Introduction }}
\setcounter{equation}{0}
In this paper, we are concerned about the existence and multiplicity of positive solutions for the following singular $(p,q)$-Laplacian equation
\begin{eqnarray}
\label{eq1}
 \begin{cases}
-\Delta_p u-\Delta_q u+a(x)u^{p-1}+b(x)u^{q-1}=f(x)u^{-\gamma}+\lambda g(x)u^{\alpha}, \;\;\;\;\hfill \mbox{in}\;\; V,\\
u>0,\;\;u\in W_a^{1,p}(V) \cap W_b^{1,q}(V),
 \end{cases}
\end{eqnarray}
on a weighted locally finite graph $G=(V,E)$, where $\Delta_s u, s=p,q$ is the discrete $s$-Laplacian defined by (\ref{Eq9}) below on graphs, $0<\gamma<1<q\leq p<\alpha+1$, $f> 0$, $g\geq0$, $g\not\equiv0$, $f\in L^1(V)\cap L^{\frac{p}{p-1+\gamma}}(V) \cap L^{\frac{q}{q-1+\gamma}}(V)$ and $g\in L^1(V)\cap L^\infty(V)$, $a(x)$ and $b(x)$ are positive continuous function and satisfy the following condition:\\
$(A)$\;\;there exists positive constants $a_0$ and $b_0$ such that $a(x)\geq a_0>0$ and $b(x)\geq b_0>0$ for all $x\in V$.
\par
Graphs represent a natural framework that is particularly well-suited to model various phenomena observed in real-world scenarios.
Some analysis on graphs have been applied to the investigation of machine learning, data analysis, neural network, image processing etc (for example, see \cite{Alkama S2014, Arnaboldi V 2015, Elmoataz A2015, Ta V T 2010, Ta V T 2008}). Recently, some researchers begin to study the partial differential equations on graphs. In \cite{Grigor'yan 2017}, Grigor-Lin-Yang studied the following nonlinear equation with $p=2$ on a locally finite graph $G=(V,E)$:
\begin{eqnarray}
\label{c3}
 -\Delta_p u+h(x)|u|^{p-2}u=f(x,u), \;\;\mbox{$x \in V$},
\end{eqnarray}
where $p>1$, $h: V\rightarrow \R$ and $f: V\times \R \rightarrow \R$. Under superquadratic conditions on $f$,
they demonstrated the existence of strictly positive solutions for equation (\ref{c3}). In \cite{Chang logarithmic 2023}, Chang-Wang-Yan considered (\ref{c3}) with $p=2$ and  $f(x,u)= u\log u^2$. They employed the Nehari manifold method and the mountain pass theorem to establish the existence of ground state solutions. By setting $p=2, h(x)= \lambda a(x)+1$ and $f(x,u)= |u|^{q-2}u$ (where $\lambda>1$ and $q>2$), and incorporating the term $\Delta^2 u$ to the left-hand side of equation (\ref{c3}), the resulting equation is referred to as the biharmonic equation. When $a$ is a potential function satisfying the coercive condition, Han-Shao-Zhao in \cite{Han 2020} studied the biharmonic equation by the method of Nehari manifold. They proved the existence of a ground state solution and showed that this solution converges to the ground state solution of the corresponding Dirichlet problem as $\lambda \rightarrow +\infty$. In \cite{Yang P 2023}, Yang-Zhang investigated a $(p,q)$-Laplacian coupled system. When $p=q$ and $u=v$, the system they studied reduces to equation (\ref{c3}). They generalized two Sobolev embedding theorems and, utilizing the mountain pass theorem and Ekeland's variational principle, obtained results on the existence and multiplicity  when the nonlinear term satisfies the sub-$(p,q)$ conditions and super-$(p,q)$ conditions, respectively. Pan-Ji in \cite{Pan 2023} extended equation (\ref{c3}) to the nonlinear Kirchhoff equations, and they obtained the existence of a least energy sign-changing solution via constrained variational method, and proved that the energy of the least energy sign-changing solution is strictly larger than twice that of the least energy solutions. Additionally, they set up the convergence property of the least energy sign-changing solution. More results can be referred to \cite{Hu 2024, Ou 2024, Zhang  2024}. However, up to now, there are few results on existence of positive solutions for the singular $(p,q)$-Laplacian equation like (\ref{eq1}) on a weighted locally finite graph.
The $(p,q)$-Laplacian problem (\ref{eq1}) arises from a general reaction-diffusion system:
\begin{eqnarray*}
u_t=\mbox{div}[D(u)\nabla u]+c(x,u),\;\; x\in \R^N.
\end{eqnarray*}
This system finds extensive applications across physical and related scientific disciplines, including biophysics \cite{Fife 2013,Schneider 1989}, plasma physics \cite{Herrera 1992,Smoller 2012} and chemical reaction design \cite{Aris 1994}. In these applications, the function $u$ describes a concentration, $\mbox{div}[D(u)\nabla u]$ denotes the diffusion process with a potentially non-constant diffusion coefficient $D(u)$ and $c(x,u)$ corresponds to the reaction and relates to source and loss processes. Typically, in chemical and biological applications, the reaction term $c(x,u)$ takes on a polynomial form with respect to the concentration $u$.
\par
In our case, we discuss weak solutions of (\ref{eq1}) with a diffusion coefficient having a power law dependency $D(u) = (|\nabla u|^{p-2}+ |\nabla u|^{q-2})$ and reaction terms $c(x,u)= -a(x)|u|^{p-2}u-b(x) |u|^{q-2} u +f(x)u^{-\gamma}+\lambda g(x)u^\alpha$. The $(p,q)$-Laplacian equation in the Euclidean setting has been noticed in the last years,  and the existence, multiplicity, nonexistence and concentration of weak solutions were studied in \cite{Ambrosio 2021, Baldelli 2022, Cherfils 2005, Figueiredo 2011, Li 2009} and the references therein.
\par
Besides, the nonlinear term in (\ref{eq1}) is singular at $u=0$, i.e. $f(x)u^{-\gamma} \rightarrow \infty$ as $u \rightarrow 0$. The singular problem were considered in the differential equation for the first time by Fulks-Maybee \cite{Fulks 1960}, where they illustrated the physical problem on a region $\Omega$ of $\R^3$ occupied by an electrical conductor. $u(x,t)$ is the temperature at the point $x\in \Omega$ and at time $t$, and it satisfies the equation
\begin{eqnarray*}
c u_t- \kappa\Delta u=\frac{E^2(x,t)}{\sigma(u)},\;\;x\in\Omega.
\end{eqnarray*}
$E(x,t)$ represents the local voltage drop in $\Omega$ as a function of both position and time. The parameters $c$ and $\kappa$ denote the specific heat and thermal conductivity, respectively. Generally, $\sigma(u)$ can be characterized as a a positive function of $u$ which increases and tends to zero as $u$ approaches zero, exemplified by expressions like $u^\gamma$ with $\gamma>0$. Then the differential equation is singular at $u = 0$. The existence and multiplicity of solutions to these singular problems have been studied extensively. In \cite{Sun Y 2001}, Sun-Wu-Long considered the following equation with $W(x)=1$:
\begin{eqnarray}
\label{c4}
 \begin{cases}
 \begin{array}{ll}
-\Delta u + \lambda W(x) u^\beta+P(x)u^{-\gamma}=0,\;\;&\mbox{in $\Omega$,}\\
u>0,\;\;\;\;\;\;\;\;\;\;\;\;\;\;\;\;\;\;\;\;\;\;\;\;\;\;\;\;\;\;\;\;\;\;\;\;\;\;
&\mbox{in $\Omega$,}\\
u=0,\;\;\;\;\;\;\;\;\;\;\;\;\;\;\;\;\;\;\;\;\;\;\;\;\;\;\;\;\;\;\;\;\;\;\;\;\;\;
&\mbox{on $\partial\Omega$,}
 \end{array}
 \end{cases}
\end{eqnarray}
where $\Omega \subset \R^N$ is a bounded domain, $1<\beta<2^*-1, 2^*=\frac{2N}{N-2},0<\gamma<1$, $\lambda>0$ has an appropriate range and $P: \Omega \rightarrow \R$ is a given non-negative non-trival function. They decomposed the solution set into three submanifolds, with the definition of the submanifolds being related to $\beta$ and $\gamma$. Subsequently, they obtained the existence and multiplicity results of weak solutions by the method of Nehari manifold and the Ekeland's variational principle, with $P\in L^2(\Omega)$. Specifically, the energy value of one solution is less than zero and is a local minimum, whereas the energy value of the other solution is greater than or equal to zero and is a minimum on the Nehari submanifold. In \cite{Sun 2008}, Sun-Li studied (\ref{c4}) with $P(x)\geq 0, P \not \equiv0$, $W$ being a given function with the set $\{x \in \Omega:W(x) >0 \}$ of positive measures, and $P, W\in C(\overline{\Omega})$. The primary objective of this paper is to provide an explicit estimate for the lower bound of the critical parameter $\lambda$. It offers a detailed formula that correlates with $\gamma$, $p$ and $\Omega$. Furthermore, the paper broadens the scope of $\lambda$ beyond what was discussed in \cite{Sun Y 2001}, demonstrating that equation (\ref{c4}) admits a positive solution within each of the two non-zero submanifolds. It also establishes that the energy values of these two positive solutions are minimized on their respective submanifolds.
There are also a few studies on (\ref{c4}) with $\Omega=\R^N$. In \cite{Alves  2022}, Alves-Ricardo-Kaye considered the following equation:
\begin{eqnarray}
\label{c5}
-\Delta u + V(x)u= \lambda A(x) u^{-\gamma}+ B(x)u^{p},\;\;\mbox{in $\R^N$},
\end{eqnarray}
where $1<p<2^*-1, \lambda>0$, $0<A\in L^{\frac{2}{1+\gamma}}(\R^N)$, $B\in L^\infty(\R^N)$ may change sign and $V$ is a positive continuous function. They obtained the multiplicity and non-existence results of negative-energy solutions by Nehari manifold method. Besides, they also characterized a $\lambda$-behavior of the energy functional of the solutions. For more references on singular problems, one can see \cite{Corra 2023, Goncalves 2007, Jose Valdo 2007, Liu 2009}.
\par
To the best of our knowledge, there appears to be no prior research investigating the existence of solutions for a singular Laplacian equation on locally finite graphs. Inspired by \cite{Alves  2022, Sun Y 2001, Sun 2008}, in this paper, we consider whether the Ekeland's variational principle can be applied to the investigation on the multiplicity and uniqueness of positive solutions for equation (\ref{eq1}) involving superlinear-singular term on locally finite graphs. There are three key steps in our proofs: $(1)$ Choose the suitable working space. Due to the singular term, the energy functional is not differentiable, preventing the direct application of classical critical point theory. We circumvent this difficulty by constraining the problem to the cone of non-negative functions $W_+$ (see (\ref{c6})) and using Ekeland's variational principle in conjunction with compactness analysis. Furthermore, we decompose the Nehari manifold $D_{\lambda}\subset W_+$ into three submanifolds $D_{\lambda}^+, D_{\lambda}^-$ and $D_{\lambda}^0$ when the parameter $\lambda$ has an appropriate range. We prove that  the sets   $D_{\lambda}^+, D_{\lambda}^-$ are non-empty, while the set $D_{\lambda}^0$ contains only the zero element (Lemma 2.4).
$(2)$ Prove that each of the two submanifolds $D_{\lambda}^+, D_{\lambda}^-$  contains one solution. At first, utilizing Ekeland's variational principle, we construct minimizing sequences $\{u_k\} \subset D_{\lambda}^+$ and $\{v_k\} \subset D_{\lambda}^-$, respectively. Through a proof by contradiction combined with certain calculations, we demonstrate that the sequence $\{u_k\}$ converges strongly to a $u_{\lambda} \in D_{\lambda}^+$, and by leveraging the continuous embedding property (Lemma 2.1) and the Brezis-Lieb Lemma (Lemma 2.2), we obtain the sequence $\{v_k\}$ converges strongly to a $v_{\lambda} \in D_{\lambda}^-$.
Moreover, it is proved that the submanifolds $D_{\lambda}^+$ and $D_{\lambda}^-$ are closed (Lemma 2.6), thereby ensuring that these two solutions are non-coincidence.
$(3)$ Prove the uniqueness of the solution for $\lambda<0$.
 Since the energy functional is coercive and weakly lower semi-continuous if  $\lambda<0$,  then it follows from  the least action principle that there exists a global minimizer $\omega_\lambda$. Furthermore, with the help of the inequality (\ref{s7}), we can establish the uniqueness of the solution.
\par
  Our work develops those results in \cite{Alves  2022, Sun Y 2001, Sun 2008} from the following aspects:\\
  (I) \;  Different from \cite{Alves 2022, Sun Y 2001, Sun 2008}, our work is conducted on a locally finite graph $V$, rather than on a bounded domain $\Omega\subset \R^N$ or $\R^N$. Furthermore, we concentrate on the more general quasilinear equation (\ref{eq1}) featuring a $(p,q)$-Laplacian operator $(1<q\leq p)$, instead of equations (\ref{c4}) and (\ref{c5}) with the classical Laplacian operator.
\vskip2mm
\noindent
(II) \; Our proof strategy is not entirely identical to that developed in \cite{Sun Y 2001, Sun 2008}, which stems from two fundamental aspects: (i)  the quasilinear nature  of the $(p,q)$-Laplacian operator $-\Delta_p-\Delta_q$, and  (ii) the unique definition of the gradient $\nabla u$ on graph, which differs significantly from its Euclidean counterpart in $\R^N$. This distinction is concretely reflected in the following three aspects:
\begin{itemize}
 \item The equality $J_\lambda(|u|)=J_\lambda(u)$, as established in \cite{Sun Y 2001, Sun 2008}, cannot be achieved in our case, and consequently, the non-negativity assumption
$u\geq 0$ employed in \cite{Sun Y 2001, Sun 2008} is no longer applicable. Inspired by the work of \cite{Alves 2022}, we resolve this problem by introducing the set
$W_+$(see (\ref{c6})) and establishing the existence of positive solutions in $W_+$, instead of the set
$W$  adopted in \cite{Sun Y 2001, Sun 2008}.
 \item The approach we employ to demonstrate the formulas (\ref{x21}) and (\ref{x22}) diverges from the methodology presented in \cite{Alves 2022, Sun Y 2001, Sun 2008}. Since the analysis in \cite{Alves 2022, Sun Y 2001, Sun 2008} was conducted on the domain $\Omega\subseteq \R^N$ and the  compactness of the Sobolev embedding $W^{1,2}(
     \Omega)\hookrightarrow L^{\theta}(\Omega)$ $(2\le \theta<2^*)$ hold, they were able to deduce the validity of formulas (\ref{x21}) and (\ref{x22}) by applying the Dominated Convergence Theorem straightforwardly. However, we consider the equation (\ref{eq1})  on a locally finite graph $V$ with $|V|$ being infinite and our assumption $(A)$ can not ensure the  compactness of the Sobolev embedding $W^{1,s}(V)\hookrightarrow L^{\theta}(V)$ $(s\le \theta<\infty)$, so that we are  unable to demonstrate these two formulas using the idea in \cite{Alves 2022, Sun Y 2001, Sun 2008}. Instead of relying on the  compactness of the Sobolev embedding, we just make the simple constraints  $f\in L^1(V)$ and $g\in L^1(V)$, and then by leveraging point-wise convergence that $u_k(x) \rightarrow u(x)$  in  for all $x\in V$ and the Dominated convergence theorem, we can effectively tackle this issue.
 \item Our approach  to prove the strong convergence of a minimizing sequence  $\{u_k\}$ in  $D_{\lambda}^+$ differs from that in \cite{Sun Y 2001, Sun 2008}. In \cite{Sun Y 2001, Sun 2008}, establishing such an inequality like (\ref{x38}) is crucial for demonstrating the strong convergence. Based on the Riesz representation theorem in Hilbert spaces, they simply derived the inequality  from the definition of weak convergence.
     However, our work is conducted in a Banach space where the Riesz representation theorem is not available. Consequently, we have to pursue an alternative strategy. We observe that, due to the discrete structure of the graph, it is obvious that $\nabla u_k(x)$ converges  to $\nabla u_{\lambda}(x)$ for all $x\in V$ if $u_k$ weakly converges to $u_{\lambda}$ as $k\to \infty$, which is not easy in the Euclidean setting in general. Then, in conjunction with Theorem 2 in \cite{Brezis 1983}, the Brezis-Lieb lemma about $|\nabla u_k|^q$ naturally holds on graph (see Lemma 2.2). Furthermore, inspired by the idea in \cite{Alves 2022},  by employing  the generalized Brezis-Lieb lemma on graphs, we can demonstrate the strong convergence without using the inequality (\ref{x38}) which is, however, still important to  prove that the limit $ u $ is a solution of the problem (\ref{eq1}). Since we have already demonstrated the strong convergence, the derivation of inequality (\ref{x38}) follows as a natural consequence.
\end{itemize}
\par
To describe our problems and results more clearly, we review  some concepts and assumptions as follows (more details can be seen in \cite{Grigor'yan 2017, Grigor'yan A 2016 yamabe, Grigor'yan A 2018 book}). Let $G=(V, E)$ be a graph and a measure $\mu:V\rightarrow \R^+$ be finite, where $V$ denotes the vertexes set and $E$ denotes the edges set. \\
$(i)$\;\;\;\; $G$ is called a locally finite graph: there are only finite $y \in V$ such that $xy \in E$ for any $x \in V$;\\
$(ii)$\;\;\; $G$ is connected: any two vertices $x$ and $y$ can be connected via finitely many edges;\\
$(iii)$\;\; $\mu$ is a uniformly positive measure: there exists a constant $\mu_{0}>0$ such that $\mu(x)\geq \mu_{0}$ for all $x\in V$;\\
$(iv)$\;\; $G$ is symmetric: $\omega_{xy}=\omega_{yx}$, for any edge $xy\in E$ with two vertexes of $x,y\in V$. Moreover assume that $\omega_{xy}>0$ and for any $x\in V$, $\sum_{y\thicksim x}\omega_{xy}$ is bounded by a positive constan $C$, where $\sum_{y\thicksim x}\omega_{xy}$ denotes the sum over all neighbors $y$ of $x$ with respect to the edge weights $\omega_{xy}$ and $y \thicksim  x$ stands for any vertex $y$ connected with $x$ by an edge $xy \in E$;\\
$(v)$\;\;\; $d(x,y)$: the distance of two vertices $x,y \in V$, that is, the minimal number of edges which connect these two vertices.

\par
For any $x\in V$, we define the Laplacian operator by $\Delta\psi:C(V)\rightarrow C(V)$ with
\begin{eqnarray*}
\label{Eq4}
\Delta \psi(x)=\frac{1}{\mu(x)}\sum\limits_{y\thicksim x}w_{xy}(\psi(y)-\psi(x)),
\end{eqnarray*}
where $C(V)$ is the set of all real functions on $V$.
Let
$$
D_{x,y}\psi(x):=\frac{1}{\sqrt{2}}(\psi(x)-\psi(y))\sqrt{\frac{\omega_{xy}}{\mu(x)}}
$$
be the directional derivative of a function $\psi: V\rightarrow \R$, and then the gradient  of $\psi$ is defined to be the vector
$$
\nabla \psi(x)=(D_{x,y}\psi(x))_{y\in V},
$$
which is indexed by the vertices $y\in V$.
Let
\begin{eqnarray*}
\label{Eq5}
\Gamma(\psi_1,\psi_2)(x):=\frac{1}{2\mu(x)}\sum\limits_{y\thicksim x}w_{xy}(\psi_1(y)-\psi_1(x))(\psi_2(y)-\psi_2(x))=\nabla \psi_1\nabla \psi_2.
\end{eqnarray*}
Denote $\Gamma(\psi)=\Gamma(\psi,\psi)$ and then the length of the gradient is defined by
\begin{eqnarray*}
\label{Eq6}
|\nabla \psi|(x)=\sqrt{\Gamma(\psi)(x)}=\left(\frac{1}{2\mu(x)}\sum\limits_{y\thicksim x}w_{xy}(\psi(y)-\psi(x))^2\right)^{\frac{1}{2}}.
\end{eqnarray*}
It is easy to verify that
\begin{eqnarray}
\label{x20}
|\nabla \psi| \geq |\nabla |\psi||.
\end{eqnarray}

For any function $\psi:V\rightarrow\mathbb{R}$, we denote
\begin{eqnarray*}\label{Eq8}
\int\limits_V \psi(x) d\mu=\sum\limits_{x\in V}\mu(x)\psi(x).
\end{eqnarray*}
When $s\geq2$,  the $s$-Laplacian operator is defined by $\Delta_s\psi:C(V)\rightarrow C(V)$ with
\begin{eqnarray}
\label{Eq9}
\Delta_s\psi(x)=\frac{1}{2\mu(x)}\sum\limits_{y\sim x}\left(|\nabla \psi|^{s-2}(y)+|\nabla \psi|^{s-2}(x)\right)\omega_{xy}(\psi(y)-\psi(x)).
\end{eqnarray}
In the distributional sense, $\Delta_s \psi$ can be written as follows. For any $\phi\in\mathcal{C}_c(V)$,
\begin{eqnarray*}\label{Eq10}
\int\limits_V(\Delta_s \psi)\phi d\mu=-\int\limits_V|\nabla \psi|^{s-2}\Gamma(\psi,\phi)d\mu,
\end{eqnarray*}
where $\mathcal{C}_c(V):=\{u:V\rightarrow \R:\{x \in V: u(x)\neq 0\} \mbox{ is of finite cardinality}\}$.

\par
When $1\leq\theta <+\infty$, we define
$$
L^{\theta }(V)=\left\{\psi:V\to\R\Big|\int_{V}|\psi(x)|^{\theta } d\mu<\infty\right\}
$$
with the norm
\begin{eqnarray*}
\label{b1}
\|\psi\|_{\theta }=\left(\int_{V}|\psi(x)|^{\theta } d\mu\right)^\frac{1}{{\theta }}.
\end{eqnarray*}
When ${\theta } =+\infty$, we define
$$
L^\infty(V)=\left\{\psi:V\to\R\Big|\sup_{x\in V}|\psi(x)|<\infty\right\}
$$
with the norm
\begin{eqnarray*}
\label{Z7}
\|\psi\|_{\infty}=\sup_{x \in V}|\psi(x)|< \infty.
\end{eqnarray*}

\vskip2mm
\par
Define
\begin{eqnarray*}
W^{1,s}(V)=\left\{u:V\to\R\Big|\int_V(|\nabla  u|^s+|u|^s)d\mu<\infty\right\}
\end{eqnarray*}
which is provided with the norm
\begin{eqnarray*}
\label{e1}
\|u\|_{W^{1,s}(V)}=\left(\int_V(|\nabla  u|^s+|u|^s)d\mu\right)^\frac{1}{s},
\end{eqnarray*}
where $s>1$, and define
\begin{eqnarray*}
W^{1,p}_a(V):=\left\{u\in W^{1,p}(V)\Big|\int_V a(x)|u|^pd\mu<\infty\right\}\;\;\;\;
\mbox{and}\;\;\;\;W^{1,q}_b(V):=\left\{u\in W^{1,q}(V)\Big|\int_V b(x)|u|^qd\mu<\infty\right\}
\end{eqnarray*}
with the following norms, respectively,
\begin{eqnarray*}
\|u\|_{W^{1,p}_a(V)}=\left(\int_V(|\nabla u|^p+a(x)|u|^p)d\mu\right)^{\frac{1}{p}}\;\;\;\;
\mbox{and}\;\;\;\;\|u\|_{W^{1,q}_b(V)}=\left(\int_V(|\nabla u|^q+b(x)|u|^q)d\mu\right)^{\frac{1}{q}}.
\end{eqnarray*}
\par
To study the problem (\ref{eq1}), it is natural to consider the function space
$W:=W_a^{1,p}(V)\cap W_b^{1,q}(V)$. Define the norm $\|u\|_{W}=\|u\|_{W^{1,p}_a(V)}+\|u\|_{W^{1,q}_b(V)}$.
Then, $(W,\|\cdot\|_{W})$ is a reflexive Banach space (see \cite{Shaomeng 2023}).
\par
The functional related to (\ref{eq1}) is defined by
\begin{eqnarray*}
\label{e2}
J_\lambda(u)=\frac{1}{p}\int_V(|\nabla u|^p+a|u|^p)d\mu+\frac{1}{q}\int_V(|\nabla u|^q+b|u|^q)d\mu-\frac{1}{1-\gamma}\int_V f(x) |u|^{1-\gamma}d\mu - \frac{\lambda}{\alpha+1}\int_Vg(x) |u|^{\alpha+1}d\mu.
\end{eqnarray*}
If $u\in W$ such that $u>0$ for all $x\in V$ and
\begin{eqnarray*}\label{c7}
&&\int_V(|\nabla u|^{p-2} \Gamma(u,\phi)+a(x)u^{p-1}\phi)d\mu
+\int_V (|\nabla u|^{q-2} \Gamma(u,\phi)+b(x)u^{q-1}\phi)d\mu\nonumber\\
&=&\int_Vf(x)u^{-\gamma} \phi d\mu+{\lambda}\int_V g(x)u^{\alpha}\phi d\mu
\end{eqnarray*}
for any $\phi \in W$, we call that  $u\in W$ is a weak solution of (\ref{eq1}) in $W$.
\par
Since we are devoted to finding the positive solutions. Motivated by \cite{Alves  2022}, we constrain $J_\lambda$ to the cone of non-negative functions in $W$, that is
\begin{eqnarray}\label{c6}
W_+=\{u\in W: u \geq 0  \}.
\end{eqnarray}
Define the Nehari manifold
$$
\mathcal{D}_\lambda:=\{v=t(u)u: u\in W_+,\varphi_u(t)=0\}
$$
where $t(u)$ are the zero point of the map
\begin{eqnarray*}
t\in (0,+\infty)\rightarrow \varphi_u(t)&:=&t^{\gamma}\frac{d}{dt}J_\lambda(tu)\nonumber\\
&=&t^{p-1+\gamma}\|u\|^p_{W^{1,p}_a (V)}+t^{q-1+\gamma}\|u\|^q_{W_b^{1,q}(V)}\nonumber-\int_V f(x)u^{1-\gamma}d\mu-\lambda t^{\alpha+\gamma}\int_V g(x)u^{\alpha+1}d\mu\\
&  =   &t^{-1+\gamma}\left[\|tu\|^p_{W_a^{1,p}(V)}+\|tu\|^q_{W_b^{1,q}(V)}
-\int_V f(x)(tu)^{1-\gamma}d\mu
-\lambda\int_V g(x)(tu)^{\alpha+1}d\mu\right].
\end{eqnarray*}
Furthermore,
\begin{eqnarray*}
\varphi'_{u}(t)
=(p-1+\gamma)t^{p-2+\gamma}\|u\|^p_{W^{1,p}_a (V)}+(q-1+\gamma)t^{q-2+\gamma}\|u\|^q_{W_b^{1,q}(V)}-(\alpha+\gamma)\lambda t^{\alpha+\gamma-1}\int_V g(x)u^{\alpha+1}d\mu.
\end{eqnarray*}
According to the definition of $\mathcal{D_\lambda}$, we know that the non-negative weak solution must lie in $\mathcal{D}_\lambda$. We split $\mathcal{D_\lambda}$ into $\mathcal{D}_\lambda^+ \cup \mathcal{D}_\lambda^0\cup \mathcal{D}^-_\lambda$, which are defined as follows:
\begin{eqnarray}
\mathcal{D}^+_\lambda  &   =    &\left\{v=t(u)u\in \mathcal{D}_\lambda: (p-1+\gamma)\|v\|^p_{W_a^{1,p}(V)}+(q-1+\gamma)\|v\|^q_{W_b^{1,q}(V)}-(\alpha+\gamma)\lambda\int_V g(x)v^{\alpha+1}d\mu>0\right\}\nonumber\\
&=&\{v=t(u)u\in \mathcal{D}_\lambda: \varphi'_{u}(t)>0\}\label{x1},\\
\mathcal{D}^0_\lambda  &   =    &\left\{v=t(u)u\in \mathcal{D}_\lambda: (p-1+\gamma)\|v\|^p_{W_a^{1,p}(V)}+(q-1+\gamma)\|v\|^q_{W_b^{1,q}(V)}-(\alpha+\gamma)\lambda\int_V g(x)v^{\alpha+1}d\mu=0\right\}\nonumber\\
&=&\{v=t(u)u\in \mathcal{D}_\lambda: \varphi'_{u}(t)=0\},\nonumber\\
\mathcal{D}^-_\lambda  &   =    &\left\{v=t(u)u\in \mathcal{D}_\lambda: (p-1+\gamma)\|v\|^p_{W_a^{1,p}(V)}+(q-1+\gamma)\|v\|^q_{W_b^{1,q}(V)}-(\alpha+\gamma)\lambda\int_V g(x)v^{\alpha+1}d\mu<0\right\} \nonumber\\
&=&\{v=t(u)u\in \mathcal{D}_\lambda: \varphi'_{u}(t)<0\}.\label{x2}
\end{eqnarray}
\vskip2mm
\noindent

\par At the end of this section, we present our main result. Let
\begin{eqnarray}
\label{w1}
  \begin{cases}
X(\lambda)
=\left( \frac{p-1+\gamma}{\lambda(\alpha+\gamma)\|g\|_\infty\mu_0^{\frac{p-1-\alpha}{p}}a_0^{-\frac{\alpha+1}{p}}}\right)
^{\frac{1}{\alpha+1-p}},\\
X_0=
 \left( \frac{\alpha+\gamma}{\alpha+1-p}
 \right)^{\frac{1}{p-1+\gamma}}a_0^{-\frac{1-\gamma}{p(p-1+\gamma)}}\|f\|^{\frac{1}{p-1+\gamma}}_{\frac{p}{p-1+\gamma}},\\
S(\lambda)
=\left( \frac{p-1+\gamma}{\lambda(\alpha+\gamma)} \right)^{\frac{1}{\alpha+1-p}}
\left(\frac{1}{\mu_0}\right)^{-\frac{1}{\alpha+1}} a_0^{\frac{1}{\alpha+1-p}} \left( \frac{1}{\|g\|_\infty} \right)^{\frac{1}{\alpha+1-p}}, \\
S_0
=\left( \frac{\alpha+\gamma}{\alpha+1-p}\right)^{\frac{1}{p-1+\gamma}} \left(\frac{1}{\mu_0}\right)^{\frac{\alpha+1-p}{p(\alpha+1)}} a_0^{-\frac{1}{p-1+\gamma}}
\|f\|^{\frac{1}{p-1+\gamma}}_{\frac{p}{p-1+\gamma}},\\
\Lambda_{*}
=\left(\frac{p+\gamma-1}{\alpha+\gamma}
\right)^{\frac{\alpha+\gamma}{p+\gamma-1}}
\left(\frac{\alpha+1-p}{p-1+\gamma}\right)^{\frac{\alpha+1-p}{p+\gamma-1}}\frac{1}{\|g\|_\infty}
\mu_0^{\frac{\alpha+1-p}{p}}a_0^{\frac{\alpha+\gamma}{p-1+\gamma}} \left(\frac{1}{\|f\|_{\frac{p}{p-1+\gamma}}}\right)^{\frac{\alpha+1-p}{p+\gamma-1}}.
\end{cases}
\end{eqnarray}

\vskip2mm
\noindent
{\bf Theorem 1.1.} {\it Let $G=(V, E)$ be a locally finite graph and it is connected, symmetric and there exists $\mu_{0}>0$ such that $\mu(x)\geq \mu_{0}$ for all $x\in V$. For any $x\in V$, $\sum_{y\thicksim x}\omega_{xy} <C$, where $C$ is a positive constant. Assume $0<\lambda < \Lambda_{*} $, and condition $(A)$ hold. Then equation (\ref{eq1}) has at least two positive solutions $u_\lambda$ and $v_\lambda$, with $u_\lambda$ belonging to $\mathcal{D}^+_\lambda$ and $v_\lambda$ to $\mathcal{D}^-_\lambda$. Moreover, $\|v_\lambda\|_{W_a^{1,p}(V)} > X(\lambda)> X_0> \|u_\lambda\|_{W_a^{1,p}(V)}$ and $\|v_\lambda\|_{\alpha+1} > S(\lambda) > S_0>\|u_\lambda\|_{\alpha+1} $, and $\|v_\lambda\|_W>\|v_\lambda\|_{W_a^{1,p}(V)} =X(\lambda)\rightarrow +\infty$ as $p\rightarrow \alpha+1$.
}

\vskip2mm
\noindent
{\bf Theorem 1.2.} {\it Let $G=(V, E)$ be a locally finite graph and it is connected, symmetric and there exists $\mu_{0}>0$ such that $\mu(x)\geq \mu_{0}$ for all $x\in V$. For any $x\in V$, $\sum_{y\thicksim x}\omega_{xy} <C$, where $C$ is a positive constant. Assume $\lambda < 0$,  and condition $(A)$ holds. Then  equation (\ref{eq1}) has one and only one positive solution.}
\vskip2mm
\noindent
\par
This paper is organized as follows. In section 2 we present Sobolev embedding results and a generalized Brezis-Lieb lemma on locally finite graphs. Besides, we clarify some technical results about the energy functional $J_\lambda$ and the Nehari manifold associated to $J_\lambda$. In section 3, we demonstrate the multiplicity of solutions and analyze the behavior of the solution $v_\lambda$ in $\mathcal{D}^-_\lambda$ as $p\rightarrow \alpha+1$
for $0< \lambda< \Lambda_*$ (Theorem 1.1). In section 4, we show the uniqueness of a solution for $\lambda< 0$ (Theorem 1.2).

\vskip2mm
{\section{Preliminaries}}
  \setcounter{equation}{0}
  \par
In this section, we present and prove a Sobolev embedding theorem and some lemmas which are important in the proofs of our main results.
 \vskip2mm
\noindent
{\bf Lemma 2.1.} {\it Let $G=(V,E)$ be a locally finite graph. Assume $\mu(x)\ge \mu_0>0$, $a,b$ satisfies $(A)$ and $f \in L^1(V)\cap L^{\frac{p}{p-1+\gamma}}(V)  \cap L^{\frac{q}{q-1+\gamma}}(V)$. Then there hold\\
$(i)$\;\;
\begin{eqnarray*}
\|u\|_\infty \leq  \frac{1}{a_0^{\frac{1}{p}}\mu_0^{\frac{1}{p}}}\|u\|_{W_a^{1,p}(V)},\;\;\;
\|u\|_\infty \leq \frac{1}{b_0^{\frac{1}{q}}\mu_0^{\frac{1}{q}}}\|u\|_{W_b^{1,q}(V)}
\end{eqnarray*}
and
\begin{eqnarray}
\label{x13}
\|u\|_\theta \leq  \mu_0^{\frac{p-\theta}{p\theta}}a_0^{-\frac{1}{p}}\|u\|_{W_a^{1,p}(V)},\;\;\|u\|_\theta \leq  \mu_0^{\frac{q-\theta}{q\theta}}b_0^{-\frac{1}{q}}\|u\|_{W_b^{1,q}(V)},\;\;\mbox{for all $1<q\leq p\leq \theta < \infty$}.
\end{eqnarray}
$(ii)$\;\;
\begin{eqnarray}
\label{x14}
\int_V f(x)|u|^{1-\gamma}d\mu \leq a_0^{-\frac{1-\gamma}{p}}\|f\|_{\frac{p}{p-1+\gamma}}\cdot\|u\|^{1-\gamma}_{W_a^{1,p}(V)},\;\;
\int_V f(x)|u|^{1-\gamma}d\mu  \leq b_0^{-\frac{1-\gamma}{q}} \|f\|_{\frac{q}{q-1+\gamma}}\|u\|^{1-\gamma}_{W_b^{1,q}(V)}.
\end{eqnarray}
}
\vskip2mm
\noindent
{\bf Proof}. The conclusion $(i)$ has been proved in \cite{Yang P 2023}. Next, we prove the conclusion (ii). In fact, by H\"{o}lder's inequality and (\ref{x13}), we infer that
\begin{eqnarray*}
\int_V f(x)|u|^{1-\gamma}d\mu &\leq & \|f\|_{\frac{p}{p-1+\gamma}} \|u\|_p^{1-\gamma}\\
&\leq & a_0^{-\frac{1-\gamma}{p}}\|f\|_{\frac{p}{p-1+\gamma}}\|u\|^{1-\gamma}_{W_a^{1,p}(V)}.
\end{eqnarray*}
Similarly, we also obtain $\int_V f(x)|u|^{1-\gamma}d\mu  \leq b_0^{-\frac{1-\gamma}{q}} \|f\|_{\frac{q}{q-1+\gamma}}\|u\|^{1-\gamma}_{W_b^{1,q}(V)}$.
\qed

 \vskip2mm
\noindent
{\bf Lemma 2.2.}  {\it (Brezis-Lieb Lemma on locally finite graphs) Assume that condition $(A)$ holds. Then for a bounded sequence $\{u_k\}\subset W$, there exists $u \in W$ such that
\begin{eqnarray*}
\lim_{k\rightarrow \infty} \left(\|u_k\|^{p}_{W_a^{1,p}(V)}-\|u_k-u\|^{p}_{W_a^{1,p}(V)} \right)=\|u\|^{p}_{W_a^{1,p}(V)}\;\;\mbox{and}\;\;\lim_{k\rightarrow \infty} \left(\|u_k\|^{q}_{W_b^{1,q}(V)}-\|u_k-u\|^{q}_{W_b^{1,q}(V)} \right)=\|u\|^{q}_{W_b^{1,q}(V)}.
\end{eqnarray*}
}
\vskip2mm
\noindent
{\bf Proof.}\ \ Since $\{u_k\}\subset W$ is bounded, there exists a $u \in W$ such that, up to a subsequence, $u_k \rightharpoonup u$ in $W$. In particular,
$$
\lim_{k\rightarrow \infty} \int_V u_k \varphi d\mu=\int_V u  \varphi d\mu,\;\;\forall\;\varphi\in C_c(V),
$$
which implies
\begin{eqnarray}
\label{xx1}
\lim_{k\rightarrow \infty} u_k(x)=u(x)\;\;\mbox{for any fixed $x\in V$},
\end{eqnarray}
if we choose $\varphi\in C_c(V)$ defined by
\begin{eqnarray*}
\varphi(y)=
\begin{cases}
1,\;\;y=x,\\
0,\;\;y\neq x.
\end{cases}
\end{eqnarray*}
For every sufficiently small $\varepsilon>0$,  there exists $C_\varepsilon>0$ such that
\begin{eqnarray*}\label{z2}
\left |a(x) |u_k|^p -a(x)|u_k-u|^{p} \right|   \leq \varepsilon a(x) |u_k-u|^{p} +  C_\varepsilon  a(x) |u|^{p}.
\end{eqnarray*}
Let $u_k =u+ (u_k-u)$ be a sequence of measurable functions such that:\\
$(i)$\;\;$u_k-u \rightarrow 0$ for all $x\in V$;\\
$(ii)$\;\;By the definition of $W_a^{1,p}(V)$, we can get $\int_V a(x)|u|^p d\mu < \infty;$\\
$(iii)$\;\;Since $\{u_k\}\subset W$ is bounded, we obtain $\int_V a(x)|u_k-u|^{p} d\mu \leq \|u_k-u\|_{W_a^{1,p}(V)} \leq C <\infty$. The constant $C$ is independent of $\varepsilon$ and $k$;\\
$(iv)$\;\;$\int_V C_\varepsilon a(x)|u|^{p} d\mu \leq C_\varepsilon \|u\|_{W_a^{1,p}(V)} <\infty$ for all $\varepsilon>0$.\\
Then it follows from  Theorem 2 in \cite{Brezis 1983} that
\begin{eqnarray}\label{qw2}
\lim_{k\rightarrow \infty} \left(\int_V( a(x)| u_k|^p - a(x) | (u_k-u)|^{p}) d\mu \right)=\int_V  a(x)| u|^p d\mu.
\end{eqnarray}
Moreover, by the definition of gradient and (\ref{xx1}), it is easy to see that  $\lim_{k\rightarrow \infty} \nabla u_k(x) = \nabla u(x)$ for all $x\in V$.  Similar to the argument of (\ref{qw2}) with $a(x)\equiv 1$ and replacing $u_k$ and $u$ with $\nabla u_k$ and $\nabla u$, respectively, it holds that
\begin{eqnarray}\label{q2}
\lim_{k\rightarrow \infty} \left(\int_V (|\nabla u_k|^p -|\nabla (u_k-u)|^{p}) d\mu \right)=\int_V |\nabla u|^p d\mu.
\end{eqnarray}
Combining with $(\ref{q2})$, we conclude
\begin{eqnarray*}
\lim_{k\rightarrow \infty} \left(\|u_k\|^{p}_{W_a^{1,p}(V)}-\|u_k-u\|^{p}_{W_a^{1,p}(V)} \right)=\|u\|^{p}_{W_a^{1,p}(V)}.
\end{eqnarray*}
Similarly, we also obtain that
\begin{eqnarray*}
\lim_{k\rightarrow \infty} \left(\|v_k\|^{q}_{W_b^{1,q}(V)}-\|v_k-v\|^{q}_{W_b^{1,q}(V)} \right)=\|v\|^{q}_{W_b^{1,q}(V)}.
\end{eqnarray*}
The proof is complete.
\qed
 \vskip2mm
\noindent
{\bf Lemma 2.3.}  {\it If $u\in W_+$ is a weak solution of (\ref{eq1}), $u$ is also a point-wise solution of equation (\ref{eq1}).
}
\vskip2mm
\noindent
{\bf Proof.}\ \ Since $u\in W_+$ is a weak solution of equation  (\ref{eq1}), for any $\phi \in W$, there holds
\begin{eqnarray}
\label{s3}
&&\int_V[|\nabla u|^{p-2} \Gamma(u,\phi)+a(x)u^{p-1}\phi]d\mu
+\int_V [|\nabla u|^{q-2} \Gamma(u,\phi)+b(x)u^{q-1}\phi]d\mu\nonumber\\
&=&\int_Vf(x)u^{-\gamma} \phi d\mu+{\lambda}\int_V g(x)u^{\alpha}\phi d\mu.
\end{eqnarray}
By Lemma 2.1 in \cite{Han X L 2021}, we obtain
\begin{eqnarray}
\label{s2}
\int_V|\nabla u|^{s-2} \Gamma(u,\phi)d\mu
=-\int_V (\Delta_s u)\phi d\mu, \ \ \mbox{ for }s=p,q.
\end{eqnarray}
Inserting (\ref{s2}) into (\ref{s3}), it follows that
\begin{eqnarray}
\label{s1}
&&\int_V [-(\Delta_p u)\phi+a(x)u^{p-1}\phi]d\mu
+\int_V [- (\Delta_q u)\phi+b(x)u^{q-1}\phi]d\mu\nonumber\\
&=&\int_Vf(x)u^{-\gamma}\phi d\mu+{\lambda}\int_V g(x)u^{\alpha}\phi d\mu,\;\;\;\;\mbox{$\forall \phi \in C_c(V)$}.
\end{eqnarray}
For any fixed $x_0\in V$, taking a test function $\phi: V\rightarrow \R$ in (\ref{s1}) with
\begin{eqnarray*}
\phi(x)=
\begin{cases}
\begin{array}{ll}
  1\;\;\;& \mbox{$x=x_0$},\\
  0\;\;\;& \mbox{$x\neq x_0$},
  \end{array}
   \end{cases}
\end{eqnarray*}
we have
\begin{eqnarray*}
-\Delta_p u(x_0)-\Delta_q u(x_0)+a(x_0)u(x_0)^{p-1}
+b(x)u(x_0)^{q-1}
=f(x)u(x_0)^{-\gamma} +{\lambda} g(x)u(x_0)^{\alpha}.
\end{eqnarray*}
Since $x_0$ is arbitrary, we conclude that $u$ is a point-wise solution of equation (\ref{eq1}).
\qed

 \vskip2mm
\noindent
{\bf Lemma 2.4.}  {\it Assume that condition $(A)$ holds. Then for all $\lambda\in (0,\Lambda_{\alpha,p,\gamma})$, the following conclusions hold:\\
$(i)$\;\;$\mathcal{D}^+_\lambda$ and $\mathcal{D}^-_\lambda$ are non-empty;\\
$(ii)$\;\;$\mathcal{D}^0_\lambda=\{0\}$.
}
\vskip2mm
\noindent
{\bf Proof.}\ \ $(i)$\;\;For any given $u\in W_+\backslash \{0\}$. Note that $\varphi_u(t):(0,\infty)\mapsto \R$ defined by
\begin{eqnarray}
\label{x4}
\varphi_u(t)=
t^{p-1+\gamma}\|u\|^p_{W^{1,p}_a (V)}+t^{q-1+\gamma}\|u\|^q_{W_b^{1,q}(V)}-\int_V f(x)u^{1-\gamma}d\mu-\lambda t^{\alpha+\gamma}\int_V g(x)u^{\alpha+1}d\mu.
\end{eqnarray}
and
\begin{eqnarray*}
\label{e4}
\varphi'_{u}(t)
&=&(p-1+\gamma)t^{p-2+\gamma}\|u\|^p_{W^{1,p}_a (V)}+(q-1+\gamma)t^{q-2+\gamma}\|u\|^q_{W_b^{1,q}(V)}-(\alpha+\gamma)\lambda t^{\alpha+\gamma-1}\int_V g(x)u^{\alpha+1}d\mu\\
&=&t^{\alpha+\gamma-1}\left[(p-1+\gamma)t^{p-1-\alpha}\|u\|^p_{W^{1,p}_a (V)}+(q-1+\gamma)t^{q-1-\alpha}\|u\|^q_{W_b^{1,q}(V)}-(\alpha+\gamma)\lambda\int_V g(x)u^{\alpha+1}d\mu\right]\\
&:=&t^{\alpha+\gamma-1}\left[M(t)-(\alpha+\gamma)\lambda\int_V g(x)u^{\alpha+1}d\mu\right],
\end{eqnarray*}
where $M(t)=(p-1+\gamma)t^{p-1-\alpha}\|u\|^p_{W^{1,p}_a (V)}+(q-1+\gamma)t^{q-1-\alpha}\|u\|^q_{W_b^{1,q}(V)}$,
and $M'(t)=(p-1-\alpha)(p+\gamma-1)t^{p-\alpha-2}\|u\|^p_{W_a^{1,p}(V)}+(q-1-\alpha)(q+\gamma-1)t^{q-\alpha-2}\|u\|^q_{W_b^{1,q}(V)}$. Since $0<\gamma <1<q\leq p<\alpha+1$, we can infer that $M'(t)<0$, that is, $M(t)$ is strictly decreasing and $M(t)>0$ for $t\in(0,\infty)$. Actually, for any given $u\in W$, $g\not\equiv 0$ and $\lambda>0$, there exists $(\alpha+\gamma)\lambda\int_V g(x)u^{\alpha+1}d\mu>0$. Then we can find a unique $\tilde{t}\in(0,\infty)$ such that $M(t)=(\alpha+\gamma)\lambda\int_V g(x)u^{\alpha+1}d\mu$. Thus there is a unique $\tilde{t}\in(0,\infty)$ such that $\varphi'_{u}(\tilde{t})=0$. Then
\begin{eqnarray}
\varphi'_{u}(\tilde{t})=0
&     \Longleftrightarrow     & (p-1+\gamma)\tilde{t}^{p-1-\alpha}\|u\|^p_{W^{1,p}_a (V)}+(q-1+\gamma)\tilde{t}^{q-1-\alpha}\|u\|^q_{W_b^{1,q}(V)}\nonumber\\
&&=(\alpha+\gamma)\lambda\int_V g(x)u^{\alpha+1}d\mu\nonumber\label{xx3}\\
&     \Rightarrow             & (p-1+\gamma)\tilde{t}^{p-1-\alpha}\|u\|^p_{W^{1,p}_a (V)}< (\alpha+\gamma)\lambda\int_V g(x)u^{\alpha+1}d\mu\nonumber\\
&     \Rightarrow             & {\tilde{t}}<\left(\frac{(\alpha+\gamma)\lambda\int_V g(x)u^{\alpha+1}d\mu}{(p-1+\gamma)\|u\|^p_{W^{1,p}_a (V)}}\right)^{\frac{1}{p-\alpha-1}}:= {t^*}\label{xx3}
\end{eqnarray}
Inserting (\ref{xx3}) into (\ref{x4}) and applying (\ref{x14}), we derive
\begin{eqnarray}
&&\varphi_{u}(t^*)\nonumber\\
&=&{(t^*)}^{p-1+\gamma}\|u\|^p_{W_a^{1,p}(V)}+{(t^*)}^{q-1+\gamma}\|u\|^q_{W_b^{1,q}(V)}
-\int_V f(x)u^{1-\gamma}d\mu-\lambda{(t^*)}^{\alpha+\gamma} \int_V g(x)u^{\alpha+1}d\mu\nonumber\\
&>&\left(\frac{(\alpha+\gamma)\lambda\int_V g(x)u^{\alpha+1}d\mu}{(p-1+\gamma)\|u\|^p_{W_a^{1,p}(V)}}\right)^{\frac{p+\gamma-1}{p-\alpha-1}}\|u\|^p_{W_a^{1,p}(V)}\nonumber\\
&&-\int_V f(x)u^{1-\gamma}d\mu-\left(\frac{(\alpha+\gamma)\lambda\int_V g(x)u^{\alpha+1}d\mu}{(p-1+\gamma)\|u\|^p_{W_a^{1,p}(V)}}\right)^{\frac{\alpha+\gamma}{p-\alpha-1}}\lambda\int_V g(x)u^{\alpha+1}d\mu\nonumber\\
&=&
\left(\frac{\alpha+\gamma}{p-1+\gamma}\right)^{\frac{\alpha+\gamma}{p-\alpha-1}}\cdot\left(\frac{\alpha+\gamma}{p-1+\gamma}-1\right)\cdot\left(\frac{\left(\lambda\int_V g(x)u^{\alpha+1}d\mu\right)^{p+\gamma-1}}{\|u\|^{p(\alpha+\gamma)}_{W_a^{1,p}(V)}}\right)^{\frac{1}{p-\alpha-1}}-\int_V f(x)u^{1-\gamma}d\mu\nonumber\\
&=&
\left(\frac{\alpha+1-p}{p-1+\gamma}\right)\cdot\left(\frac{\alpha+\gamma}{p-1+\gamma}\right)^{\frac{\alpha+\gamma}{p-\alpha-1}}\left(\frac{\|u\|^{p(\alpha+\gamma)}_{W_a^{1,p}(V)}}{\left(\lambda\int_V g(x)u^{\alpha+1}d\mu\right)^{p+\gamma-1}}\right)^{\frac{1}{\alpha+1-p}}-\int_V f(x)u^{1-\gamma}d\mu\label{x7}\\
&\geq&
\left(\frac{\alpha+1-p}{p-1+\gamma}\right)\cdot\left(\frac{\alpha+\gamma}{p-1+\gamma}\right)^{\frac{\alpha+\gamma}{p-\alpha-1}}\left(\frac{\|u\|^{p(\alpha+\gamma)}_{W_a^{1,p}(V)}}
{\left(\lambda \|g\|_\infty \mu_0^{\frac{p-\alpha-1}{p}}a_0^{-\frac{\alpha+1}{p}}\|u\|^{\alpha+1}_{W_a^{1,p}(V)}\right)^{p+\gamma-1}}\right)^{\frac{1}{\alpha+1-p}}\nonumber\\
&&-a_0^{-\frac{1-\gamma}{p}} \|f\|_{\frac{p}{p-1-\gamma}}\|u\|^{1-\gamma}_{W_a^{1,p}(V)}\nonumber\\
&=&
\left[\left(\frac{\alpha+1-p}{p-1+\gamma}\right)\cdot\left(\frac{\alpha+\gamma}{p-1+\gamma}\right)^{\frac{\alpha+\gamma}{p-\alpha-1}}
\left(\frac{1}{\lambda}\right)^{\frac{p-1+\gamma}{\alpha+1-p}}\left(\frac{1}
{\|g\|_\infty \mu_0^{\frac{p-\alpha-1}{p}}a_0^{-\frac{\alpha+1}{p}}}\right)^{\frac{p+\gamma-1}{\alpha+1-p}}
-a_0^{-\frac{1-\gamma}{p}} \|f\|_{\frac{p}{p-1-\gamma}}\right]\|u\|^{1-\gamma}_{W_a^{1,p}(V)}\nonumber\\
&\equiv&Z(\lambda)\|u\|^{1-\gamma}_{W_a^{1,p}(V)}\nonumber
\end{eqnarray}
where $Z(\lambda):=\left(\frac{\alpha+1-p}{p-1+\gamma}\right)\cdot\left(\frac{\alpha+\gamma}{p-1+\gamma}\right)^{\frac{\alpha+\gamma}{p-\alpha-1}}
\left(\frac{1}{\lambda}\right)^{\frac{p-1+\gamma}{\alpha+1-p}}\left(\frac{1}
{\|g\|_\infty \mu_0^{\frac{p-\alpha-1}{p}}a_0^{-\frac{\alpha+1}{p}}}\right)^{\frac{p+\gamma-1}{\alpha+1-p}}
-a_0^{-\frac{1-\gamma}{p}} \|f\|_{\frac{p}{p-1-\gamma}}$.
Since $\varphi_u(t)\rightarrow -\int_V f(x)u^{1-\gamma}d\mu$ as $t\rightarrow 0$ and $\varphi_u(t)\rightarrow -\infty$ as $t\rightarrow \infty$ with $\tilde{t}$ being the unique stationary point of $\varphi_u(t)$, it follows that $\tilde{t}\in(0,\infty)$ must be the unique maximum value point of $\varphi_u(t)$. From (\ref{xx3}), we have $0<\tilde{t}<t^*$, which implies $\varphi_u(\tilde{t})>\varphi_u(t^*)$. Furthermore, for $0<\lambda<\Lambda_{*}$, we deduce $\varphi_u(t^*)>0$, and consequently $\varphi_u(\tilde{t})>0$.
Since there is a unique $\tilde{t}$ such that $\varphi'_{u}(\tilde{t})=0$, we infer that $\varphi_u(t)$ is increasing and then decreasing along $t>0$. Therefore, for each $\lambda\in (0,\Lambda_{*})$, there exist exactly two points $0<t_1<\tilde{t}<t_2 $ such that $\varphi_u(t_i)=0$, for $i=1,2$. This shows that $t_i(u)u\in \mathcal{D}_\lambda$.
Besides, we also have $\varphi_{u}'({t_1})>0$ and $\varphi_{u}'({t_2})<0$(see Figure 1). From (\ref{x1}) and (\ref{x2}), we get that $t_1(u)u\in \mathcal{D}^+_\lambda$ and $t_2(u)u\in \mathcal{D}^-_\lambda$. In conclusion, $\mathcal{D}^+_\lambda$ and $\mathcal{D}^-_\lambda$ are non-empty for $\lambda\in (0,\Lambda_{*})$.
\par
$(ii)$\;\;We claim that $\mathcal{D}^0_\lambda=\{0\}$. Otherwise, we assume that there exists $u_0 \not\equiv 0$ and $u_0 \in \mathcal{D}^0_\lambda\subset\mathcal{D}_\lambda$. It follows that
\begin{eqnarray}
\|u_0\|^p_{W_a^{1,p}(V)}+\|u_0\|^q_{W_b^{1,q}(V)}-\int_V f(x)u_0^{1-\gamma}d\mu-\lambda\int_V g(x)u_0^{\alpha+1}d\mu=0\label{x5}
\end{eqnarray}
and
\begin{eqnarray}
(p-1+\gamma)\|u_0\|^p_{W_a^{1,p}(V)}+(q-1+\gamma)\|u_0\|^q_{W_b^{1,q}(V)}-(\alpha+\gamma)\lambda\int_V g(x)u_0^{\alpha+1}d\mu=0.\label{x6}
\end{eqnarray}
From (\ref{x7}), (\ref{x5}) and (\ref{x6}), we obtain that
\begin{eqnarray*}
&&0\\
&   <  &
Z(\lambda)\|u_0\|^{1-\gamma}_{W_a^{1,p}(V)}\\
&\leq &
\left(\frac{\alpha+1-p}{p-1+\gamma}\right)\cdot\left(\frac{\alpha+\gamma}{p-1+\gamma}\right)^{\frac{\alpha+\gamma}{p-\alpha-1}}\cdot
\left(\frac{\|u_0\|^{p(\alpha+\gamma)}_{W_a^{1,p}(V)}}{\left(\lambda\int_V g(x)u_0^{\alpha+1}d\mu\right)^{p+\gamma-1}}\right)^{\frac{1}{\alpha+1-p}}
-\int_V f(x)u_0^{1-\gamma}d\mu\\
&  =     &
\left(\frac{\alpha+1-p}{p-1+\gamma}\right)\cdot\left(\frac{\alpha+\gamma}{p-1+\gamma}\right)^{\frac{\alpha+\gamma}{p-\alpha-1}}\cdot\\
&&\left(\frac{(\alpha+\gamma)^{p-1+\gamma}\|u_0\|^{p(\alpha+\gamma)}_{W_a^{1,p}(V)}}{\left((p-1+\gamma)\|u_0\|^{p}_{W_a^{1,p}(V)}+(q-1+\gamma)\|u_0\|^q_{W_b^{1,q}(V)}\right)^{p-1+\gamma}}\right)^{\frac{1}{\alpha+1-p}}
-\|u_0\|^{p}_{W_a^{1,p}(V)}\\
&&-\|u_0\|^q_{W_b^{1,q}(V)}+\frac{(p-1+\gamma)\|u_0\|^{p}_{W_a^{1,p}(V)}+(q-1+\gamma)\|u_0\|^q_{W_b^{1,q}(V)}}{\alpha+\gamma}\\
&  <     &
\left(\frac{\alpha+1-p}{p-1+\gamma}\right)\cdot\left(\frac{\alpha+\gamma}{p-1+\gamma}\right)^{\frac{\alpha+\gamma}{p-\alpha-1}}\cdot
\left(\frac{(\alpha+\gamma)^{p-1+\gamma}\|u_0\|^{p(\alpha+\gamma)}_{W_a^{1,p}(V)}}{\left((p-1+\gamma)\|u_0\|^{p}_{W_a^{1,p}(V)}\right)^{p-1+\gamma}}\right)^{\frac{1}{\alpha+1-p}}\\
&&-\left(1-\frac{p-1+\gamma}{\alpha+\gamma}\right)\|u_0\|^{p}_{W_a^{1,p}(V)}-\left(1-\frac{q-1+\gamma}{\alpha+\gamma}\right)\|u_0\|^q_{W_b^{1,q}(V)}\\
&   <     &
\left(\frac{\alpha+1-p}{p-1+\gamma}\right)\cdot\left(\frac{\alpha+\gamma}{p-1+\gamma}\right)^{\frac{\alpha+\gamma}{p-\alpha-1}}\cdot
\left(\frac{\alpha+\gamma}{p-1+\gamma}\right)^{\frac{p-1+\gamma}{\alpha+1-p}}\|u_0\|^{p}_{W_a^{1,p}(V)}
-\left(\frac{\alpha-p+1}{\alpha+\gamma}\right)\|u_0\|^{p}_{W_a^{1,p}(V)}\\
&   =     &
\left[ \left(\frac{\alpha+1-p}{p-1+\gamma}\right)\cdot\left(\frac{p-1+\gamma}{\alpha+\gamma}\right)^{\frac{p-1+\gamma}{\alpha+1-p}}\cdot\frac{p-1+\gamma}{\alpha+\gamma}\cdot
\left(\frac{\alpha+\gamma}{p-1+\gamma}\right)^{\frac{p-1+\gamma}{\alpha+1-p}}-\frac{\alpha-p+1}{\alpha+\gamma}
\right]\|u_0\|^{p}_{W_a^{1,p}(V)}\\
&=&0.
\end{eqnarray*}
It's a contradiction. Thus $u_0 \equiv 0$.
The proof is complete.
\qed

\vskip2mm
\noindent
{\bf Remark 2.1.}
To establish the positivity of $\varphi_u(\tilde{t})$, the methodology employed in \cite{Sun 2008} becomes inapplicable for our case. Their approach relied on explicit computation of $\tilde{t}$ through direct substitution into $\varphi_u(t)$, this strategy proves inadequate for our generalized equation (\ref{eq1}) with $p\geq q>1$ (as opposed to the special case $p= q=2$ in \cite{Sun 2008}), which creates substantial computational barriers for determining $\tilde{t}$ explicitly. To resolve this challenge, we initiate our analysis by investigating the monotonicity properties of $\varphi_u(t)$. Through rigorous examination, we establish that $\varphi_u(t)$ possesses a unique stationary point at $\tilde{t} \in(0,\infty)$, which corresponds to its global maximum. Building upon this characterization, we subsequently identify a value $t^* \in(0,\infty)$ satisfying the crucial inequality:
$0<\varphi_u(t^*)<\varphi_u(\tilde{t})$.
\vskip2mm
\noindent
{\bf Remark 2.2.} For any given $u\in W_+\setminus\{0\}$, since $\varphi_u(t)=t^{\gamma}\frac{d}{dt}J_\lambda(tu)$, we know $\frac{d}{dt}J_\lambda(tu) = t^{-\gamma}\varphi_u(t)$ for all $t>0$. By the proof of Lemma 2.4 $(i)$, we get that $\varphi'_{u}(t_1)>0$, $\varphi'_{u}(t_2)<0$, $\varphi_u(t_1)=\varphi_u(t_2)=0$, $\tilde{t}$ is the only stationary point and $\varphi_u(\tilde{t})>0$ (see Figure 1). Then we can deduce that $J_\lambda(tu)$ is increasing over the intervals $[t_1, t_2]$ and decreasing over the intervals $(0,t_1]$ and $[t_2, \infty]$ (see Figure 2 and Figure 3).
\captionsetup{justification=centering, singlelinecheck=false, position=bottom}
\begin{figure}[H]
  \centering
  \begin{minipage}[c]{0.3\textwidth}
  \includegraphics[width=0.85\textwidth]{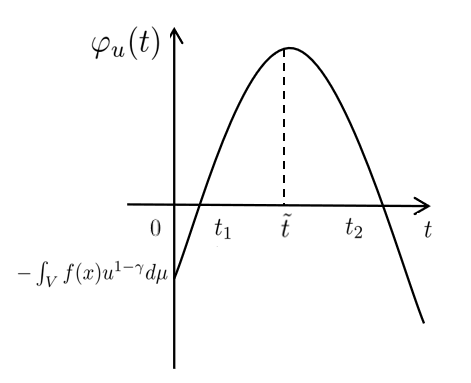}\\
    \vspace{-0.27cm}
  \caption{the  behaviour of $\varphi_u(t)$.}\label{a1}
   \end{minipage}
   \begin{minipage}[c]{0.3\textwidth}
  \includegraphics[width=0.85\textwidth]{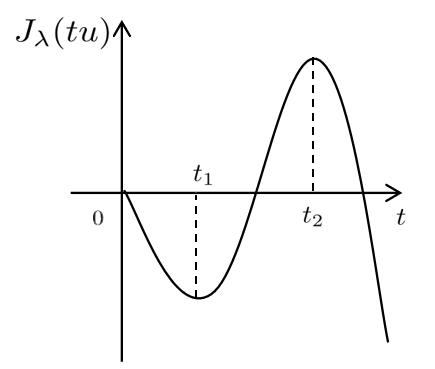}\\
    \vspace{-0.4cm}
  \caption{the first kind of behaviour of $J_\lambda(tu)$.}\label{a2}
   \end{minipage}
    \begin{minipage}[c]{0.3\textwidth}
  \includegraphics[width=0.85\textwidth]{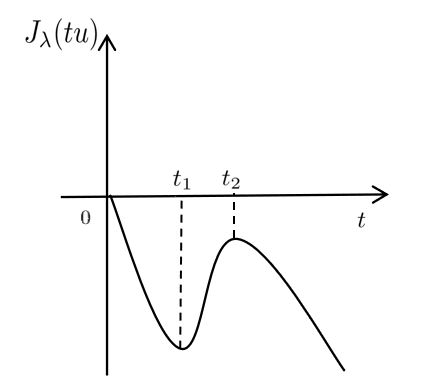}\\
    \vspace{-0.62cm}
  \caption{the second kind of behaviour of $J_\lambda(tu)$.}\label{a3}
   \end{minipage}
\end{figure}

 \vskip2mm
\noindent
{\bf Lemma 2.5.}  {\it
Assume that condition $(A)$ holds. Then for each $u\in \mathcal{D}^-_\lambda$ and $v\in \mathcal{D}^+_\lambda$ for $\lambda\in(0, \Lambda_*)$, the following conclusion hold:\\
$(i)$\;\; $\|u\|_{W_a^{1,p}(V)} > X(\lambda)> X_0> \|v\|_{W_a^{1,p}(V)}$;\\
$(ii)$\;\;$\|u\|_{\alpha+1} > S(\lambda) > S_0>\|v\|_{\alpha+1} $.
}
\vskip2mm
\noindent
{\bf Proof.}\ \ $(i)$\;\;Since $u\in \mathcal{D}^-_\lambda \subset \mathcal{D}_\lambda$, applying  (\ref{x2}) and (\ref{x13}) yields
\begin{eqnarray}
\label{x11}
0&  >   &(p-1+\gamma)\|u\|^p_{W_a^{1,p}(V)}+(q-1+\gamma)\|u\|^q_{W_b^{1,q}(V)}-(\alpha+\gamma)\lambda\int_V g(x)u^{\alpha+1}d\mu\nonumber\\
&   >   &(p-1+\gamma)\|u\|^p_{W_a^{1,p}(V)}-\lambda(\alpha+\gamma)\int_V g(x)u^{\alpha+1}d\mu\nonumber\\
&  \geq &(p-1+\gamma)\|u\|^p_{W_a^{1,p}(V)}-\lambda(\alpha+\gamma)\|g\|_\infty\mu_0^{\frac{p-1-\alpha}{p}}a_0^{-\frac{\alpha+1}{p}}\|u\|^{\alpha+1}_{W_a^{1,p}(V)}\nonumber\\
 \Rightarrow  &&(p-1+\gamma)\|u\|^p_{W_a^{1,p}(V)} < \lambda(\alpha+\gamma)\|g\|_\infty\mu_0^{\frac{p-1-\alpha}{p}}a_0^{-\frac{\alpha+1}{p}}\|u\|^{\alpha+1}_{W_a^{1,p}(V)}\nonumber\\
\Rightarrow  && \|u\|_{W_a^{1,p}(V)} > \left( \frac{p-1+\gamma}{\lambda(\alpha+\gamma)\|g\|_\infty\mu_0^{\frac{p-1-\alpha}{p}}a_0^{-\frac{\alpha+1}{p}}}\right)
^{\frac{1}{\alpha+1-p}}=X(\lambda).
\end{eqnarray}
For any $v\in \mathcal{D}^+_\lambda \subset \mathcal{D}_\lambda$, it follows that
\begin{eqnarray}
\label{x8}
(p-1+\gamma)\|v\|^p_{W_a^{1,p}(V)}+(q-1+\gamma)\|v\|^q_{W_b^{1,q}(V)}-(\alpha+\gamma)\lambda\int_V g(x)v^{\alpha+1}d\mu>0
\end{eqnarray}
and
\begin{eqnarray}
\label{x9}
\|v\|^p_{W_a^{1,p}(V)}+\|v\|^q_{W_b^{1,q}(V)}-\int_V f(x)v^{1-\gamma}d\mu-\lambda\int_V g(x)v^{\alpha+1}d\mu=0.
\end{eqnarray}
Combining (\ref{x9}) with (\ref{x8}) through (\ref{x14}), we derive
\begin{eqnarray}
\label{x12}
0&  <   &(p-1+\gamma)\|v\|^p_{W_a^{1,p}(V)}+(q-1+\gamma)\|v\|^q_{W_b^{1,q}(V)}-(\alpha+\gamma)\lambda\int_V g(x)v^{\alpha+1}d\mu\nonumber\\
&   =   & (p-1+\gamma)\|v\|^p_{W_a^{1,p}(V)}+(q-1+\gamma)\|v\|^q_{W_b^{1,q}(V)}\nonumber\\
&&-(\alpha+\gamma)\|v\|^p_{W_a^{1,p}(V)}-(\alpha+\gamma)\|v\|^q_{W_b^{1,q}(V)}+(\alpha+\gamma)\int_V f(x)v^{1-\gamma}d\mu\nonumber\\
&  =   &(p-1-\alpha)\|v\|^p_{W_a^{1,p}(V)}+(q-1-\alpha)\|v\|^q_{W_b^{1,q}(V)}+(\alpha+\gamma)\int_V f(x)v^{1-\gamma}d\mu\nonumber\\
&   <   &(p-1-\alpha)\|v\|^p_{W_a^{1,p}(V)}+(\alpha+\gamma)\int_V f(x)v^{1-\gamma}d\mu\nonumber\\
&   \leq  &(p-1-\alpha)\|v\|^p_{W_a^{1,p}(V)}+(\alpha+\gamma) a_0^{-\frac{1-\gamma}{p}}\|f\|_{\frac{p}{p-1+\gamma}}\|v\|^{1-\gamma}_{W_a^{1,p}(V)}\nonumber\\
\Rightarrow  &&  \|v\|_{W_a^{1,p}(V)} <  \left( \frac{\alpha+\gamma}{\alpha+1-p}\right)^{\frac{1}{p-1+\gamma}}a_0^{-\frac{1-\gamma}{p(p-1+\gamma)}}\|f\|^{\frac{1}{p-1+\gamma}}_{\frac{p}{p-1+\gamma}} :=X_0.
\end{eqnarray}
By Appendix A.1, we have $X(\lambda)=X_0$ as $\lambda=\Lambda_{*}$. Then there exists $w \in W_+\setminus \{0\}$ such that
\begin{eqnarray*}
\begin{cases}
\|w\|_{W_a^{1,p}(V)}> X(\Lambda_{*})\;\;\;\;\;\;\mbox{if}\;\;w\in \mathcal{D}^-_\lambda;\\
\|w\|_{W_a^{1,p}(V)}<X(\Lambda_{*})\;\;\;\;\;\;\mbox{if}\;\;w\in \mathcal{D}^+_\lambda.
\end{cases}
\end{eqnarray*}
We refer to $\lambda=\Lambda_{*}$ as the dividing line between $\mathcal{D}^+_\lambda$ and $\mathcal{D}^-_\lambda$.
We also find that $X(\lambda) \rightarrow +\infty$ as $\lambda \rightarrow 0$. In this sense, the range of the parameter $\lambda$ is optimal. This tells us that the smaller the perturbation term, the larger the norm of the solution.
Besides, we derive that $\|u\|_{W_a^{1,p}(V)} - \|v\|_{W_a^{1,p}(V)} > X(\lambda)- X_0 > 0$ by the monotonicity of $X(\lambda)$.\\
$(ii)$\;\;For any $u\in \mathcal{D}^-_\lambda$, by the process of calculating (\ref{x11}) and (\ref{x13}), it follows that
\begin{eqnarray}
\label{x15}
0&   >   &(p-1+\gamma)\|u\|^p_{W_a^{1,p}(V)}-\lambda(\alpha+\gamma)\int_V g(x)u^{\alpha+1}d\mu\nonumber\\
&    \geq   &(p-1+\gamma)\left(\frac{1}{\mu_0}\right)^{\frac{p-\alpha-1}{\alpha+1}}a_0 \|u\|^p_{\alpha+1}-\lambda(\alpha+\gamma)\|g\|_\infty \|u\|^{\alpha+1}_{\alpha+1}\nonumber\\
&    =   &\left[(p-1+\gamma)\left(\frac{1}{\mu_0}\right)^{\frac{p-\alpha-1}{\alpha+1}}a_0 - \lambda(\alpha+\gamma)\|g\|_\infty \|u\|^{\alpha+1-p}_{\alpha+1}\right]\|u\|^{p}_{\alpha+1}\nonumber\\
\Rightarrow  &&\|u\|_{\alpha+1} > \left( \frac{p-1+\gamma}{\lambda(\alpha+\gamma)} \right)^{\frac{1}{\alpha+1-p}}
\left(\frac{1}{\mu_0}\right)^{-\frac{1}{\alpha+1}} a_0^{\frac{1}{\alpha+1-p}} \left( \frac{1}{\|g\|_\infty} \right)^{\frac{1}{\alpha+1-p}} =S(\lambda).
\end{eqnarray}
For any $v\in \mathcal{D}^+_\lambda$, by (\ref{x13}) and (\ref{x12}), we infer that
\begin{eqnarray*}
\label{x16}
& &\left( \frac{\alpha+\gamma}{\alpha+1-p}\right)^{\frac{1}{p-1+\gamma}}a_0^{-\frac{1-\gamma}{p(p-1+\gamma)}}
\|f\|^{\frac{1}{p-1+\gamma}}_{\frac{p}{p-1+\gamma}} >\|v\|_{W_a^{1,p}(V)}
 \geq   \left(\frac{1}{\mu_0}\right)^{\frac{p-\alpha-1}{p(\alpha+1)}} a_0^{\frac{1}{p}} \|v\|_{\alpha+1}\nonumber\\
&   \Rightarrow   &  \|v\|_{\alpha+1} < \left( \frac{\alpha+\gamma}{\alpha+1-p}\right)^{\frac{1}{p-1+\gamma}} \left(\frac{1}{\mu_0}\right)^{\frac{\alpha+1-p}{p(\alpha+1)}} a_0^{-\frac{1}{p-1+\gamma}}
\|f\|^{\frac{1}{p-1+\gamma}}_{\frac{p}{p-1+\gamma}}= S_0.
\end{eqnarray*}
By Appendix A.2, we have $S(\lambda)=S_0$ as $\lambda=\Lambda_{*}$ and $S(\lambda) \rightarrow +\infty$ as $\lambda \rightarrow 0$. Besides, we derive that $\|u\|_{\alpha+1} - \|v\|_{\alpha+1} > S(\lambda)- S_0 > 0$ by the monotonicity of $S(\lambda)$. The proof is complete.
\qed

 \vskip2mm
\noindent
{\bf Lemma 2.6.}  {\it Assume that condition $(A)$ holds. Then $\mathcal{D}^-_\lambda$ is closed in $W_+$ and $\mathcal{D}^+_\lambda$ is also closed in $W_+$ for $\lambda \in (0,\Lambda_{*})$.
}
\vskip2mm
\noindent
{\bf Proof.}\ \ Let $\{u_k\} \subset\mathcal{D}^-_\lambda$ be a sequence such that $u_k \rightarrow u_0$ in $W$ as $k\rightarrow \infty$. Then $u_k \rightarrow u_0$ in ${W_a^{1,p}(V)}$ and $u_k \rightarrow u_0$ in ${W_ b^{1,q}(V)}$ as $k\rightarrow \infty$. From Lemma 2.2 it follows that
\begin{eqnarray}\label{ax3}
\|u_k\|_{W_a^{1,p}(V)}\rightarrow \|u_0\|_{W_a^{1,p}(V)}\;\;\;\;\mbox{and}\;\;\;\; \|u_k\|_{W_b^{1,q}(V)}\rightarrow \|u_0\|_{W_b^{1,q}(V)}.
\end{eqnarray}
Moreover Lemma 2.1 implies that
\begin{eqnarray}\label{x18}
u_k(x)\rightarrow u_0(x),& \forall \;x\in V,\;k \rightarrow \infty.
\end{eqnarray}
It is easy to verify that $W_+$ is closed and convex. Hence $u_0 \in W_+$ and then $u_0 \geq 0$. Since $\{u_k\}$ is bounded in $W_+ \subset W$, there exists a constant $C_1$ such that
\begin{eqnarray}
\label{x26}
C_1 &   \geq    & \left(\int_V (|\nabla u_k|^p+a(x)|u_k|^p)d\mu\right)^{\frac{1}{p}}+\left(\int_V (|\nabla u_k|^q+b(x)|u_k|^q)d\mu\right)^{\frac{1}{q}}\nonumber\\
&   \geq    &
\mu_0^{\frac{\alpha+1-p}{p(\alpha+1)}} a_0^{\frac{1}{p}}\left( \int_V |u_k|^{\alpha+1}d\mu \right)^{\frac{1}{\alpha+1}}
+\mu_0^{\frac{\alpha+1-q}{q(\alpha+1)}} b_0^{\frac{1}{q}}\left( \int_V |u_k|^{\alpha+1} d\mu\right)^{\frac{1}{\alpha+1}}\;\;\;\mbox{(by (\ref{x13}))}\nonumber\\
&   \geq    &
\left(\mu_0^{\frac{\alpha+1-p}{p(\alpha+1)}} a_0^{\frac{1}{p}} +\mu_0^{\frac{\alpha+1-q}{q(\alpha+1)}} b_0^{\frac{1}{q}} \right) \mu_0^{\frac{1}{\alpha+1}} \|u_k\|_\infty\nonumber\\
&  =    &
\left(\mu_0^{\frac{1}{p}} a_0^{\frac{1}{p}} +\mu_0^{\frac{1}{q}} b_0^{\frac{1}{q}} \right) \|u_k\|_\infty\nonumber\\
\Rightarrow    &&  \|u_k\|_\infty \leq \frac{C_1}{ \mu_0^{\frac{1}{p}} a_0^{\frac{1}{p}} +\mu_0^{\frac{1}{q}} b_0^{\frac{1}{q}}}.
\end{eqnarray}
Then we obtain
\begin{eqnarray*}
f(x)(u_k^{1-\gamma}-u_0^{1-\gamma}) \leq f(x)u_k^{1-\gamma}\leq f(x)\left( \frac{C_1}{ \mu_0^{\frac{1}{p}} a_0^{\frac{1}{p}} +\mu_0^{\frac{1}{q}} b_0^{\frac{1}{q}}}\right) ^{1-\gamma}.
\end{eqnarray*}
Thus, let $G_1(x):= f(x) \left( \frac{C_1}{ \mu_0^{\frac{1}{p}} a_0^{\frac{1}{p}} +\mu_0^{\frac{1}{q}} b_0^{\frac{1}{q}}}\right) ^{1-\gamma}$, we have that $f(x)(u_k^{1-\gamma}-u_0^{1-\gamma})\leq G_1(x)$. Since $G_1(x) \in L^1(V)$, the Dominated Convergence Theorem with (\ref{x18}) yields
\begin{eqnarray}
\label{x21}
\lim_{k\rightarrow \infty}\int_V f(x)u_k^{1-\gamma}d\mu =\int_V f(x)u_0^{1-\gamma}d\mu.
\end{eqnarray}
Similarly, let $G_2(x):= g(x) \left( \frac{C_1}{ \mu_0^{\frac{1}{p}} a_0^{\frac{1}{p}} +\mu_0^{\frac{1}{q}} b_0^{\frac{1}{q}}}\right) ^{\alpha+1}$. Then we have that $g(x) |u_k|^{\alpha+1} \leq G_2(x)$. Since $G_2(x) \in L^1(V)$, the Dominated Convergence Theorem with (\ref{x18}) yields
\begin{eqnarray}
\label{x22}
\lim_{k\rightarrow \infty}\int_V g(x)u_k^{\alpha+1}d\mu = \int_V g(x)u_0^{\alpha+1}d\mu.
\end{eqnarray}
Since $u_k \in \mathcal{D}^-_\lambda \subset \mathcal{D}_\lambda$, from (\ref{ax3}), (\ref{x21}) and (\ref{x22}), we obtain
\begin{eqnarray*}
\label{x17}
 0
&=&   \|u_k\|^p_{W_a^{1,p}(V)}+\|u_k\|^q_{W_b^{1,q}(V)}-\int_V f(x)u_k^{1-\gamma}d\mu-\lambda\int_V g(x)u_k^{\alpha+1}d\mu\nonumber\\
&    =   &\lim_{k\rightarrow \infty} \left[\|u_k\|^p_{W_a^{1,p}(V)}+\|u_k\|^q_{W_b^{1,q}(V)}-\int_V f(x)u_k^{1-\gamma}d\mu-\lambda\int_V g(x)u_k^{\alpha+1}d\mu\right]\nonumber\\
&    =   &\|u_0\|^p_{W_a^{1,p}(V)}+\|u_0\|^q_{W_b^{1,q}(V)}-\int_V f(x)u_0^{1-\gamma}d\mu-\lambda\int_V g(x)u_0^{\alpha+1}d\mu.
\end{eqnarray*}
Then $u_0 \in \mathcal{D}_\lambda$. From (\ref{x22}) and the fact that $u_k\in\mathcal{D}^-_\lambda$, we derive that
\begin{eqnarray*}
&& 0\\
&    \geq    &
\lim_{k \rightarrow \infty}\left[(p-1+\gamma)\|u_k\|^p_{W_a^{1,p}(V)}+(q-1+\gamma)\|u_k\|^q_{W_b^{1,q}(V)}-(\alpha+\gamma)\lambda\int_V g(x)u_k^{\alpha+1}d\mu\right]\\
&    =    & (p-1+\gamma)\|u_0\|^p_{W_a^{1,p}(V)}+(q-1+\gamma)\|u_0\|^q_{W_b^{1,q}(V)}-(\alpha+\gamma)\lambda\int_V g(x)u_0^{\alpha+1}d\mu.
\end{eqnarray*}
Thus $u_0 \in \mathcal{D}^-_\lambda \cup \mathcal{D}^0_\lambda $. Besides, from Lemma 2.5, we know that $\|u_k\|_{\alpha+1}>S(\lambda)$. Then from Lemma 2.2, we get $\lim_{k \rightarrow \infty} \|u_k\|_{\alpha+1}= \|u_0\|_{\alpha+1} \geq S(\lambda)>0$. Therefore $u_0 \neq 0$. By Lemma 2.4 $(ii)$, it follows that $u_0 \in \mathcal{D}^-_\lambda$. Similarly, we can prove that $\mathcal{D}^+_\lambda$ is also closed in $W_+$.
\qed

 \vskip2mm
\noindent
{\bf Lemma 2.7.}  {\it Assume that condition $(A)$ holds. Then given $u\in \mathcal{D}^+_\lambda$ (resp. $\mathcal{D}^-_\lambda$), there exists a positive constant $\varepsilon$ and a continuous function $l=l(w)>0, w\in W_+, \|w\|_{W} \leq \varepsilon$ satisfying that
\begin{eqnarray*}
l(0)=1,\;\;l(w)(u+w) \in \mathcal{D}^+_\lambda\;\;(\mbox{resp.}\;\;\mathcal{D}^-_\lambda),\;\;\forall w\in W_+.
\end{eqnarray*}
}
\noindent
{\bf Proof.}\ \
We define that $L(t,w) = \varphi_{u+w}(t): (0, \infty) \times W_+ \mapsto \R $ as follows:
\begin{eqnarray*}
L(t,w) = t^{p-1+\gamma} \|u+w\|^p_{W_a^{1,p}(V)} + t^{q-1+\gamma}\|u+w\|^q_{W_b^{1,q}(V)}
- \int_V f(x)(u+w)^{1-\gamma}d\mu -\lambda t^{\alpha+\gamma}\int_V g(x)(u+w)^{\alpha+1}d\mu
\end{eqnarray*}
and
\begin{eqnarray*}
L'_t(t,w)&=& (p-1+\gamma) t^{p-2+\gamma} \|u+w\|^p_{W_a^{1,p}(V)} + (q-1+\gamma) t^{q-2+\gamma} \|u+w\|^q_{W_b^{1,q}(V)}\\
&&- (\alpha+\gamma)\lambda t^{\alpha+\gamma-1} \int_V g(x)(u+w)^{\alpha+1}d\mu .
\end{eqnarray*}
Since $u\in \mathcal{D}^+_\lambda \subset \mathcal{D}_\lambda$ (resp. $u\in \mathcal{D}^-_\lambda \subset \mathcal{D}_\lambda$), it holds that
\begin{eqnarray*}
L(1,0)  =     \|u\|^p_{W_a^{1,p}(V)} +  \|u\|^q_{W_b^{1,q}(V)}+  \int_V f(x)u^{1-\gamma}d\mu -\lambda \int_V g(x)u^{\alpha+1}d\mu
  =     0
\end{eqnarray*}
and
\begin{eqnarray*}
&&L'_t(1,0)  =      (p-1+\gamma) \|u\|^p_{W_a^{1,p}(V)} + (q-1+\gamma)  \|u\|^q_{W_b^{1,q}(V)}
-(\alpha+\gamma)\lambda\int_V g(x)u^{\alpha+1}d\mu
 >   0.\\
&&\mbox{(resp. $L'_t(1,0)<0$)}.
\end{eqnarray*}
Therefore we can apply the implicit function theorem at the point $(1,0)$. It follows that there exists a $\varepsilon_1>0$ and a continuous differential function $l=l(w):B(0, \varepsilon_1)\subset W_+ \rightarrow \R^+ $ satisfying that $l(0)=1$ and $L(l(w),w)=0$ for all $w \in W_+$ with $\|w\|_{W}< \varepsilon_1$. Thus we obtain that
\begin{eqnarray*}
0 &    =     & L(l(w),w)\\
&     =     &  (l(w))^{p-1+\gamma} \|u+w\|^p_{W_a^{1,p}(V)} + (l(w))^{q-1+\gamma} \|u+w\|^q_{W_b^{1,q}(V)}\\
&& - \int_V f(x)(u+w)^{1-\gamma}d\mu -(l(w))^{\alpha+\gamma}\lambda \int_V g(x)(u+w)^{\alpha+1}d\mu\\
&     =     &  (l(w))^{\gamma-1}\left[\|l(w)(u+w)\|^p_{W_a^{1,p}(V)} + \|l(w)(u+w)\|^q_{W_b^{1,q}(V)}\right.\\
&& \left. -  \int_V f(x)(l(w) (u+w))^{1-\gamma}d\mu -\lambda \int_V g(x)(l(w)(u+w))^{\alpha+1}d\mu\right]\\
\Rightarrow    && l(w) (u+w)\in \mathcal{D}_\lambda.
\end{eqnarray*}
From the local sign-preservation of $L'_t(t,w)$ and $L'_t(1,0)>0$ (resp. $L'_t(1,0)<0$), there exists a $\varepsilon_2>0$ such that $L'_t(l(w),w)>0$ (resp. $L'_t(l(w),w)<0$ ) as $(l(w),w) \in [1-\varepsilon_2,1+\varepsilon_2] \times [-\varepsilon_1,\varepsilon_2].$
Hence, taking $\varepsilon=\min\{\varepsilon_1,\varepsilon_2\}>0$, we obtain that
\begin{eqnarray*}
&&l(w)(u+w) \in \mathcal{D}^+_\lambda\;\;\mbox{for $\forall w\in W_+$ with $\|w\|_{W}<\varepsilon$}.\\
&&\mbox{(resp. $l(w)(u+w) \in \mathcal{D}^-_\lambda$  for $\forall w\in W_+$ with $\|w\|_{W}<\varepsilon$.)}
\end{eqnarray*}
We complete the proof.
\qed

 \vskip2mm
\noindent
{\bf Lemma 2.8.}  {\it Assume that condition $(A)$ holds. The following properties hold:\\
$(i)$\;\;$J_\lambda(u)$ is coercive on $\mathcal{D}_\lambda$;\;\;
$(ii)$\;\;$\inf_{\mathcal{D}^+_\lambda \cup \mathcal{D}^0_\lambda} J_\lambda(u) \in (-\infty, 0)$.
}
\vskip2mm
\noindent
{\bf Proof.}\ \ $(i)$\;\;Since $u \in \mathcal{D}_\lambda$, then
\begin{eqnarray}
\label{x19}
\lambda\int_V g(x)u^{\alpha+1}d\mu=   \|u\|^p_{W_a^{1,p}(V)}+ \|u\|^q_{W_b^{1,q}(V)}-\int_V f(x)u^{1-\gamma}d\mu .
\end{eqnarray}
We first introduce an inequality (see \cite{Xie J 2018}):
\begin{eqnarray}
\label{s4}
|x|^p+|y|^q \geq \frac{(|x|+|y|)^{\min\{p,q\}}}{\max\{2^{q-1},2^{p-1}\}}-1,
\end{eqnarray}
for all $p,q>1$ and $(x,y) \in \R^N\times\R^N$.
Combining (\ref{x19}), (\ref{x14}) with (\ref{s4}), we derive that
\begin{eqnarray*}
&&J_\lambda(u)\\
&    =    &
\frac{1}{p}  \|u\|^p_{W_a^{1,p}(V)} +\frac{1}{q} \|u\|^q_{W_b^{1,q}(V)}-\frac{1}{1-\gamma}\int_V f(x)u^{1-\gamma}d\mu -\frac{1}{\alpha+1}\left(\|u\|^p_{W_a^{1,p}(V)}+ \|u\|^q_{W_b^{1,q}(V)}-\int_V f(x)u^{1-\gamma}d\mu \right)\\
&      =     &
\left(\frac{1}{p} - \frac{1}{\alpha+1} \right)\|u\|^p_{W_a^{1,p}(V)}+\left(\frac{1}{q} - \frac{1}{\alpha+1} \right) \|u\|^q_{W_b^{1,q}(V)}-\left(\frac{1}{1-\gamma} - \frac{1}{\alpha+1} \right)\int_V f(x)u^{1-\gamma}d\mu\\
&      \geq     &
\left(\frac{1}{p} - \frac{1}{\alpha+1} \right)\|u\|^p_{W_a^{1,p}(V)}+\left(\frac{1}{q} - \frac{1}{\alpha+1} \right) \|u\|^q_{W_b^{1,q}(V)}-\left(\frac{1}{1-\gamma} - \frac{1}{\alpha+1} \right)a_0^{-\frac{1-\gamma}{p}}\|f\|_{\frac{p}{p-1+\gamma}}\|u\|^{1-\gamma}_{W_a^{1,p}(V)}\\
&      \geq     &
\left(\frac{1}{p} - \frac{1}{\alpha+1} \right) \left(\|u\|^p_{W_a^{1,p}(V)}+\|u\|^q_{W_b^{1,q}(V)}\right)- \left(\frac{1}{1-\gamma} - \frac{1}{\alpha+1} \right)a_0^{-\frac{1-\gamma}{p}}\|f\|_{\frac{p}{p-1+\gamma}}\|u\|^{1-\gamma}_{W_a^{1,p}(V)}\\
&      \geq     &
\frac{1}{2^{p-1}}\left(\frac{1}{p} - \frac{1}{\alpha+1} \right)\|u\|_W^q-\frac{\alpha+1-p}{p(\alpha+1)}-
\left(\frac{1}{1-\gamma} - \frac{1}{\alpha+1} \right)a_0^{-\frac{1-\gamma}{p}}\|f\|_{\frac{p}{p-1+\gamma}}\|u\|^{1-\gamma}_{W_a^{1,p}(V)}.
\end{eqnarray*}
By the fact that $\gamma<1<q\leq p<\alpha$, we can conclude that $J_\lambda(u)$ is coercive on $\mathcal{D}_\lambda$.\\
$(ii)$\;\;From $(i)$, $J_\lambda(u)$ is bounded from below on $\mathcal{D}_\lambda$. Then $\inf_{\mathcal{D}^+_\lambda \cup \mathcal{D}^0_\lambda} J_\lambda(u) \neq - \infty$. Besides, for $u \in \mathcal{D}^+_\lambda $, it follows that
\begin{eqnarray*}
&&J_\lambda(u)\\
&      =     & \left(\frac{1}{p} - \frac{1}{1-\gamma} \right)\|u\|^p_{W_a^{1,p}(V)}+\left(\frac{1}{q} - \frac{1}{1-\gamma} \right)\|u\|^q_{W_b^{1,q}(V)} -\left( \frac{1}{\alpha+1}-\frac{1}{1-\gamma} \right)\lambda\int_V g(x)u^{\alpha+1}d\mu\\
&      =    &\frac{1-\gamma-p}{p(1-\gamma)} \|u\|^p_{W_a^{1,p}(V)}+\frac{1-\gamma-q}{q(1-\gamma)} \|u\|^q_{W_b^{1,q}(V)} +\frac{\alpha+\gamma}{(1-\gamma)(\alpha+1)}  \lambda\int_V g(x)u^{\alpha+1}d\mu\\
&      <    & \frac{1-\gamma-p}{(\alpha+1)(1-\gamma)} \|u\|^p_{W_a^{1,p}(V)}+\frac{1-\gamma-q}{(\alpha+1)(1-\gamma)} \|u\|^q_{W_b^{1,q}(V)} +\frac{\alpha+\gamma}{(1-\gamma)(\alpha+1)} \lambda\int_V g(x)u^{\alpha+1}d\mu\\
&      =    & \frac{1}{(1-\gamma)(\alpha+1)} \left[  (1-\gamma-p)\|u\|^p_{W_a^{1,p}(V)}+(1-\gamma-q)\|u\|^q_{W_b^{1,q}(V)}+(\alpha+\gamma)\lambda\int_V g(x)u^{\alpha+1}d\mu\right]\\
&      \leq    &  0.
\end{eqnarray*}
Note that $\mathcal{D}^0_\lambda=\{0\}$ and $\mathcal{D}^+_\lambda\neq \emptyset$. Hence $\inf_{\mathcal{D}^+_\lambda \cup \mathcal{D}^0_\lambda} J_\lambda(u)\leq \inf_{\mathcal{D}^+_\lambda } J_\lambda(u)<0$. The proof is complete.
\qed

\vskip2mm
{\section{The multiplicity of positive solutions as $\lambda\in(0,\Lambda_*)$}}
  \setcounter{equation}{0}
  \par
In this section we prove Theorem 1.1 by Ekeland's variational principle.
\vskip2mm
\textbf{3.1. Solution in $\mathcal{D}^+_\lambda$ with $\lambda \in (0,\Lambda_{*})$}
 \vskip2mm
\noindent
\par
Using Ekeland's variational principle, there exists a minimizing sequence $\{u_k\} \subset \mathcal{D}^+_\lambda \cup \mathcal{D}^0_\lambda$ such that
\begin{eqnarray}
&&(i)\;\;J_\lambda(u_k) < \inf_{\mathcal{D}^+_\lambda \cup \mathcal{D}^0_\lambda} J_\lambda(u) +\frac{1}{k};\nonumber\\
&&(ii)\;\;J_\lambda(\tilde{u})\geq J_\lambda(u_k)-\frac{1}{k}\|\tilde{u}-u_k\|_{W}\;\; \mbox{for all}\;\; \tilde{u}\in \mathcal{D}^+_\lambda \cup \mathcal{D}^0_\lambda.\label{s6}
\end{eqnarray}
Since $\{u_k\}$ is a minimizing sequence, we can obtain $\lim_{k \rightarrow \infty} J_\lambda(u_k) = \inf_{\mathcal{D}^+_\lambda \cup \mathcal{D}^0_\lambda} J_\lambda$. By Lemma 2.8 $(ii)$ and $\mathcal{D}^0_\lambda=\{0\}$, we can assume that
$\{u_k\} \subset \mathcal{D}^+_\lambda$. By coerciveness of $J_\lambda(u)$, $\{u_k\}$ is bounded in ${W_+}$. Then there exists a positive constant $C_2$ such that $\|u_k\|_{W} \leq C_2$, for all $k=1,2,\cdots.$ Thus by Lemma 2.1, there exists $u_\lambda \in {W}$ such that, up to a subsequence,
\begin{eqnarray*}
\label{x42}
 \begin{cases}
u_k \rightharpoonup u_\lambda,& \mbox{in} \;\;W_+,\;\;\mbox{as}\;\;k \rightarrow \infty,\\
u_k(x)\rightarrow u_\lambda(x),& \forall \;\;x\in V,\;\;\mbox{as}\;\;k \rightarrow \infty.
   \end{cases}
\end{eqnarray*}
Since $W_+$ is closed and convex, we obtain $u_\lambda \in W_+$ and then $u_\lambda \geq 0$. Similar to (\ref{x21}) and (\ref{x22}), this implies that
\begin{eqnarray}
\label{x29}
\lim_{k\rightarrow \infty}\int_V f(x)u_k^{1-\gamma}d\mu =\int_V f(x)u_\lambda^{1-\gamma}d\mu
\end{eqnarray}
and
\begin{eqnarray}
\label{x30}
\lim_{k\rightarrow \infty}\int_V g(x)u_k^{\alpha+1}d\mu = \int_V g(x)u_\lambda^{\alpha+1}d\mu.
\end{eqnarray}
By the weakly lower semi-continuity of norm, (\ref{x29}) and (\ref{x30}), we infer that
$$
J_\lambda(u_\lambda) \leq \liminf_{k \rightarrow \infty} J_\lambda(u_k)=\inf_{\mathcal{D}^+_\lambda \cup \mathcal{D}^0_\lambda} J_\lambda<0
$$
Thus $u_\lambda \not\equiv 0$. Specifically, we will prove the existence of a solution $u_\lambda$ on $\mathcal{D}^+_\lambda$ through the following three steps.
\par
Step (I): we will claim that there exist a positive constant $C_3$ and a subsequence of $\{u_k\}$, still referred to it as  $\{u_k\}$, such that
\begin{eqnarray}
\label{x23}
(p-1+\gamma)\|u_k\|^p_{W_a^{1,p}(V)}+ (q-1+\gamma) \|u_k\|^q_{W_b^{1,q}(V)}- \lambda(\alpha+\gamma)\int_V g(x) u_k^{\alpha+1} d\mu \geq C_3.
\end{eqnarray}
Firstly, we prove that
\begin{eqnarray}
\label{x24}
\limsup_{k \rightarrow \infty}\left[(p-1+\gamma)\|u_k\|^p_{W_a^{1,p}(V)}+ (q-1+\gamma) \|u_k\|^q_{W_b^{1,q}(V)}\right] > \lambda(\alpha+\gamma)\int_V g(x) u_\lambda^{\alpha+1} d\mu.
\end{eqnarray}
In fact, by $\{u_k\} \subset \mathcal{D}^+_\lambda \cup \mathcal{D}^0_\lambda$ and (\ref{x30}), we deduce that
\begin{eqnarray*}
\limsup_{k \rightarrow \infty}\left[(p-1+\gamma)\|u_k\|^p_{W_a^{1,p}(V)}+ (q-1+\gamma) \|u_k\|^q_{W_b^{1,q}(V)}\right] & \geq  &\limsup_{k \rightarrow \infty} \left[\lambda(\alpha+\gamma)\int_V g(x) u_k^{\alpha+1} d\mu\right]\\
&  =  &\lambda(\alpha+\gamma)\int_V g(x) u_\lambda^{\alpha+1} d\mu.
\end{eqnarray*}
Suppose by contradiction that
\begin{eqnarray*}
\label{x25}
\limsup_{k \rightarrow \infty}\left[(p-1+\gamma)\|u_k\|^p_{W_a^{1,p}(V)}+ (q-1+\gamma) \|u_k\|^q_{W_b^{1,q}(V)}\right] = \lambda(\alpha+\gamma)\int_V g(x) u_\lambda^{\alpha+1} d\mu.
\end{eqnarray*}
Combining  $\{u_k\} \subset \mathcal{D}^+_\lambda \cup \mathcal{D}^0_\lambda$ with (\ref{x30}) gives
\begin{eqnarray*}
&&\liminf_{k \rightarrow \infty}\left[(p-1+\gamma)\|u_k\|^p_{W_a^{1,p}(V)}+ (q-1+\gamma) \|u_k\|^q_{W_b^{1,q}(V)}\right]\geq\lambda(\alpha+\gamma)\liminf_{k \rightarrow \infty} \int_V g(x)u_k^{\alpha+1}d\mu\\
&   \Rightarrow    &\liminf_{k \rightarrow \infty}\left[(p-1+\gamma)\|u_k\|^p_{W_b^{1,p}(V)}+ (q-1+\gamma) \|u_k\|^q_{W_b^{1,q}(V)}\right] \geq \lambda(\alpha+\gamma) \int_V g(x)u_\lambda^{\alpha+1}d\mu.
\end{eqnarray*}
Then we obtain that
\begin{eqnarray*}
&&\lambda(\alpha+\gamma) \int_V g(x)u_\lambda^{\alpha+1}d\mu
 \leq
\liminf_{k \rightarrow \infty}\left[(p-1+\gamma)\|u_k\|^p_{W_a^{1,p}(V)}+ (p-1+\gamma)\|u_k\|^q_{W_b^{1,q}(V)}\right]\\
&  \leq    &
\limsup_{k \rightarrow \infty}\left[(p-1+\gamma)\|u_k\|^p_{W_a^{1,p}(V)}+ (p-1+\gamma) \|u_k\|^q_{W_b^{1,q}(V)}\right]
 =
 \lambda(\alpha+\gamma) \int_V g(x)u_\lambda^{\alpha+1}d\mu,
\end{eqnarray*}
that is,
\begin{eqnarray}
\label{x3}
\lim_{k \rightarrow \infty}\left[(p-1+\gamma)\|u_k\|^p_{W_a^{1,p}(V)}+ (q-1+\gamma)\|u_k\|^q_{W_b^{1,q}(V)}\right]=\lambda(\alpha+\gamma) \int_V g(x)u_\lambda^{\alpha+1}d\mu.
\end{eqnarray}
From the fact that $u_\lambda\not\equiv0$, (\ref{x3}) implies that $\|u_k\|^p_{W_a^{1,p}(V)}$ converges to a positive number $W_1$ and $\|u_k\|^q_{W_b^{1,q}(V)}$ also converges to a positive number $W_2$, where $W_1$ and $W_2$ satisfy
\begin{eqnarray}
\label{x27}
(p-1+\gamma)W_1+ (q-1+\gamma) W_2 =  \lambda (\alpha+\gamma) \int_V g(x)u_\lambda^{\alpha+1}d\mu.
\end{eqnarray}
From (\ref{x29}) and (\ref{x30}), it follows that
\begin{eqnarray}
\label{x28}
&&\int_V f(x) u_\lambda^{1-\gamma} d\mu= \lim_{k\rightarrow \infty}\int_V f(x)u_k^{1-\gamma}d\mu\nonumber\\
&     =    &
\lim_{k\rightarrow \infty}\left[\|u_k\|^p_{W_a^{1,p}(V)} +\|u_k\|^q_{W_b^{1,q}(V)} -\lambda\int_V g(x)u_k^{\alpha+1}d\mu  \right]
=W_1 + W_2 -  \lambda\int_V g(x)u_\lambda^{\alpha+1}d\mu  \nonumber\\
&     =    &
\frac{\alpha+1-p}{\alpha+\gamma}W_1+  \frac{\alpha+1-q}{\alpha+\gamma}W_2.
\end{eqnarray}
Note that $Z(\lambda)>0$. By (\ref{x7}), (\ref{x27}) and (\ref{x28}), we calculate that
\begin{eqnarray*}
&&0\\
&       <      &
Z(\lambda)\|u_k\|^{1-\gamma}_{W_a^{1,p}(V)}\\
&\leq &\left(\frac{\alpha+1-p}{p-1+\gamma}\right)\cdot\left(\frac{\alpha+\gamma}{p-1+\gamma}\right)^{\frac{\alpha+\gamma}{p-\alpha-1}}\cdot
\left(\frac{\|u_k\|^{p(\alpha+\gamma)}_{W_a^{1,p}(V)}}{\left(\lambda\int_V g(x)u_k^{\alpha+1}d\mu\right)^{p+\gamma-1}}\right)^{\frac{1}{\alpha+1-p}}
-\int_V f(x)u_k^{1-\gamma}d\mu\\
&     \rightarrow      & \left(\frac{\alpha+1-p}{p-1+\gamma}\right)\cdot\left(\frac{\alpha+\gamma}{p-1+\gamma}\right)^{\frac{\alpha+\gamma}{p-\alpha-1}}\cdot
\left(\frac{(\alpha+\gamma)^{p-1+\gamma}W_1^{\alpha+\gamma}}{\left((p-1+\gamma)W_1+ (q-1+\gamma) W_2\right)^{p+\gamma-1}}\right)^{\frac{1}{\alpha+1-p}}\\
&&-
\frac{\alpha+1-p}{\alpha+\gamma}W_1- \frac{\alpha+1-q}{\alpha+\gamma}W_2\\
&       <      &
\left(\frac{\alpha+1-p}{p-1+\gamma}\right)\cdot\left(\frac{\alpha+\gamma}{p-1+\gamma}\right)^{\frac{\alpha+\gamma}{p-\alpha-1}}\cdot
\left(\frac{(\alpha+\gamma)^{p-1+\gamma}W_1^{\alpha+\gamma}}{\left((p-1+\gamma)W_1\right)^{p+\gamma-1}}\right)^{\frac{1}{\alpha+1-p}}-
\frac{\alpha+1-p}{\alpha+\gamma}W_1- \frac{\alpha+1-q}{\alpha+\gamma}W_2\\
&       =      &
- \frac{\alpha+1-q}{\alpha+\gamma}W_2<0.
\end{eqnarray*}
It is a contradiction. So (\ref{x24}) holds. Then we may pass to a subsequence such that (\ref{x23}) holds.
\par
Step (II): Applying Lemma 2.7 with $u=u_k$, $w=t\varphi, \varphi\in W^+$ and $t>0$ small, we can find $l_k(t) := l_k(t\varphi)>0$, $l_k(0)=1$ and $l_k(t)(u_k + t\varphi)\in \mathcal{D}^+_\lambda \subset \mathcal{D}_\lambda$ for all $\|t\varphi\|_{W} < \varepsilon$. We will show that $\{l'_{k^+}(0)\}$ is bounded for $k$ large enough.
\par
Since
\begin{eqnarray}
0  &=&   \|l_k(t)(u_k + t\varphi)\|^p_{W_a^{1,p}(V)}+\|l_k(t)(u_k + t\varphi)\|^q_{W_b^{1,q}(V)}\nonumber\\
&& -\int_V f(x)(l_k(t)(u_k + t\varphi))^{1-\gamma}d\mu-\lambda\int_V g(x)(l_k(t)(u_k + t\varphi))^{\alpha+1}d\mu\label{x31}
\end{eqnarray}
and
\begin{eqnarray}
0  &=&   \|u_k\|^p_{W_a^{1,p}(V)} +  \|u_k\|^q_{W_b^{1,p}(V)}-  \int_V f(x)u_k^{1-\gamma}d\mu -\lambda \int_V g(x)u_k^{\alpha+1}d\mu\label{x32}.
\end{eqnarray}
Combining (\ref{x31}) and (\ref{x32}) and noting that $u_k\geq 0$,  we obtain that $u_k+t\varphi\geq u_k$ and then
\begin{eqnarray*}
&&0\\
&     =      &
(l_k^p(t)-1)\|u_k + t\varphi\|^p_{W_a^{1,p}(V)} + \left(\|u_k + t\varphi\|^p_{W_a^{1,p}(V)}- \|u_k\|^p_{W_a^{1,p}(V)}\right)\\
&&+ (l_k^q(t)-1)\|u_k + t\varphi\|^q_{W_b^{1,q}(V)}+\left(\|u_k + t\varphi\|^q_{W_b^{1,q}(V)}- \|u_k\|^q_{W_b^{1,q}(V)}\right)\\
&&-(l_k^{1-\gamma}(t)-1)\int_V f(x)(u_k + t\varphi)^{1-\gamma}d\mu-\int_V f(x)[(u_k + t\varphi)^{1-\gamma}-u_k ^{1-\gamma}]d\mu\\
&&-\lambda (l_k^{\alpha+1}(t)-1)\int_V g(x)(u_k+ t\varphi)^{\alpha+1}d\mu - \lambda\int_V g(x)[(u_k+ t\varphi)^{\alpha+1}-u_k^{\alpha+1}]d\mu\\
&     \leq      &
(l_k^p(t)-1)\|u_k + t\varphi\|^p_{W_a^{1,p}(V)} + \left(\|u_k + t\varphi\|^p_{W_a^{1,p}(V)}- \|u_k\|^p_{W_a^{1,p}(V)}\right)\\
&&+ (l_k^q(t)-1)\|u_k + t\varphi\|^q_{W_b^{1,q}(V)}+\left(\|u_k + t\varphi\|^q_{W_b^{1,q}(V)}- \|u_k\|^q_{W_b^{1,q}(V)}\right)\\
&&-(l_k^{1-\gamma}(t)-1)\int_V f(x)(u_k + t\varphi)^{1-\gamma}d\mu\\
&&-\lambda (l_k^{\alpha+1}(t)-1)\int_V g(x)(u_k+ t\varphi)^{\alpha+1}d\mu -\lambda \int_V g(x)[(u_k+ t\varphi)^{\alpha+1}-u_k^{\alpha+1}]d\mu.
\end{eqnarray*}
Dividing by $t>0$ and letting $t \rightarrow 0^+$, by (\ref{x32}), we obtain that there exist $\theta_i\in (0,1), i=1,2\cdots7$ such that
\begin{eqnarray*}
\label{z1}
&&0 \nonumber\\
&     \leq     &
\lim_{t \rightarrow 0^+} \left[ \frac{l_k^p(t)-1}{t} \|u_k + t\varphi\|^p_{W_a^{1,p}(V)} + \frac{\|u_k + t\varphi\|^p_{W_a^{1,p}(V)}- \|u_k\|^p_{W_a^{1,p}(V)}}{t}\right.\nonumber\\
&&\left.+\frac{l_k^q(t)-1}{t} \|u_k + t\varphi\|^q_{W_b^{1,q}(V)}+\frac{\|u_k + t\varphi\|^q_{W_b^{1,q}(V)}- \|u_k\|^q_{W_b^{1,q}(V)}}{t}  -\frac{l_k^{1-\gamma}(t)-1}{t}  \int_V f(x)(u_k + t\varphi)^{1-\gamma}d\mu  \right.\nonumber\\
&&\left.-\frac{\lambda (l_k^{\alpha+1}(t)-1)}{t} \int_V g(x)(u_k+ t\varphi)^{\alpha+1}d\mu - \frac{ \lambda\int_V g(x)((u_k+ t\varphi)^{\alpha+1}-u_k^{\alpha+1})d\mu}{t} \right]\nonumber\nonumber\\
&     =     &
\lim_{t \rightarrow 0^+} \left[ pl_k^{p-1}(\theta_1 t)\cdot l'_{k}(\theta_1 t)\|u_k + t\varphi\|^p_{W_a^{1,p}(V)} + p\int_V [|\nabla (u_k +\theta_2 t\varphi)|^{p-2}\Gamma(u_k +\theta_2 t\varphi, \varphi)\right.\nonumber\\
&&\left.+a(x)(u_k +\theta_2 t\varphi)^{p-1}\varphi]d\mu+ ql_k^{q-1}(\theta_3 t)\cdot l'_{k}(\theta_3 t)\|u_k + t\varphi\|^q_{W_b^{1,q}(V)}\right.\nonumber\\
&&\left.+q\int_V  [|\nabla (u_k +\theta_4 t\varphi)|^{q-2}\Gamma(u_k +\theta_4 t\varphi, \varphi)+b(x)(u_k +\theta_4 t\varphi)^{q-1}\varphi]d\mu\right.\nonumber\\
&&\left.-(1-\gamma)l_k^{-\gamma}(\theta_5 t)\cdot l'_{k}(\theta_5 t)\int_V f(x)(u_k + t\varphi)^{1-\gamma}d\mu
\right.\nonumber\\
&&\left.-\lambda(\alpha+1)l_k^{\alpha}(\theta_6 t)\cdot l'_{k}(\theta_6 t)\int_V g(x)(u_k+ t\varphi)^{\alpha+1}d\mu
-\lambda(\alpha+1)\int_V g(x)(u_k+\theta_7 t\varphi)^{\alpha}\varphi d\mu
\right]\nonumber\\
&     =     &
pl'_{k^+}(0)\|u_k \|^p_{W_a^{1,p}(V)}+p\int_V (|\nabla u_k|^{p-2}\Gamma(u_k , \varphi)+a(x)u_k ^{p-1} \varphi)d\mu\nonumber\\
&&+ql'_{k^+}(0)\|u_k \|^q_{W_b^{1,q}(V)}+q\int_V (|\nabla u_k|^{q-2}\Gamma(u_k, \varphi)+b(x)u_k ^{q-1} \varphi)d\mu\nonumber\\
&&-(1-\gamma)l'_{k^+}(0)\int_V f(x)u_k ^{1-\gamma}d\mu - \lambda(\alpha+1) l'_{k}(0)\int_V g(x)u_k^{\alpha+1}d\mu
-\lambda(\alpha+1)\int_V g(x)u_k^{\alpha}\varphi d\mu\nonumber\\
&     =     &
l'_{k^+}(0)\left[p\|u_k \|^p_{W_a^{1,p}(V)}+q \|u_k \|^q_{W_b^{1,q}(V)}-(1-\gamma) \int_V f(x)u_k ^{1-\gamma}d\mu - \lambda(\alpha+1) \int_V g(x)u_k^{\alpha+1}d\mu \right]\nonumber\\
&&+p\int_V (|\nabla u_k|^{p-2}\Gamma(u_k , \varphi)+a(x)u_k ^{p-1}  \varphi)d\mu+q\int_V (|\nabla u_k|^{q-2}\Gamma(u_k, \varphi)+b(x)u_k ^{q-1}  \varphi)d\mu\nonumber\\
&&-\lambda(\alpha+1)\int_V g(x)u_k^{\alpha}\varphi d\mu\nonumber\\
&     =     &
l'_{k^+}(0)\left[(p-1+\gamma)\|u_k\|^p_{W_a^{1,p}(V)}+ (q-1+\gamma) \|u_k\|^q_{W_b^{1,q}(V)}- \lambda(\alpha+\gamma)\int_V g(x) u_k^{\alpha+1} d\mu\right]\nonumber\\
&&+p\int_V (|\nabla u_k|^{p-2}\Gamma(u_k , \varphi)+a(x)u_k ^{p-1}  \varphi)d\mu+q\int_V  (|\nabla u_k|^{q-2}\Gamma(u_k, \varphi)+b(x)u_k ^{q-1}  \varphi)d\mu\nonumber\\
&&-\lambda(\alpha+1)\int_V g(x)u_k^{\alpha} \varphi d\mu,
\end{eqnarray*}
where $l'_{k^+}(0)$ is the right-hand derivative with respect to $t$. \\
Let
\begin{eqnarray*}
H_1(u_k)&=&p\int_V (|\nabla u_k|^{p-2}\Gamma(u_k , \varphi)+a(x)u_k ^{p-1}  \varphi)d\mu+q\int_V  (|\nabla u_k|^{q-2}\Gamma(u_k, \varphi)+b(x)u_k ^{q-1}  \varphi)d\mu\\
&&-\lambda(\alpha+1)\int_V g(x)u_k^{\alpha} \varphi d\mu.
\end{eqnarray*}
The combination of Appendix A.3, Lemma 2.1 and $\|u_k\|_W\leq C_2$ leads to
\begin{eqnarray}\label{q1}
&&H_1(u_k)\nonumber\\
&\leq &p\int_V (|\nabla u_k|^{p-2}\Gamma(u_k , \varphi)+a(x)u_k ^{p-1}  \varphi)d\mu+q\int_V  (|\nabla u_k|^{q-2}\Gamma(u_k, \varphi)+b(x)u_k ^{q-1}  \varphi)d\mu\nonumber\\
&&+\lambda(\alpha+1)\int_V g(x)u_k^{\alpha} \varphi d\mu\nonumber\\
&\leq &p\|u_k\|^{p-1}_{W_a^{1,p}(V)}\|\varphi\|_{W_a^{1,p}(V)}+ q\|u_k\|^{q-1}_{W_b^{1,q}(V)}\|\varphi\|_{W_b^{1,q}(V)} +\lambda(\alpha+1)\|g\|_\infty \left(\frac{C_2}{ \mu_0^{\frac{1}{p}} a_0^{\frac{1}{p}} +\mu_0^{\frac{1}{q}}b_0^{\frac{1}{q}}}\right)^{\alpha-p+1}\|u_k\|_p^{p-1}\|\varphi\|_p\nonumber\\
&\leq &pC_2^{p-1}\|\varphi\|_{W_a^{1,p}(V)}+ qC_2^{q-1}\|\varphi\|_{W_b^{1,q}(V)} +\lambda(\alpha+1)\|g\|_\infty \left(\frac{C_2}{ \mu_0^{\frac{1}{p}} a_0^{\frac{1}{p}} +\mu_0^{\frac{1}{q}}b_0^{\frac{1}{q}}}\right)^{\alpha-p+1}a_0^{-\frac{p-1}{p}}C_2^{p-1}\|\varphi\|_p\nonumber\\
&:=&\widetilde{C_2}.
\end{eqnarray}
Let
\begin{eqnarray*}
H_2(u_k)=(p-1+\gamma)\|u_k\|^p_{W_a^{1,p}(V)}+ (q-1+\gamma) \|u_k\|^q_{W_b^{1,q}(V)}- (\alpha+\gamma)\lambda\int_V g(x) u_k^{\alpha+1} d\mu.
\end{eqnarray*}
Then (\ref{x23}) and (\ref{q1}) imply that
\begin{eqnarray}
\label{x33}
l'_{k^+}(0) \geq -\frac{H_1(u_k)}{H_2(u_k)}\geq -\frac{\widetilde{C_2}}{C_3}.
\end{eqnarray}
Thus, we conclude that $\{l'_{k}(0)\}$ is bounded from below.
\par
Next, we claim that there exists $C_4 >0$ such that $l'_{k^+}(0)<C_4$. Without loss of generality, we assume $l'_{k^+}(0) \geq0$.
Since $u_k\in \mathcal{D}^+_\lambda \subset \mathcal{D}_\lambda$, $l_k(t)(u_k + t\varphi)\in \mathcal{D}^+_\lambda \subset \mathcal{D}_\lambda$ and (\ref{s6}), we obtain that
\begin{eqnarray*}
&&[l_k(t)-1]\frac{\|u_k\|_{W}}{k}+tl_k(t)\frac{\|\varphi\|_{W}}{k}\\
&    \geq   &
\frac{1}{k}\|l_k(t)(u_k+t\varphi)-u_k\|_{W}\\
&    \geq   &
J_\lambda(u_k)-J_\lambda(l_k(t)(u_k+t\varphi))\\
&    =   &
\frac{1}{p}\|u_k\|^p_{W_a^{1,p}(V)}+\frac{1}{q}\|u_k\|^q_{W_b^{1,q}(V)}-\frac{1}{1-\gamma}\int_V f(x)u_k^{1-\gamma}d\mu-\frac{\lambda}{\alpha+1}\int_V g(x)u_k^{\alpha+1}d\mu\\
&& -\frac{1}{p}\|l_k(t)(u_k+t\varphi)\|^p_{W_a^{1,p}(V)}-\frac{1}{q}\|l_k(t)(u_k+t\varphi)\|^q_{W_b^{1,q}(V)}\\
&&+\frac{1}{1-\gamma}\int_V f(x)(l_k(t)(u_k+t\varphi))^{1-\gamma}d\mu+\frac{\lambda}{\alpha+1}\int_V g(x)(l_k(t)(u_k+t\varphi))^{\alpha+1}d\mu\\
&    =   &
\frac{1}{p}\|u_k\|^p_{W_a^{1,p}(V)}+\frac{1}{q}\|u_k\|^q_{W_b^{1,q}(V)}
-\frac{1}{1-\gamma}\left[ \|u_k\|^p_{W_a^{1,p}(V)} + \|u_k\|^q_{W_b^{1,q}(V)}- \lambda\int_V g(x)u_k^{\alpha+1}d\mu \right]\\
&&-\frac{\lambda}{\alpha+1}\int_V g(x)u_k^{\alpha+1}d\mu-\frac{1}{p}\|l_k(t)(u_k+t\varphi)\|^p_{W_a^{1,p}(V)}-\frac{1}{q}\|l_k(t)(u_k+t\varphi)\|^q_{W_b^{1,q}(V)}\\
&&+\frac{\lambda}{\alpha+1}\int_V g(x)(l_k(t)(u_k+t\varphi))^{\alpha+1}d\mu
+\frac{1}{1-\gamma} \left[\|l_k(t)(u_k+t\varphi)\|^p_{W_a^{1,p}(V)}\right.\\
&&\left. + \|l_k(t)(u_k+t\varphi)\|^q_{W_b^{1,q}(V)}- \lambda\int_V g(x)(l_k(t)(u_k+t\varphi))^{\alpha+1}d\mu\right]\\
&    =   &
\left(\frac{1}{p} - \frac{1}{1-\gamma} \right)\|u_k\|^p_{W_a^{1,p}(V)}+\left(\frac{1}{q} - \frac{1}{1-\gamma}\right)\|u_k\|^q_{W_b^{1,q}(V)}-\lambda\left(\frac{1}{\alpha+1}- \frac{1}{1-\gamma} \right)\int_V g(x)u_k^{\alpha+1}d\mu\\
&&
-\left(\frac{1}{p} - \frac{1}{1-\gamma} \right)\|l_k(t)(u_k+t\varphi)\|^p_{W_a^{1,p}(V)}-\left(\frac{1}{q} - \frac{1}{1-\gamma}\right)\|l_k(t)(u_k+t\varphi)\|^q_{W_b^{1,q}(V)}\\
&&+\lambda\left(\frac{1}{\alpha+1}- \frac{1}{1-\gamma} \right)\int_V g(x)(l_k(t)(u_k+t\varphi))^{\alpha+1}d\mu\\
&    =   &
\left(\frac{1}{p} - \frac{1}{1-\gamma} \right) \left(\|u_k\|^p_{W_a^{1,p}(V)}- \|u_k+t\varphi\|^p_{W_a^{1,p}(V)}\right) + \left(\frac{1}{p} - \frac{1}{1-\gamma}\right)(1-l_k^p(t))\|u_k+t\varphi\|^p_{W_a^{1,p}(V)}\\
&&+\left(\frac{1}{q} - \frac{1}{1-\gamma} \right)\left(\|u_k\|^q_{W_b^{1,q}(V)}- \|u_k+t\varphi\|^q_{W_b^{1,q}(V)}\right)+ \left(\frac{1}{q} - \frac{1}{1-\gamma}\right)(1-l_k^q(t))\|u_k+t\varphi\|^q_{W_b^{1,q}(V)}\\
&&-\lambda\left(\frac{1}{\alpha+1}- \frac{1}{1-\gamma} \right)\int_V g(x)(u_k^{\alpha+1}-(u_k+t\varphi)^{\alpha+1})d\mu\\
&&-\lambda\left(\frac{1}{\alpha+1}- \frac{1}{1-\gamma} \right)(1-l_k^{\alpha+1}(t))\int_V g(x)(u_k+t\varphi)^{\alpha+1}d\mu\\
&    =   &
\frac{p+\gamma-1}{p(1-\gamma)} \left(\|u_k+t\varphi\|^p_{W_a^{1,p}(V)}- \|u_k\|^p_{W_a^{1,p}(V)}\right)
+\frac{p+\gamma-1}{p(1-\gamma)}(l_k^p(t)-1)\|u_k+t\varphi\|^p_{W_a^{1,p}(V)}\\
&&+\frac{q+\gamma-1}{q(1-\gamma)} \left( \|u_k+t\varphi\|^q_{W_b^{1,q}(V)}-\|u_k\|^q_{W_b^{1,q}(V) }\right)+\frac{q+\gamma-1}{q(1-\gamma)} (l_k^q(t)-1)\|u_k+t\varphi\|^q_{W_b^{1,q}(V)}\\
&&-\lambda\frac{\alpha+\gamma}{(\alpha+1)(1-\gamma)} \int_V g(x)((u_k+t\varphi)^{\alpha+1}-u_k^{\alpha+1})d\mu
-\lambda\frac{\alpha+\gamma}{(\alpha+1)(1-\gamma)}(l_k^{\alpha+1}(t)-1) \int_V g(x)(u_k+t\varphi)^{\alpha+1}d\mu.
\end{eqnarray*}
Dividing by $t>0$ and passing to the limit as $t \rightarrow 0^+$ gives
\begin{eqnarray*}
&& \lim_{t \rightarrow 0^+} \left[ \frac{l_k(t)-l_k(0)}{t}\frac{\|u_k\|_{W}}{k}+ l_k(t)\frac{\|\varphi\|_{W}}{k} \right]\\
&     =     &
l'_{k^+}(0) \frac{\|u_k\|_{W}}{k} + \frac{\|\varphi\|_{W}}{k}\\
&     \geq     &
\frac{p+\gamma-1}{1-\gamma} \int_V(|\nabla u_k|^{p-2}\Gamma(u_k,\varphi)+a(x)|u_k|^{p-2}u_k \varphi)d\mu + \frac{p+\gamma-1}{1-\gamma}l'_{k^+}(0)\|u_k\|_{W_a^{1,p}(V)}\\
&&+\frac{q+\gamma-1}{1-\gamma} \int_V (|\nabla u_k|^{q-2}\Gamma(u_k,\varphi)+b(x)|u_k|^{q-2}u_k \varphi)d\mu+ \frac{q+\gamma-1}{1-\gamma} l'_{k^+}(0) \|u_k\|_{W_b^{1,q}(V)}\\
&&-\lambda\frac{\alpha+\gamma}{1-\gamma} \int_V g(x)u_k^{\alpha}\varphi d\mu
-\lambda\frac{\alpha+\gamma}{1-\gamma}l'_{k^+}(0) \int_V g(x)u_k^{\alpha+1}d\mu.
\end{eqnarray*}
According to the above inequality, (\ref{x23}), (\ref{x32}), $\{u_k\} \subset \mathcal{D}^+_\lambda$ and the fact $\frac{(1-\gamma)C_2}{k} < \frac{C_3}{2}$ for $k$ large enough, we deduce that
\begin{eqnarray*}
&&\frac{\|\varphi\|_{W}}{k}\\
&    \geq     &
\frac{l'_{k^+}(0)}{1-\gamma}\left[ (p+\gamma-1)\|u_k\|^p_{W_a^{1,p}(V)}+ (q+\gamma-1) \|u_k\|^q_{W_b^{1,q}(V)}-\lambda (\alpha+\gamma) \int_V g(x)u_k^{\alpha+1}d\mu\right.\\
&&\left.- (1-\gamma)\frac{\|u_k\|_{W}}{k} \right]+\frac{p+\gamma-1}{1-\gamma} \int_V(|\nabla u_k|^{p-2}\Gamma(u_k,\varphi)+a(x)u_k^{p-1} \varphi)d\mu\\
&&+\frac{q+\gamma-1}{1-\gamma} \int_V (|\nabla u_k|^{q-2}\Gamma(u_k,\varphi)+b(x)u_k^{q-1} \varphi)d\mu-\lambda\frac{\alpha+\gamma}{1-\gamma} \int_V g(x)u_k^{\alpha}\varphi d\mu\\
&    \geq    &
\frac{C_3}{2(1-\gamma)}l'_{k^+}(0)+\frac{p+\gamma-1}{1-\gamma} \int_V(|\nabla u_k|^{p-2}\Gamma(u_k,\varphi)+a(x)u_k^{p-1} \varphi)d\mu\\
&&+\frac{q+\gamma-1}{1-\gamma}  \int_V (|\nabla u_k|^{q-2}\Gamma(u_k,\varphi)+b(x)u_k^{q-1} \varphi)d\mu-\lambda\frac{\alpha+\gamma}{1-\gamma} \int_V g(x)u_k^{\alpha} \varphi d\mu.
\end{eqnarray*}
It follows from Appendix A.3 that
\begin{eqnarray*}
\frac{C_3}{2(1-\gamma)}l'_{k^+}(0) &\leq&
\frac{\|\varphi\|_{W}}{k}-
\frac{p+\gamma-1}{1-\gamma}\int_V(|\nabla u_k|^{p-2}\Gamma(u_k,\varphi)+a(x)u_k^{p-1} \varphi)d\mu\\
&&-\frac{q+\gamma-1}{1-\gamma}\int_V (|\nabla u_k|^{q-2}\Gamma(u_k,\varphi)+b(x)u_k^{q-1} \varphi)d\mu+\lambda\frac{\alpha+\gamma}{1-\gamma}\int_V g(x)u_k^{\alpha} \varphi d\mu\\
&    \leq     &\frac{\|\varphi\|_{W}}{k}+ \frac{p+\gamma-1}{1-\gamma}\|u_k\|^{p-1}_{W_a^{1,p}(V)}\|\varphi\|_{W_a^{1,p}(V)}+\frac{q+\gamma-1}{1-\gamma} \|u_k\|^{q-1}_{W_b^{1,q}(V)}\|\varphi\|_{W_b^{1,q}(V)}\\
&&+\lambda\frac{\alpha+\gamma}{1-\gamma} \|g\|_\infty \left(\frac{C_2}{ \mu_0^{\frac{1}{p}} a_0^{\frac{1}{p}} +\mu_0^{\frac{1}{q}}b_0^{\frac{1}{q}}}\right)^{\alpha-p+1}\|u_k\|_p^{p-1}\|\varphi\|_p.
\end{eqnarray*}
Due to the boundedness of $\{u_k\}$, we conclude that $\{l'_{k^+}(0)\} $ is bounded from above. Thus, combining with (\ref{x33}), we obtain $\{l'_{k^+}(0)\}$ is bounded for $k$ large enough.
\par
Step (III): we will prove $u_\lambda$ is a weak solution of (\ref{eq1}).
Firstly, we claim $u_k \rightarrow u_\lambda$ in $W_+$. By Lemma 2.3, we can find a unique $t_1(u_\lambda)\in(0,+\infty)$ such that $t_1(u_\lambda)u_\lambda\in \mathcal{D}^+_\lambda$, i.e, $\varphi_{u_\lambda}(t_1(u_\lambda))=0$. Suppose $\|u_\lambda\|_W < \liminf_{k\rightarrow \infty}\|u_k\|_W$. Thus we have $\liminf_{k\rightarrow \infty} J_\lambda(u_k)>J_\lambda(u_\lambda)$ and
\begin{eqnarray}
\label{c2}
\liminf_{k\rightarrow \infty} \varphi_{u_k}(t_1(u_\lambda)) > \varphi_{u_\lambda}(t_1(u_\lambda))=0.
\end{eqnarray}
Since $u_k\in \mathcal{D}^+_\lambda$, we get that $t_1(u_k)=1$, $\varphi_{u_k}(t_1(u_k))=0$ and $\varphi_{u_k}'(t_1(u_k))>0$. Then combining (\ref{c2}), we obtain that $t_1(u_k)< t_1(u_\lambda)$ for all large $k$. Since $\varphi_{u_\lambda}(t_1(u_\lambda))=0$, from Figure 1 in Remark 2.2, we get $\varphi_{u_\lambda}(t)\leq0$ in $[t_1(u_k), t_1(u_\lambda)]$, that is, $J_\lambda (t_1(u_\lambda)u_\lambda)< J_\lambda (u_\lambda)$. Therefore
\begin{eqnarray*}
\inf_{\mathcal{D}^+_\lambda \cup \mathcal{D}^0_\lambda} J_\lambda \leq J_\lambda (t_1(u_\lambda)u_\lambda)<J_\lambda (u_\lambda)<
\liminf_{k\rightarrow \infty}J_\lambda (u_k)=\inf_{\mathcal{D}^+_\lambda \cup \mathcal{D}^0_\lambda} J_\lambda.
\end{eqnarray*}
It is an absurd. Thus $\|u_\lambda\|_W = \liminf_{k\rightarrow \infty}\|u_k\|_W$. From Lemma 2.2, we conclude $u_k \rightarrow u_\lambda$ in $W_+$. Furthermore, by the fact that $u_k \in \mathcal{D}^+_\lambda $, (\ref{x1}), (\ref{x29}), (\ref{x30}) and (\ref{x24}), it is easy to obtain that $u_\lambda \in \mathcal{D}^+_\lambda$.
\par
Applying (\ref{s6}) again, we infer that
\begin{eqnarray*}
&&\frac{1}{k}\left[|l_k(t)-1|\|u_k\|_{W}+tl_k(t)\|\varphi\|_{W}\right]\\
&     \geq      &
J_\lambda(u_k)-J_\lambda(l_k(t)(u_k+t\varphi))\\
&     =       &
-\left(\frac{l^p_k(t)-1}{p}\right)\|u_k\|^p_{W_a^{1,p}(V)}-\left(\frac{l^q_k(t)-1}{q}\right)\|u_k\|^q_{W_b^{1,q}(V)}\\
&&+\lambda\left(\frac{l^{\alpha+1}_k(t)-1}{\alpha+1}\right)\int_V g(x)(u_k+t\varphi)^{\alpha+1}d\mu+\frac{l^{1-\gamma}_k(t)-1}{1-\gamma}\int_V f(x)(u_k+t\varphi)^{1-\gamma}d\mu\\
&&+\frac{l^p_k(t)}{p}(\|u_k\|^p_{W_a^{1,p}(V)}-\|u_k+t\varphi\|^p_{W_a^{1,p}(V)})       +\frac{l^q_k(t)}{q}(\|u_k\|^q_{W_b^{1,q}(V)}-\|u_k+t\varphi\|^q_{W_b^{1,q}(V)})\\
&&+\frac{\lambda}{\alpha+1}\int_V g(x)[(u_k+t\varphi)^{\alpha+1}-u_k^{\alpha+1}]d\mu+\frac{1}{1-\gamma}\int_Vf(x)[(u_k+t\varphi)^{1-\gamma}-u_k^{1-\gamma}]d\mu.
\end{eqnarray*}
Dividing by $t>0$ and passing to the limit as $t \rightarrow 0^+$. Since $u_k \in \mathcal{D}_\lambda$, $l_k(0)=1$ and $f(x)[(u_k+t\varphi)^{1-\gamma}-u_k^{1-\gamma}] \geq 0$ for all $x \in V$ and $t>0$, it follows from Fatou's Lemma that there exists $\theta \in (0,1)$ such that
\begin{eqnarray}
\label{x34}
&&\frac{1}{k}\left[|l'_{k^+}(0)|\|u_k\|_{W}+\|\varphi\|_{W}\right]\nonumber\\
&    \geq     &
-l'_{k^+}(0)\left[\|u_k\|^p_{W_a^{1,p}(V)} + \|u_k\|^q_{W_b^{1,q}(V)} - \lambda\int_V g(x)u_k^{\alpha+1}d\mu-\int_V f(x)u_k^{1-\gamma}d\mu \right]\nonumber\\
&&-\int_V(|\nabla u_k|^{p-2}\Gamma(u_k,\varphi)+a(x)u_k^{p-1}\varphi)d\mu
-\int_V (|\nabla u_k|^{q-2}\Gamma(u_k,\varphi)+b(x)u_k^{q-1}\varphi)d\mu\nonumber\\
&&+\lambda\int_V g(x)u_k^{\alpha} \varphi d\mu
+\liminf_{t\rightarrow 0^+} \frac{1}{1-\gamma}\int_V \frac{f(x)[(u_k+t\varphi)^{1-\gamma}-u_k^{1-\gamma}]}{t}d\mu\nonumber\\
&    =       &
-\int_V(|\nabla u_k|^{p-2}\Gamma(u_k,\varphi)+a(x)u_k^{p-1}\varphi)d\mu
-\int_V (|\nabla u_k|^{q-2}\Gamma(u_k,\varphi)+b(x)u_k^{q-1} \varphi)d\mu\nonumber\\
&&+\lambda\int_V g(x)u_k^{\alpha}\varphi d\mu
+\liminf_{t\rightarrow 0^+} \frac{1}{1-\gamma}\int_V \frac{f(x)[(u_k+t\varphi)^{1-\gamma}-u_k^{1-\gamma}]}{t}d\mu\nonumber\\
&    \geq     &
-\int_V(|\nabla u_k|^{p-2}\Gamma(u_k,\varphi)+a(x)u_k^{p-1}\varphi)d\mu
-\int_V (|\nabla u_k|^{q-2}\Gamma(u_k,\varphi)+b(x)u_k^{q-1} \varphi)d\mu\nonumber\\
&&+\lambda\int_V g(x)u_k^{\alpha}\varphi d\mu
+\int_V \liminf_{t\rightarrow 0^+} f(x)(u_k+\theta t\varphi)^{-\gamma}\varphi d\mu\nonumber\\
&    =       &
-\int_V(|\nabla u_k|^{p-2}\Gamma(u_k,\varphi)+a(x)u_k^{p-1}\varphi)d\mu
-\int_V (|\nabla u_k|^{q-2}\Gamma(u_k,\varphi)+b(x)u_k^{q-1} \varphi)d\mu\nonumber\\
&&+\lambda\int_V g(x)u_k^{\alpha}\varphi d\mu
+\int_V f(x) u_k^{-\gamma} \varphi d\mu.
\end{eqnarray}
Since $\{l'_{k^+}(0)\}$ and $\{u_k\}$ are bounded, we deduce that
\begin{eqnarray*}
\label{c1}
\int_V f(x) u_k^{-\gamma} \varphi d\mu
&  \leq    &\frac{C_2\cdot \max\{\frac{\widetilde{C_2}}{C_3}, C_4\}+\|\varphi\|_{W}}{k}
+\int_V(|\nabla u_k|^{p-2}\Gamma(u_k,\varphi)+a(x)u_k^{p-1}\varphi)d\mu\nonumber\\
&&+\int_V (|\nabla u_k|^{q-2}\Gamma(u_k,\varphi)+b(x)u_k^{q-1} \varphi)d\mu
-\lambda\int_V g(x)u_k^{\alpha} \varphi d\mu.
\end{eqnarray*}
Passing to the limit as $k \rightarrow \infty$ and applying Appendix A.4 along with Fatou's Lemma, we obtain, for any $\varphi\in W_+$, that
\begin{eqnarray}\label{x38}
&&\int_V f(x) u_\lambda^{-\gamma}\varphi d\mu\nonumber\\
&  \leq  &\liminf_{k\rightarrow \infty} \int_V f(x) u_k^{-\gamma}\varphi d\mu\nonumber\\
&\leq&
\int_V(|\nabla u_\lambda|^{p-2}\Gamma(u_\lambda,\varphi)+a(x)u_\lambda^{p-1} \varphi ) d\mu+\int_V (|\nabla u_\lambda|^{q-2}\Gamma(u_\lambda,\varphi)+b(x) u_\lambda^{q-1} \varphi)d\mu
-\lambda\int_V g(x) u_\lambda^{\alpha} \varphi d\mu.
\end{eqnarray}
For any given $\phi \in W$ and $\varepsilon>0$, we define $\Phi\in W_+$ by $\Phi= (u_\lambda+ \varepsilon \phi)^+$. Inserting $\Phi$ into (\ref{x38}) and utilizing the fact that $u_\lambda\in \mathcal{D}^+_\lambda$, it follows that
\begin{eqnarray}\label{w2}
0&  \leq   &
\int_V(|\nabla u_\lambda|^{p-2}\Gamma(u_\lambda,\Phi)+a(x)u_\lambda^{p-1} \Phi ) d\mu+
\int_V(|\nabla u_\lambda|^{q-2}\Gamma(u_\lambda,\Phi)+ b(x) u_\lambda^{q-1}\Phi)d\mu\nonumber\\
&&-\int_V f(x) u_\lambda^{-\gamma}\Phi d\mu
-\lambda\int_V g(x) u_\lambda^{\alpha}\Phi d\mu\nonumber\\
&    =     &
\int_{\{x\in V: u_\lambda+ \varepsilon \phi\geq0\}}\left[|\nabla u_\lambda|^{p-2}\Gamma(u_\lambda,u_\lambda+ \varepsilon \phi)+a(x)u_\lambda^{p-1}(u_\lambda+ \varepsilon \phi)- f(x) u_\lambda^{-\gamma}(u_\lambda+ \varepsilon \phi) \right.\nonumber\\
&& \left. +
(|\nabla u_\lambda|^{q-2}\Gamma(u_\lambda,u_\lambda+ \varepsilon \phi)+ b(x) u_\lambda^{q-1} (u_\lambda+ \varepsilon \phi))
-\lambda  g(x) u_\lambda^{\alpha} (u_\lambda+ \varepsilon \phi) \right]  d\mu\nonumber\\
&    =     &
\left(\int_{V}-\int_{\{x\in V: u_\lambda+ \varepsilon \phi<0\}}\right)\left[|\nabla u_\lambda|^{p-2}\Gamma(u_\lambda,u_\lambda+ \varepsilon \phi)+a(x)u_\lambda^{p-1}(u_\lambda+ \varepsilon \phi)- f(x) u_\lambda^{-\gamma}(u_\lambda+ \varepsilon \phi) \right.\nonumber\\
&& \left. +
 (|\nabla u_\lambda|^{q-2}\Gamma(u_\lambda,u_\lambda+ \varepsilon \phi)+ b(x) u_\lambda^{q-1} (u_\lambda+ \varepsilon \phi))
-\lambda  g(x) u_\lambda^{\alpha} (u_\lambda+ \varepsilon \phi) \right]  d\mu\nonumber\\
&    =     &
\|u_\lambda\|^p_{W_a^{1,p}(V)} +\|u_\lambda\|^q_{W_b^{1,q}(V)}
- \int_V f(x) u_\lambda^{1-\gamma}d\mu- \lambda\int_V g(x) u_\lambda^{\alpha+1} d\mu \nonumber\\
&&+ \varepsilon \int_V \left[ |\nabla u_\lambda|^{p-2}\Gamma(u_\lambda,\phi)+a(x) u_\lambda^{p-1}\phi
+ |\nabla u_\lambda|^{q-2}\Gamma(u_\lambda,\phi)+b(x) u_\lambda^{q-1} \phi - f(x) u_\lambda^{-\gamma} \phi -\lambda g(x) u_\lambda^{\alpha} \phi \right]d\mu\nonumber\\
&&-\int_{\{x\in V: u_\lambda+ \varepsilon \phi<0\}}\left[|\nabla u_\lambda|^{p-2}\Gamma(u_\lambda,u_\lambda+ \varepsilon \phi)+a(x)u_\lambda^{p-1}(u_\lambda+ \varepsilon \phi)\right.\nonumber\\
&&\left. +
|\nabla u_\lambda|^{q-2}\Gamma(u_\lambda,u_\lambda+ \varepsilon \phi)+ b(x) u_\lambda^{q-1} (u_\lambda+ \varepsilon \phi)-f(x) u_\lambda^{-\gamma}(u_\lambda+ \varepsilon \phi)
-\lambda g(x) u_\lambda^{\alpha} (u_\lambda+ \varepsilon \phi)\right] d\mu\nonumber\\
&    \leq     &
\varepsilon \int_V \left[ |\nabla u_\lambda|^{p-2}\Gamma(u_\lambda,\phi)+a(x) u_\lambda^{p-1}\phi
+ |\nabla u_\lambda|^{q-2}\Gamma(u_\lambda,\phi)+b(x) u_\lambda^{q-1} \phi  - f(x) u_\lambda^{-\gamma} \phi -\lambda g(x) u_\lambda^{\alpha} \phi \right]d\mu\nonumber\\
&&-\varepsilon \int_{\{x\in V: u_\lambda+ \varepsilon \phi<0\}} \left[ |\nabla u_\lambda|^{p-2}\Gamma(u_\lambda,\phi)+a(x) u_\lambda^{p-1}\phi + |\nabla u_\lambda|^{q-2}\Gamma(u_\lambda,\phi)+b(x) u_\lambda^{q-1} \phi
\right]d\mu.
\end{eqnarray}
Since $u_\lambda \geq 0$, the measure of the domain of integration $\{x\in V: u_\lambda+ \varepsilon \phi<0\}$ tend to zero as $\varepsilon \rightarrow 0$. Then $\int_{\{x\in V: u_\lambda+ \varepsilon \phi<0\}} \left[ |\nabla u_\lambda|^{p-2}\Gamma(u_\lambda,\phi) +a(x) u_\lambda^{p-1}\phi + |\nabla u_\lambda|^{q-2}\Gamma(u_\lambda,\phi)+b(x) u_\lambda^{q-1} \phi
\right]d\mu\rightarrow 0$ as $\varepsilon \rightarrow 0$. Thus dividing by $\varepsilon$ and letting $\varepsilon \rightarrow 0$ in (\ref{w2}), we obtain
\begin{eqnarray*}\label{b5}
\int_V \left[ |\nabla u_\lambda|^{p-2}\Gamma(u_\lambda,\phi)+a(x) u_\lambda^{p-1}\phi
+ |\nabla u_\lambda|^{q-2}\Gamma(u_\lambda,\phi)+b(x) u_\lambda^{q-1} \phi  - f(x) u_\lambda^{-\gamma} \phi -\lambda g(x) u_\lambda^{\alpha} \phi \right]d\mu \geq 0.
\end{eqnarray*}
Since $\phi$ is arbitrary, the inequality also holds for any given $-\phi$. Then there holds that
\begin{eqnarray}\label{b5}
\int_V \left[ |\nabla u_\lambda|^{p-2}\Gamma(u_\lambda,\phi)+a(x) u_\lambda^{p-1}\phi
+ |\nabla u_\lambda|^{q-2}\Gamma(u_\lambda,\phi)+b(x) u_\lambda^{q-1} \phi  - f(x) u_\lambda^{-\gamma} \phi -\lambda g(x) u_\lambda^{\alpha} \phi \right]d\mu = 0.
\end{eqnarray}
Thus $u_\lambda$ is a weak solution of (\ref{eq1}).
\par
For any fixed $x_0\in V$, let
\begin{eqnarray}\label{z1}
e_{x_0}=
\begin{cases}
\begin{array}{ll}
1\;\;\;\;\;\;\;x=x_0,\\
0\;\;\;\;\;\;\;x\neq x_0.
\end{array}
\end{cases}
\end{eqnarray}
Replacing $u_k$ in Appendix 3 with $u_\lambda$ for any given $\phi \in W$, we can obtain that
\begin{eqnarray}\label{q3}
\int_V \left[ |\nabla u_\lambda|^{p-2}\Gamma(u_\lambda,\phi)+a(x) u_\lambda^{p-1}\phi
+ |\nabla u_\lambda|^{q-2}\Gamma(u_\lambda,\phi)+b(x) u_\lambda^{q-1} \phi -\lambda g(x) u_\lambda^{\alpha} \phi \right]d\mu <+\infty .
\end{eqnarray}
Substituting (\ref{z1}) into (\ref{b5}) and taking into account (\ref{q3}), we infer that
\begin{eqnarray*}
&&\int_V f(x) u_\lambda^{-\gamma} e_{x_0} d\mu=f(x_0) u_\lambda^{-\gamma}(x_0)\\
&=&\int_V \left[ |\nabla u_\lambda|^{p-2}\Gamma(u_\lambda,e_{x_0})+a(x) u_\lambda^{p-1}e_{x_0}
+ |\nabla u_\lambda|^{q-2}\Gamma(u_\lambda,e_{x_0})+b(x) u_\lambda^{q-1} e_{x_0} -\lambda g(x) u_\lambda^{\alpha} e_{x_0} \right]d\mu <+\infty.
\end{eqnarray*}
Therefore, due to the arbitrariness of $x_0$, we conclude that $u_\lambda>0$, and hence $u_\lambda$ is a strictly positive solution.
\qed

\vskip2mm
\textbf{3.2. Solution in $\mathcal{D}^-_\lambda$ with $\lambda \in (0,\Lambda_{*})$}
 \vskip2mm
\noindent
\par
Note that $\mathcal{D}^-_\lambda$ is closed as $\lambda \in (0,\Lambda_{*})$. Then using Ekeland's variational principle, there exists a minimizing sequence $v_k \in \mathcal{D}^-_\lambda$ such that\\
$(i)$\;\;$J_\lambda(v_k) < \inf_{\mathcal{D}^-_\lambda } J_\lambda +\frac{1}{k}$;\\
$(ii)$\;\;$J_\lambda(\tilde{v})\geq J_\lambda(v_k)-\frac{1}{k}\|\tilde{v}-v_k\|_{W} $ for all $\tilde{v}\in \mathcal{D}^-_\lambda$.\\
By coerciveness of $J_\lambda(u)$ and Lemma 2.1, there exists $v_\lambda \in {W_+}$ such that, up to a subsequence,
\begin{eqnarray}
\label{x37}
 \begin{cases}
v_k \rightharpoonup v_\lambda,& \mbox{in} \;W_+,\;k \rightarrow \infty,\\
v_k(x)\rightarrow v_\lambda(x),& \forall \;x\in V,\;k \rightarrow \infty.
   \end{cases}
\end{eqnarray}
Moreover, by (\ref{x37}), Lemma 2.2 and Lemma 2.5 (ii), we get that $\|v_\lambda\|_{\alpha+1}\geq S(\lambda)$. Then $v_\lambda \not\equiv 0$. Next, we prove that $v_k \rightarrow v_\lambda$ in $W_+$. Suppose that $\|v_k-v_\lambda\|_{W_a^{1,p}(V)}\rightarrow c_1$ and $\|v_k-v_\lambda\|_{W_b^{1,q}(V)}\rightarrow c_2$, where $c_1$ and $c_2$ are not equal to zero at the same time. Applying Lemma 2.2, we deduce that
\begin{eqnarray}
&&(i)\;\;\inf_{\mathcal{D}^-_\lambda}J_\lambda=\lim_{k\rightarrow\infty}J_\lambda(v_k)\nonumber\\
&    =    &
\lim_{k\rightarrow\infty}\left[\frac{1}{p}\|v_k\|^p_{W_a^{1,p}(V)}+\frac{1}{q}\|v_k\|^q_{W_b^{1,q}(V)}-\frac{1}{1-\gamma}\int_V f(x) v_k^{1-\gamma}d\mu - \frac{\lambda}{\alpha+1}\int_V g(x) v_k^{\alpha+1}d\mu
\right]\nonumber\\
&    =    &
\frac{1}{p}(\|v_\lambda\|^p_{W_a^{1,p}(V)}+c_1^p)+\frac{1}{q}(\|v_\lambda\|^q_{W_b^{1,q}(V)}+c_2^q)-\frac{1}{1-\gamma}\int_V f(x) v_\lambda^{1-\gamma}d\mu - \frac{\lambda}{\alpha+1}\int_Vg(x) v_\lambda^{\alpha+1}d\mu\nonumber\\
&    =    &J_\lambda(v_\lambda)+ \frac{c_1^p}{p}+\frac{c_2^q}{q}.\nonumber\\
&&(ii)\;\;v_k\in\mathcal{D}_\lambda \Rightarrow \varphi_{v_k}(1)=0\nonumber\\
&    \Rightarrow   & \|v_k\|^p_{W_a^{1,p}(V)}+\|v_k\|^q_{W_b^{1,q}(V)}-\int_V f(x)v_k^{1-\gamma}d\mu - \int_Vg(x) v_k^{\alpha+1}d\mu=0\nonumber\\
&    \Rightarrow    &
0=\lim_{k\rightarrow\infty}\left[\|v_k\|^p_{W_a^{1,p}(V)}+\|v_k\|^q_{W_b^{1,q}(V)}-\int_V f(x) v_k^{1-\gamma}d\mu - \int_Vg(x) v_k^{\alpha+1}d\mu
\right]\nonumber\\
&=&\|v_\lambda\|^p_{W_a^{1,p}(V)}+c_1^p+\|v_\lambda\|^q_{W_b^{1,q}(V)}+c_2^q-\int_V f(x) v_\lambda^{1-\gamma}d\mu-\int_Vg(x) v_\lambda^{\alpha+1}d\mu\nonumber\\
&=&\varphi_{v_\lambda}(1)+c_1^p+c_2^q.\label{a5}\\
&&(iii)\;\;v_k\in\mathcal{D}^-_\lambda\Rightarrow  \varphi_{v_k}'(1)<0\nonumber\\
&    \Rightarrow   & (p-1+\gamma)\|v_k\|^p_{W_a^{1,p}(V)}+(q-1+\gamma)\|v_k\|^q_{W_b^{1,q}(V)}-\lambda(\alpha+\gamma)\int_V g(x)v_k^{\alpha+1}d\mu<0\nonumber\\
&   \Rightarrow   &
0\geq \lim_{k\rightarrow\infty}\left[ (p-1+\gamma)\|v_k\|^p_{W_a^{1,p}(V)}+(q-1+\gamma)\|v_k\|^q_{W_b^{1,q}(V)}-\lambda(\alpha+\gamma)\int_V g(x)v_k^{\alpha+1}d\mu \right]\nonumber\\
&=& (p-1+\gamma)(\|v_\lambda\|^p_{W_a^{1,p}(V)}+c_1^p)+(q-1+\gamma)(\|v_\lambda\|^q_{W_b^{1,q}(V)}+c_2^q)- \lambda(\alpha+\gamma)\int_Vg(x) v_\lambda^{\alpha+1}d\mu\nonumber\\
&=&\varphi_{v_\lambda}'(1)+c_1^p+c_2^q.\label{a6}
\end{eqnarray}
So $\varphi_{v_\lambda}(1)<0$ and $\varphi_{v_\lambda}'(1)<0$. According to Remark 2.2, there exists a $t_2(v_\lambda)\in (0,1)$ satisfying $\varphi_{v_\lambda}(t_2)=0$ with $\varphi'_{v_\lambda}(t_2)<0$, which implies that $t_2 v_\lambda \in \mathcal{D}^-_\lambda$. Let $j(t)=J_\lambda(tv_\lambda)+\frac{c_1^pt^p}{p}+\frac{c_2^qt^q}{q}$. Then by (\ref{a5}) and (\ref{a6}), we can obtain that $j'(1)=0$ and $j'(t_2)=\frac{d}{dt}J_\lambda(t_2v_\lambda)+c_1^p t_2^{p-1}+c_2^q t_2^{q-1}=t_2^{-\gamma} \varphi_{ v_\lambda}(t_2)+c_1^p t_2^{p-1}+c_2^q t_2^{q-1}=c_1^p t_2^{p-1}+c_2^q t_2^{q-1}>0$ for $t_2\in(0,1)$. These results collectively establish that $j(t)$ is strictly increasing on $[t_2,1]$. Therefore
$$
\inf_{\mathcal{D}^-_\lambda}J_\lambda(u)=\lim_{k\rightarrow\infty}J_\lambda(v_k)=j(1)
>j(t_2)>J_\lambda(t_2v_\lambda)\geq\inf_{\mathcal{D}^-_\lambda}J_\lambda.
$$
It is a contradiction. Thus $c_1=c_2=0$ and then $v_k \rightarrow v_\lambda$ in $W_+$. Since $\mathcal{D}^-_\lambda$ is closed, we obtain $v_\lambda\in \mathcal{D}^-_\lambda$. Furthermore, similar to the argument in Section 3.1, it can be shown that $v_\lambda$ is a strictly positive solution.

\qed

\vskip2mm
\noindent
{\bf Proof of Theorem 1.1}
In Section 3.1 and 3.2, we conclude that $u_\lambda\in \mathcal{D}^+_\lambda$ and $v_\lambda\in \mathcal{D}^-_\lambda$ are the strictly positive solutions.
Furthermore, given $u\in W_+$, we can find $t_1(u)$ and $t_2(u)$ $(t_i>0, i=1,2)$ such that $t_1(u)u\in \mathcal{D}^+_\lambda$ and $t_2(u)u\in \mathcal{D}^-_\lambda$. According to Remark 2.2, we get that $\varphi'_{u}(t_1)>0$ and $\varphi'_{u}(t_2)<0$. From Figure 2 and Figure 3, we conclude $J_\lambda(t_1(u)u)<J_\lambda(t_2(u)u)$. Since $v_\lambda \in \mathcal{D}^-_\lambda$, we know $t_2(v_\lambda)=1$ and there exists a $t_1(v_\lambda)>0$ such that $t_1(v_\lambda)v_\lambda\in \mathcal{D}^+_\lambda$. By setting $u=v_\lambda$, it follows that $J_\lambda(v_\lambda)>J_\lambda(t_1(v_\lambda)v_\lambda)$. Since $J_\lambda(u_\lambda)=\inf_{\mathcal{D}^+_\lambda\cup \mathcal{D}^0_\lambda}J_\lambda$, we obtain $J_\lambda(t_1(v_\lambda)v_\lambda)\geq J_\lambda(u_\lambda)$. Hence $J_\lambda(u_\lambda)$ is the least energy level in $W_+$.
\par
From (\ref{w1}) and (\ref{x11}), we present the further calculation about $v_\lambda \in \mathcal{D}^-_\lambda$:
\begin{eqnarray*}
&&\|v_\lambda\|_{W} >\|v_\lambda\|_{W_a^{1,p}(V)}>X(\lambda) \\
& = &\left( \frac{p-1+\gamma}{\lambda(\alpha+\gamma)\|g\|_\infty\mu_0^{\frac{p-1-\alpha}{p}}a_0^{-\frac{\alpha+1}{p}}}\right)
^{\frac{1}{\alpha+1-p}}\\
& = &\left[ \frac{\Lambda_{*}}{\lambda}\left(\frac{\alpha+\gamma}{\alpha+1-p}\right)^{\frac{\alpha+1-p}{p-1+\gamma}}a_0^{\frac{\alpha+1}{p}} \frac{1}{a_0^{\frac{\alpha+\gamma}{p-1+\gamma}}}\|f\|_{\frac{p}{p-1+\gamma}}^{\frac{\alpha+1-p}{p-1+\gamma}}\right]
^{\frac{1}{\alpha+1-p}}\\
& = &
a_0^{\frac{\gamma-1}{p(p-1+\gamma)}}\|f\|_{\frac{p}{p-1+\gamma}}^{\frac{1}{p-1+\gamma}}
\left(\frac{\alpha+\gamma}{\alpha+1-p}\right)^{\frac{1}{p-1+\gamma}}
\left(\frac{\Lambda_{*}}{\lambda}\right)^{\frac{1}{\alpha+1-p}}.
\end{eqnarray*}
Next, we consider the case that $p\rightarrow \alpha+1$. Let $\alpha+1=p+\varepsilon$. Then
\begin{eqnarray*}
&&\|v_\lambda\|_{W} >\|v_\lambda\|_{W_a^{1,p}(V)} :=X(\lambda) \\
& > &
a_0^{\frac{\gamma-1}{p(p-1+\gamma)}}\|f\|_{\frac{p}{p-1+\gamma}}^{\frac{1}{p-1+\gamma}}
\left(\frac{p+\varepsilon-1+\gamma}{\varepsilon}\right)^{\frac{1}{p-1+\gamma}}
\left(\frac{\Lambda_{*}}{\lambda}\right)^{\frac{1}{\varepsilon}}\\
& = &
a_0^{\frac{\gamma-1}{p(p-1+\gamma)}}\|f\|_{\frac{p}{p-1+\gamma}}^{\frac{1}{p-1+\gamma}}
\left(1+\frac{p-1+\gamma}{\varepsilon}\right)^{\frac{1}{p-1+\gamma}}
\left(\frac{\Lambda_{*}}{\lambda}\right)^{\frac{1}{\varepsilon}}\\
&:=&
C_\varepsilon\left(\frac{\Lambda_{*}}{\lambda}\right)^{\frac{1}{\varepsilon}},
\end{eqnarray*}
where $C_\varepsilon=a_0^{\frac{\gamma-1}{p(p-1+\gamma)}}\|f\|_{\frac{p}{p-1+\gamma}}^{\frac{1}{p-1+\gamma}}
\left(1+\frac{p-1+\gamma}{\varepsilon}\right)^{\frac{1}{p-1+\gamma}}$.
Then we find that $C_\varepsilon\rightarrow +\infty$ as $ \varepsilon\rightarrow 0^+$. The proof is complete.
\qed

\vskip2mm
{\section{The existence and uniqueness of positive solution as $\lambda<0$}}
  \setcounter{equation}{0}
  \par
In this section we shall prove Theorem 1.2.

 \vskip2mm
\noindent
{\bf Lemma 4.1.}  {\it Assume that condition $(A)$ holds and $\lambda<0$. Then the functional $J_\lambda(u,v)$ has the global minimizer in $W$, that is, there exists $ w_\lambda\in W$ such that $J_\lambda(w_\lambda)=\inf_W J_\lambda<0$.
}
\vskip2mm
\noindent
{\bf Proof.}\ \ For $w\in W$, it follows from (\ref{x14}) that if $\lambda<0$ and $0<\gamma<1$, then
\begin{eqnarray*}
J_\lambda(w)
&    =     &\frac{1}{p}\|w\|^p_{W^{1,p}_a(V)}+\frac{1}{q}\|w\|^q_{W^{1,q}_b(V)}-\frac{1}{1-\gamma}\int_V f(x) |w|^{1-\gamma}d\mu - \frac{\lambda}{\alpha+1}\int_Vg(x) |w|^{\alpha+1}d\mu\\
&    \geq     &\frac{1}{p}\|w\|^p_{W^{1,p}_a(V)}+\frac{1}{q}\|w\|^q_{W^{1,q}_b(V)}-\frac{1}{1-\gamma}\cdot a_0^{-\frac{1-\gamma}{p}}\cdot \|f\|_{\frac{p}{p-1+\gamma}}\cdot\|w\|^{1-\gamma}_{W^{1,p}_a(V)}.
\end{eqnarray*}
Since $1-\gamma<1$ and $p,q>1$, $J_\lambda$ is coercive and bounded from below on $W$ for any $\lambda<0$. Hence $m_\lambda=\inf_W J_\lambda$ is well defined. For $\eta>0$ and given $w\in W\setminus\{0\},$ the functional $J_\lambda(\eta w)$ expands as:
\begin{eqnarray*}
J_\lambda(\eta w)
=     \frac{\eta^p}{p}\|w\|^p_{W^{1,p}_a(V)}+\frac{\eta^q}{q}\|w\|^q_{W^{1,q}_b(V)}-\frac{\eta^{1-\gamma}}{1-\gamma}\int_V f(x) |w|^{1-\gamma}d\mu - \lambda\frac{\eta^{\alpha+1}}{\alpha+1}\int_Vg(x) |w|^{\alpha+1}d\mu.
\end{eqnarray*}
Then $J_\lambda(\eta w)<0$ for $\eta>0$ small enough. Thus $m_\lambda= \inf_W J_\lambda<0$ and there exists a minimizing sequence $\{w_k\}\subset W$ such that $\lim_{k\rightarrow \infty} J_\lambda(w_k)=m_\lambda<0.$ The coerciveness of $J_\lambda$ on $W$ shows that $\{w_k\}$ is bounded in $W$. Then there exists $w_\lambda\in W$, up to a subsequence, such that
\begin{eqnarray*}
\label{x42}
 \begin{cases}
w_k \rightharpoonup w_\lambda,& \mbox{in} \;W,\;k \rightarrow \infty,\\
w_k(x)\rightarrow w_\lambda(x),& \forall \;x\in V,\;k \rightarrow \infty.
   \end{cases}
\end{eqnarray*}
By the weakly lower semi-continuity of the norm, and similar to the argument of (\ref{x21}) and (\ref{x22}), it can be deduced that
\begin{eqnarray*}
J_\lambda(w_\lambda)
&    =     & \frac{1}{p}\|w_\lambda\|^p_{W^{1,p}_a(V)}+\frac{1}{q}\|w_\lambda\|^q_{W^{1,q}_b(V)}-\frac{1}{1-\gamma}\int_V f(x) |w_\lambda|^{1-\gamma}d\mu - \frac{\lambda}{\alpha+1}\int_Vg(x) |w_\lambda|^{\alpha+1}d\mu\\
&    \leq     & \liminf_{k\rightarrow \infty}\left[ \frac{1}{p}\|w_k\|^p_{W^{1,p}_a(V)}+\frac{1}{q}\|w_k\|^q_{W^{1,q}_b(V)}-\frac{1}{1-\gamma}\int_V f(x) |w_k|^{1-\gamma}d\mu - \frac{\lambda}{\alpha+1}\int_Vg(x) |w_k|^{\alpha+1}d\mu\right]\\
&    =     & \liminf_{k\rightarrow \infty}J_\lambda(w_k)=m_\lambda.
\end{eqnarray*}
On the other hand, obviously, there holds $m_\lambda \leq J_\lambda(w_\lambda)$. Thus $m_\lambda= J_\lambda(w_\lambda)$.
The proof is complete.
\qed

\vskip2mm
\noindent
{\bf Proof of Theorem 1.2. }
Firstly, we claim that for any $\psi\in W$ with $\psi\geq 0$, it holds $\int_V(|\nabla w_\lambda|^{p-2}\Gamma(w_\lambda,\psi)+a(x)w_\lambda^{p-1} \psi ) d\mu+
\int_V(|\nabla w_\lambda|^{q-2}\Gamma(w_\lambda,\psi)+ b(x) w_\lambda^{q-1}\psi)d\mu
-\int_V f(x) w_\lambda^{-\gamma}\psi d\mu
-\lambda\int_V g(x) w_\lambda^{\alpha}\psi d\mu \geq 0$.
\par
By (\ref{x20}), we get $|\nabla |w_\lambda||^p \leq |\nabla w_\lambda|^p$. Then $J_\lambda(|w_\lambda|)\leq J_\lambda(w_\lambda)$. On the other hand, from Lemma 4.1, we get that $J_\lambda(|w_\lambda|)\geq J_\lambda(w_\lambda)$. Thus $J_\lambda(|w_\lambda|)= J_\lambda(w_\lambda)$. Hence, we may assume $w_\lambda \geq 0$ with $w_\lambda\not\equiv0$. Given that $J_\lambda(w_\lambda)= \inf_W J_\lambda$, it follows that for all non-negative $\psi\in W$ and $0<\eta<1$,
\begin{eqnarray}
\label{b2}
0&\leq&
J_\lambda(w_\lambda+\eta\psi)-J_\lambda(w_\lambda)\nonumber\\
&    =     & \frac{1}{p}\left[\|w_\lambda+\eta\psi\|^p_{W^{1,p}_a(V)}-\|w_\lambda\|^p_{W^{1,p}_a(V)}\right]
+\frac{1}{q}\left[\|w_\lambda+\eta\psi\|^q_{W^{1,q}_b(V)}-\|w_\lambda\|^q_{W^{1,q}_b(V)}\right]\nonumber\\
&&-\frac{1}{1-\gamma}\int_V f(x) [(w_\lambda+\eta\psi)^{1-\gamma}-w_\lambda^{1-\gamma}]d\mu - \frac{\lambda}{\alpha+1}\int_Vg(x) [(w_\lambda+\eta\psi)^{\alpha+1}-w_\lambda^{\alpha+1}]d\mu.
\end{eqnarray}
Similar to (\ref{x34}), we can also obtain that $\frac{1}{1-\gamma}\liminf_{\eta \rightarrow 0^+}\int_V \frac{f(x) [(w_\lambda+\eta\psi)^{1-\gamma}-w_\lambda^{1-\gamma}]}{\eta}d\mu \geq \int_Vf(x)w_\lambda^{-\gamma}\psi d\mu$.
Besides, there exists a $\theta'\in (0,1)$, such that
\begin{eqnarray*}
&&\frac{g(x)[(w_\lambda+\eta\psi)^{\alpha+1}-w_\lambda^{\alpha+1}]}{\eta}\\
&     \leq      &
\|g\|_\infty\cdot\frac{(w_\lambda+\eta\psi)^{\alpha+1}-w_\lambda^{\alpha+1}}{\eta}
=\|g\|_\infty\cdot(\alpha+1)(w_\lambda+\theta'\eta\psi)^\alpha\psi\\
&     \leq      &
(\alpha+1)\|g\|_\infty(w_\lambda+\psi)^\alpha \psi \leq (\alpha+1)\|g\|_\infty2^{\alpha-1}(w_\lambda^\alpha +\psi^\alpha) \psi\\
 &=&  (\alpha+1)\|g\|_\infty2^{\alpha-1}(w_\lambda^\alpha\psi+\psi^{\alpha+1}).
\end{eqnarray*}
Since
$\int_V w_\lambda^\alpha\psi d\mu
\leq \left(\int_V |w_\lambda|^{\alpha+1}d\mu \right)^{\frac{\alpha}{\alpha+1}}\cdot \left (\int_V |\psi|^{\alpha+1}d\mu \right)^{\frac{1}{\alpha+1}} < \infty$ and $\int_V |\psi|^{\alpha+1}d\mu <\infty$, we obtain $(\alpha+1)\|g\|_\infty2^{\alpha-1}(w_\lambda^\alpha\psi+\psi^{\alpha+1}) \in L^1(V)$. Thus, the Dominated convergence theorem implies that
\begin{eqnarray*}
&&\lim_{\eta \rightarrow 0^+}\int_V\frac{g(x)[(w_\lambda+\eta\psi)^{\alpha+1}-w_\lambda^{\alpha+1}]}{\eta}d\mu=\int_V\lim_{\eta \rightarrow 0^+}\frac{g(x)[(w_\lambda+\eta\psi)^{\alpha+1}-w_\lambda^{\alpha+1}]}{\eta}d\mu\\
&  =  &
\int_V\lim_{\eta \rightarrow 0^+} g(x)(\alpha+1)(w_\lambda+\theta'\eta\psi)^{\alpha}\psi d\mu =(\alpha+1)\int_V g(x) w_\lambda^{\alpha}\psi d\mu.
\end{eqnarray*}
Dividing (\ref{b2}) by $\eta>0$ and passing to the liminf as $\eta\rightarrow0^+$, we deduce
\begin{eqnarray}
\label{b3}
0
& \leq  &
\liminf_{\eta \rightarrow 0^+}\left[\frac{1}{p}\frac{\|w_\lambda+\eta\psi\|^p_{W^{1,p}_a(V)}-\|w_\lambda\|^p_{W^{1,p}_a(V)}}{\eta}
+\frac{1}{q}\frac{\|w_\lambda+\eta\psi\|^q_{W^{1,q}_b(V)}-\|w_\lambda\|^q_{W^{1,q}_b(V)}}{\eta}\right.\nonumber\\
&&\left.-\frac{1}{1-\gamma}\int_V \frac{f(x) [(w_\lambda+\eta\psi)^{1-\gamma}-w_\lambda^{1-\gamma}]}{\eta}d\mu
- \frac{\lambda}{\alpha+1}\int_Vg(x) \frac{(w_\lambda+\eta\psi)^{\alpha+1}-w_\lambda^{\alpha+1}}{\eta}d\mu\right]\nonumber\\
& \leq  &
\int_V(|\nabla w_\lambda|^{p-2}\Gamma(w_\lambda,\psi)+a(x)w_\lambda^{p-1} \psi ) d\mu+
\int_V(|\nabla w_\lambda|^{q-2}\Gamma(w_\lambda,\psi)+ b(x)w_\lambda^{q-1}\psi)d\mu\nonumber\\
&&-\int_V f(x) w_\lambda^{-\gamma}\psi d\mu
-\lambda\int_V g(x) w_\lambda^{\alpha}\psi d\mu .
\end{eqnarray}
\par
Next, we prove $w_\lambda>0$ is a solution of (\ref{eq1}). Let $h:[-\sigma, \sigma]\mapsto \R$ and $h(\eta):=J_\lambda(w_\lambda+\eta w_\lambda)$ for a given $\sigma >0$. Then from Lemma 4.1, we know that $h(\eta)$ attains its minimum at $\eta=0$. Thus it follows that
\begin{eqnarray}
\label{b4}
h'(0)=\|w_\lambda\|^p_{W^{1,p}_a(V)}+\|w_\lambda\|^q_{W^{1,q}_b(V)}-\int_V f(x)w_\lambda^{1-\gamma}d\mu- \lambda\int_Vg(x) w_\lambda^{\alpha+1}d\mu=0.
\end{eqnarray}
Similar to the proof of (\ref{b5}), we also suppose $\phi \in W$ and $\varepsilon>0$, and define $\Phi\in W$ by $\Phi= (w_\lambda+ \varepsilon \phi)^+$. Inserting $\Phi$ into (\ref{b3}) and combining (\ref{b4}), we deduce that
\begin{eqnarray*}\label{a7}
\int_V \left[ |\nabla w_\lambda|^{p-2}\Gamma(w_\lambda,\phi)+a(x) w_\lambda^{p-1}\phi
+ |\nabla w_\lambda|^{q-2}\Gamma(w_\lambda,\phi)+b(x) w_\lambda^{q-1} \phi  - f(x) w_\lambda^{-\gamma} \phi -\lambda g(x) w_\lambda^{\alpha} \phi \right]d\mu \geq 0.
\end{eqnarray*}
The above inequality also holds for $-\phi$. Then it is concluded that
\begin{eqnarray}\label{b7}
\int_V \left[ |\nabla w_\lambda|^{p-2}\Gamma(w_\lambda,\phi)+a(x) w_\lambda^{p-1}\phi
+ |\nabla w_\lambda|^{q-2}\Gamma(w_\lambda,\phi)+b(x) w_\lambda^{q-1} \phi  - f(x) w_\lambda^{-\gamma} \phi -\lambda g(x) w_\lambda^{\alpha} \phi \right]d\mu = 0.
\end{eqnarray}
Thus $w_\lambda$ is a solution of (\ref{eq1}). Furthermore, similar to the final argument in Section 3.1, we can prove that $w_\lambda$ is a positive solution.
\par
Finally, we claim the solution $w_\lambda$ is unique. We first introduce an elementary inequality (see \cite{Peral I 1997}):
\begin{eqnarray}
\label{s7}
\langle|b|^{p-2}b-|a|^{p-2}a,b-a\rangle \geq
\begin{cases}
c_p|b-a|^p, \;\;\;\;\;\;\;\;\;\;\;\;\;\;\;\;\;\;\;\;\;\;\;\;\;\;\;\;\;\;\; p\geq2,\\
c_p|b-a|^2(|b|+|a|)^{p-2}, \;\;\;\;\;\;\;\;\;\ 1<p<2,
 \end{cases}
\end{eqnarray}
where $a$ and $b$ are two vectors in $\R^N(N\geq 3)$ with the standard inner product $\langle\cdot,\cdot\rangle$ and $c_p$ is a positive constant depending only on $p$.
\par
Suppose $w_*$ is also a positive solution. Then for any $\phi\in W$, there exists
\begin{eqnarray}\label{b6}
\int_V \left[ |\nabla w_*|^{p-2}\Gamma(w_*,\phi)+a(x) w_*^{p-1}\phi
+ |\nabla w_*|^{q-2}\Gamma(w_*,\phi)+b(x) w_*^{q-1} \phi  - f(x) w_*^{-\gamma} \phi -\lambda g(x) w_*^{\alpha} \phi \right]d\mu = 0.
\end{eqnarray}
Subtracting (\ref{b7}) and (\ref{b6}) and choosing $\phi=w_\lambda-w_*$, it follows from $\lambda<0$ that
\begin{eqnarray*}
0  &\geq&   \int_V f(x)(w_\lambda^{-\gamma}-w_*^{-\gamma})(w_\lambda-w_*)d\mu\\
&  =   &  \int_V \left[\left(|\nabla w_\lambda|^{p-2}\Gamma(w_\lambda,w_\lambda-w_*)-|\nabla w_*|^{p-2}\Gamma(w_*,w_\lambda-w_*)\right)
+a(x)\left( w_\lambda^{p-1}(w_\lambda-w_*)- w_*^{p-1} (w_\lambda-w_*) \right)\right.\\
&&\left.+   \left(|\nabla w_\lambda|^{q-2}\Gamma(w_\lambda,w_\lambda-w_*)-|\nabla w_*|^{q-2}\Gamma(w_*,w_\lambda-w_*)\right)
+b(x)\left( w_\lambda^{q-1}(w_\lambda-w_*)- w_*^{q-1} (w_\lambda-w_*) \right)\right.\\
&&\left.-\lambda g(x)\left(w_\lambda^{\alpha} (w_\lambda-w_*)- w_*^{\alpha} (w_\lambda-w_*)\right)\right]d\mu\\
&  =   &  \int_V \left[(|\nabla w_\lambda|^{p-2}\nabla w_\lambda- |\nabla w_*|^{p-2} \nabla w_* ) \nabla (w_\lambda-w_*)
+a(x)(w_\lambda^{p-1}-w_*^{p-1})(w_\lambda-w_*)\right.\\
&&\left.+(|\nabla w_\lambda|^{q-2}\nabla w_\lambda- |\nabla w_*|^{q-2} \nabla w_* ) \nabla (w_\lambda-w_*)
+b(x)(w_\lambda^{q-1}-w_*^{q-1})(w_\lambda-w_*)\right.\\
&&\left.-\lambda g(x)(w_\lambda^{\alpha}-w_*^{\alpha})(w_\lambda-w_*)\right]d\mu\;\;\;\;\;(\lambda<0)\\
&\geq&
\left\{
\begin{aligned}
\int_V &\left[c_p(|\nabla w_\lambda-\nabla w_*|^p+ a(x) | w_\lambda-w_*|^p+ |\nabla w_\lambda-\nabla w_*|^q+ b(x) | w_\lambda-w_*|^q)- c_{\alpha+1}\lambda g(x)| w_\lambda-w_*|^{\alpha+1}\right ]d\mu,\\
&\;\;\;\;\;\;\;\;\;\;\;\;\;\;\;\;\;\;\;\;\;\;\;\;\;\;\;\;\;\;\;\;\;\;\;\;\;\;\;\;\;\;\;\;\;\;\;\;\;\;\;\;\;\;\;\;\;\;\;\;\;\;\;\;\;\;\;\;\;\;\;\;\;\;\;\;\;\;\;
\;\;\;\;\;\;\;\;\;\;\;\;\;\;\;\;\;\;\;\;\;\;\;\;\;\;\;\;\;\;\;\;\;\;\;\;\;\;\;\;\;\; p\geq2\\
\int_V &\left[c_p(|\nabla w_\lambda-\nabla w_*|^2(|\nabla w_\lambda|+|\nabla w_*|)^{p-2} +a(x)| w_\lambda- w_*|^2( w_\lambda+ w_*)^{p-2}\right.\\
&\left.+|\nabla w_\lambda-\nabla w_*|^2(|\nabla w_\lambda|+|\nabla w_*|)^{q-2} +b(x)| w_\lambda- w_*|^2( w_\lambda+ w_*)^{q-2})\right.\\
&\left.-c_{\alpha+1}\lambda g(x)| w_\lambda-w_*|^{2}( w_\lambda+ w_*)^{\alpha-1}\right.]d\mu,  \;\;\;\;\;\;\;\;\;\;\;\;\;\;\;\;\;\;\;\;\;\;\;\;\;\;\;\;\;\;\;\;\;\;\;\;\;\;\;\;\;\;\;\;\;\;\;\;\;\;\;\;\;\;  1<p<2
\end{aligned}
\right.
\\
& \geq & 0.
\end{eqnarray*}
Thus $w_\lambda =w_*$ and so $w_\lambda$ is the unique solution of (\ref{eq1}). The proof is complete.
\qed

\vskip2mm

{\section{Appendix A}}
  \setcounter{equation}{0}
\noindent
{\bf Appendix A.1.}  {\it If $\lambda=\Lambda_{*}$, then $X(\lambda)=X_0$.}

 \vskip0mm
 \noindent
{\bf Proof.}\ \ By (\ref{x11}), we infer that
\begin{eqnarray*}
&&X(\Lambda_{*})\\
&  =   & \left( \frac{p-1+\gamma}{\Lambda_{*}(\alpha+\gamma)\|g\|_\infty\mu_0^{\frac{p-1-\alpha}{p}}a_0^{-\frac{\alpha+1}{p}}}\right)
^{\frac{1}{\alpha+1-p}}\\
&  =    &\left( \frac{p-1+\gamma}{\left(\frac{p+\gamma-1}{\alpha+\gamma}\right)\cdot\left(\frac{\alpha+1-p}{\alpha+\gamma}\right)^{\frac{\alpha+1-p}{p+\gamma-1}}\frac{1}{\|g\|_\infty}\cdot
\mu_0^{\frac{\alpha+1-p}{p}}a_0^{\frac{\alpha+\gamma}{p-1+\gamma}} \left(\frac{1}{\|f\|_{\frac{p}{p-1+\gamma}}}\right)^{\frac{\alpha+1-p}{p+\gamma-1}}(\alpha+\gamma)\|g\|_\infty\mu_0^{\frac{p-1-\alpha}{p}}a_0^{-\frac{\alpha+1}{p}}}\right)
^{\frac{1}{\alpha+1-p}}\\
&  =    & \left( \frac{\alpha+\gamma}{\alpha+1-p} \right)
^{\frac{1}{p+\gamma-1}}\cdot \|f\|^{\frac{1}{p+\gamma-1}}_{\frac{p}{p-1+\gamma}}\cdot a_0^{\frac{\gamma-1}{p(p+\gamma-1)}}=X_0,
\end{eqnarray*}
where
\begin{eqnarray*}
\left( \frac{p-1+\gamma}{\left(\frac{p+\gamma-1}{\alpha+\gamma}\right)\cdot\left(\frac{\alpha+1-p}{\alpha+\gamma}\right)^{\frac{\alpha+1-p}{p+\gamma-1}}
(\alpha+\gamma)}\right)
^{\frac{1}{\alpha+1-p}}=\left( \frac{\alpha+\gamma}{\alpha+1-p} \right)
^{\frac{1}{p+\gamma-1}},
\end{eqnarray*}
and
\begin{eqnarray*}
&&\left({-\frac{\alpha+\gamma}{p-1+\gamma}}+{\frac{\alpha+1}{p}}\right)\cdot {\frac{1}{\alpha+1-p}} = {-\frac{\alpha+\gamma}{(p-1+\gamma)(\alpha+1-p)}+\frac{\alpha+1}{p(\alpha+1-p)}}\\
&=&{\frac{-\alpha p-\gamma p+\alpha p-\alpha +\alpha \gamma -1 +p+\gamma}{p(p-1+\gamma)(\alpha+1-p)}}
={\frac{(\gamma-1)(\alpha+1-p)}{p(p-1+\gamma)(\alpha+1-p)}}
={\frac{\gamma-1}{p(p-1+\gamma)}}
\end{eqnarray*}
Thus, the proof is complete.
\qed

 \vskip2mm
\noindent
{\bf Appendix A.2.} {\it If $\lambda=\Lambda_{*}$, then $S(\lambda)=S_0$.}

 \vskip0mm
 \noindent
{\bf Proof.}\ \ By (\ref{x15}), we infer that
\begin{eqnarray*}
&&S(\Lambda_{*})\\
&  =    &\left( \frac{p-1+\gamma}{\lambda(\alpha+\gamma)} \right)^{\frac{1}{\alpha+1-p}}
\left(\frac{1}{\mu_0}\right)^{-\frac{1}{\alpha+1}} a_0^{\frac{1}{\alpha+1-p}} \left( \frac{1}{\|g\|_\infty} \right)^{\frac{1}{\alpha+1-p}}\\
&  =   & \left( \frac{(p-1+\gamma)\|g\|_\infty}{\left(\frac{p+\gamma-1}{\alpha+\gamma}\right)\left(\frac{\alpha+1-p}{\alpha+\gamma}\right)^{\frac{\alpha+1-p}{p+\gamma-1}}
\mu_0^{\frac{\alpha+1-p}{p}}a_0^{\frac{\alpha+\gamma}{p-1+\gamma}} \left(\frac{1}{\|f\|_{\frac{p}{p-1+\gamma}}}\right)^{\frac{\alpha+1-p}{p+\gamma-1}}(\alpha+\gamma)} \right)^{\frac{1}{\alpha+1-p}}
\left(\frac{1}{\mu_0}\right)^{-\frac{1}{\alpha+1}} a_0^{\frac{1}{\alpha+1-p}} \left( \frac{1}{\|g\|_\infty} \right)^{\frac{1}{\alpha+1-p}}\\
&  =   & \left(\frac{\alpha+\gamma}{\alpha+1-p}\right)^{\frac{1}{p+\gamma-1}}
\left(\frac{1}{\mu_0}\right)^{\frac{1}{p}-\frac{1}{\alpha+1}}
\left( \frac{1}{a_0}\right)^{\frac{\alpha+\gamma}{(p-1+\gamma)(\alpha+1-p)}-\frac{1}{\alpha+1-p}} \|f\|_{\frac{p}{p-1+\gamma}}^{\frac{1}{p+\gamma-1}}\\
&  =   & \left(\frac{\alpha+\gamma}{\alpha+1-p}\right)^{\frac{1}{p+\gamma-1}}
\left(\frac{1}{\mu_0}\right)^{\frac{\alpha+1-p}{p(\alpha+1)}}
\left( \frac{1}{a_0}\right)^{\frac{1}{p-1+\gamma}} \|f\|_{\frac{p}{p-1+\gamma}}^{\frac{1}{p+\gamma-1}}
\end{eqnarray*}
Thus, the proof is complete.
\qed

 \vskip2mm
\noindent
{\bf Appendix A.3.} {\it If $\{u_k\} \subset W_+$ is bounded, then for any $\varphi \in W_+$, there hold
\begin{eqnarray*}
&&I:= \int_V (|\nabla u_k|^{p-2}\Gamma(u_k , \varphi)+a(x)u_k^{p-1}\varphi)d\mu\leq \|u_k\|^{p-1}_{W_a^{1,p}(V)}\|\varphi\|_{W_a^{1,p}(V)},\\
&&II:= \int_V (|\nabla u_k|^{q-2}\Gamma(u_k, \varphi)+b(x)u_k ^{q-1} \varphi)d\mu\leq\|u_k\|^{q-1}_{W_b^{1,q}(V)}\|\varphi\|_{W_b^{1,q}(V)} ,\\
&&III:=\int_V g(x)u_k^{\alpha} \varphi d\mu\leq\|g\|_\infty \left(\frac{C_2}{ \mu_0^{\frac{1}{p}} a_0^{\frac{1}{p}} +\mu_0^{\frac{1}{q}}b_0^{\frac{1}{q}}}\right)^{\alpha-p+1}\|u_k\|_p^{p-1}\|\varphi\|_p.
\end{eqnarray*}
}
 \vskip0mm
 \noindent
{\bf Proof.}\ \
For all $0<\beta<1$ and $c,d,e,f \geq 0$, by H\"{o}lder's inequality, it holds that
\begin{eqnarray*}
f^\beta c^{1-\beta}+d^\beta e^{1-\beta} \leq (f+d)^\beta(c+e)^{1-\beta}.
\end{eqnarray*}
Let $\beta=\frac{p-1}{p}$. Then from H\"{o}lder's inequality, $\varphi \geq 0$ and $u_k \geq 0$, we deduce that
\begin{eqnarray*}
I &=& \int_V (|\nabla u_k|^{p-2}\Gamma(u_k , \varphi)+a(x)u_k ^{p-1} \varphi)d\mu\\
&  \leq  &
\int_V (|\nabla u_k|^{p-2} |\nabla u_k|\cdot|\nabla\varphi| + a^{\frac{p-1}{p}}u_k^{p-1} a^{\frac{1}{p}}\varphi)d\mu\\
&  \leq  &
\left(\int_V |\nabla u_k|^{(p-1)\cdot \frac{p}{p-1}}d\mu\right)^{\frac{p-1}{p}} \left( \int_V |\nabla\varphi|^pd\mu\right)^{\frac{1}{p}}
+ \left(\int_V  a^{\frac{p-1}{p}\cdot\frac{p}{p-1}}|u_k |^{(p-1)\cdot\frac{p}{p-1}}d\mu\right)^{\frac{p-1}{p}} \left(\int_V a^{\frac{1}{p}\cdot p}|\varphi|^p d\mu\right)^{\frac{1}{p}}\\
&  \leq  &
\|u_k\|^{p-1}_{W_a^{1,p}(V)}\|\varphi\|_{W_a^{1,p}(V)}.
\end{eqnarray*}
Similarly, we also obtain
\begin{eqnarray*}
II &=& \int_V (|\nabla u_k|^{q-2}\Gamma(u_k, \varphi)+b(x)u_k^{q-1} \varphi)d\mu
 \leq
\|u_k\|^{q-1}_{W_b^{1,q}(V)}\|\varphi\|_{W_b^{1,q}(V)}.
\end{eqnarray*}
Similar to (\ref{x26}), we have $\|u_k\|_\infty \leq \frac{C_2}{ \mu_0^{\frac{1}{p}} a_0^{\frac{1}{p}} +\mu_0^{\frac{1}{q}} b_0^{\frac{1}{q}}}$. Then by (\ref{x13}) we obtain
\begin{eqnarray*}
III &=& \int_V g(x)u_k^{\alpha} \varphi d\mu
 \leq
\|g\|_\infty\int_V u_k^{\alpha} \varphi d\mu\\
& \leq &
\|g\|_\infty \|u_k\|^{\alpha-p+1}_\infty  \int_V u_k^{p-1} \varphi d\mu\\
&  \leq &
\|g\|_\infty \|u_k\|^{\alpha-p+1}_\infty  \left(\int_V  |u_k |^{(p-1)\cdot\frac{p}{p-1}}d\mu \right)^{\frac{p-1}{p}} \left(\int_V |\varphi|^p d\mu\right)^{\frac{1}{p}}\\
& \leq &
\|g\|_\infty \left(\frac{C_2}{ \mu_0^{\frac{1}{p}} a_0^{\frac{1}{p}} +\mu_0^{\frac{1}{q}}b_0^{\frac{1}{q}}}\right)^{\alpha-p+1}\|u_k\|_p^{p-1}\|\varphi\|_p\\
& \leq &
\|g\|_\infty \left(\frac{C_2}{ \mu_0^{\frac{1}{p}} a_0^{\frac{1}{p}} +\mu_0^{\frac{1}{q}}b_0^{\frac{1}{q}}}\right)^{\alpha-p+1}a_0^{-\frac{p-1}{p}}\|u_k\|_{W_a^{1,p}(V)}^{p-1}\|\varphi\|_p
\end{eqnarray*}
The proof is complete.
\qed

 \vskip2mm
\noindent
{\bf Appendix A.4.} {\it If $u_k\geq 0$ and $u_k\rightarrow u_\lambda$ in $W_+$, then for any $\varphi \in W_+$,
\begin{eqnarray*}
&&\lim_{k \rightarrow \infty}\left[\int_V(|\nabla u_k|^{p-2}\Gamma(u_k,\varphi)+a(x) u_k ^{p-1} \varphi)d\mu\right.\nonumber\\
&&\left.+\int_V (|\nabla u_k|^{q-2}\Gamma(u_k,\varphi)+b(x) u_k ^{q-1} \varphi)d\mu
-\lambda\int_V g(x) u_k ^{\alpha} \varphi d\mu\right.\\
&&\left.-\left(\int_V(|\nabla u_\lambda|^{p-2}\Gamma(u_\lambda,\varphi)+a(x) u_\lambda ^{p-1} \varphi)d\mu\right.\right.\nonumber\\
&&\left.\left.+\int_V (|\nabla u_\lambda|^{q-2}\Gamma(u_\lambda,\varphi)+b(x) u_\lambda ^{q-1}\varphi)d\mu
-\lambda\int_V g(x) u_\lambda ^{\alpha} \varphi d\mu\right)
 \right] \leq 0.
\end{eqnarray*}
}
 \vskip0mm
 \noindent
{\bf Proof.}\ \ Similar to $(A.8), (A.10)$ and $(A.11)$ in \cite{Yang P 2023}, we can obtain that
\begin{eqnarray*}
&&\lim_{k \rightarrow \infty}\left[\left(\int_V(|\nabla u_k|^{p-2}\Gamma(u_k,\varphi)+a(x)u_k^{p-1}\varphi)d\mu-\int_V(|\nabla u_\lambda|^{p-2}\Gamma(u_\lambda,\varphi)+a(x)u_\lambda^{p-1} \varphi)d\mu\right)\right.\nonumber\\
&&\left.+\left(\int_V (|\nabla u_k|^{q-2}\Gamma(u_k,\varphi)+b(x)u_k^{q-1} \varphi)d\mu-\int_V (|\nabla u_\lambda|^{q-2}\Gamma(u_\lambda,\varphi)+b(x)u_\lambda^{q-1} \varphi)d\mu\right)
 \right]\leq 0
\end{eqnarray*}
We just need to prove $\lim_{k \rightarrow \infty} \left(\lambda\int_V g(x) u_k^{\alpha} \varphi d\mu- \lambda\int_V g(x)u_\lambda^{\alpha}\varphi d\mu\right) \leq 0$. Then
\begin{eqnarray*}
&&\lim_{k \rightarrow \infty}\lambda\int_V g(x)\left(u_k^{\alpha}- u_\lambda^{\alpha}\right) \varphi d\mu\\
&  \leq  &
\lambda \|g\|_\infty \lim_{k \rightarrow \infty}\left(\int_V |u_k^{\alpha}- u_\lambda^{\alpha}|^{\frac{\alpha+1}{\alpha}}d\mu \right)^{\frac{\alpha}{\alpha+1}} \left(\int_V |\varphi|^{\alpha+1}d\mu \right)^{\frac{1}{\alpha+1}}\\
&  \leq  &
\lambda \|g\|_\infty \lim_{k \rightarrow \infty}
\left(\int_V \alpha |u_k- u_\lambda| (u_k^{\alpha-1}+ u_\lambda^{\alpha-1})d\mu \right)^{\frac{\alpha}{\alpha+1}} \|\varphi\|_{\alpha+1}\;\;\;\;\mbox{(Mazur inequality)}\\
&  \leq  &
\lambda \|g\|_\infty\alpha \lim_{k \rightarrow \infty}
\left(\int_V |u_k- u_\lambda|^\alpha d\mu\right)^{\frac{1}{\alpha}\cdot\frac{\alpha}{\alpha+1}} \left(\int_V |u_k + u_\lambda|^{\alpha-1\cdot \frac{\alpha}{\alpha-1}}d\mu \right)^{\frac{\alpha-1}{\alpha}\cdot\frac{\alpha}{\alpha+1}} \|\varphi\|_{\alpha+1}\\
&  =  &
\lambda \|g\|_\infty \alpha \lim_{k \rightarrow \infty}
\|u_k- u_\lambda\|_\alpha^{\frac{\alpha}{\alpha+1}} \left(\int_V |u_k + u_\lambda|^{\alpha}d\mu \right)^{\frac{\alpha-1}{\alpha+1}} \|\varphi\|_{\alpha+1}\\
&  =  &
0.
\end{eqnarray*}
The proof is complete.
\qed

 \vskip2mm
 \noindent
 {\bf Funding}\\
This project is supported by Yunnan Fundamental Research Projects (grant No: 202301AT070465) and  Xingdian Talent Support Program for Young Talents of Yunnan Province.

\vskip2mm
 \noindent
 {\bf Authors' contributions}\\
The authors contribute the manuscript equally.

 \vskip2mm
 \noindent
 {\bf Competing interests}\\
The authors declare that they have no competing interests.

\vskip2mm
\renewcommand\refname{References}
{}

\begin{thebibliography}{}
\bibitem{Ambrosio 2021}Ambrosio, V, and Repov\v{s}, D. Multiplicity and concentration results for a $(p, q)$-Laplacian problem in $\R^N $. Zeitschrift f¨¹r angewandte Mathematik und Physik, 72(1)(2021): 33.

\bibitem{Aris 1994}Aris, R. Mathematical Modelling Techniques. Research Notes in Mathematics 24. Boston, MA-London: Pitman Advanced Publishing Program, (1979).

\bibitem{Alkama S2014}Alkama, S., Desquesnes, X. Infinity Laplacian on graphs with gradient term for image and data clustering. Pattern Recognition Letters, 41(2014): 65-72.

\bibitem{Arnaboldi V 2015}Arnaboldi, V., Passarella, A. Online social networks: human cognitive constraints in Facebook and Twitter personal graphs. Elsevier Science Publishers B. V.(2015).

\bibitem{Alves  2022}Alves, R.L., Santos, C.A., Silva K. Multiplicity of negative-energy solutions for singular-superlinear Schr\"{o}dinger equations with indefinite-sign potential. Communications in Contemporary Mathematics, 24(10)(2022): 2150042.

\bibitem{Brezis 1983}Br\'{e}zis, H., Lieb, E. A relation between pointwise convergence of functions and convergence of functionals. Proceedings of the American Mathematical Society, 88(3)(1983): 486-490.

\bibitem{Baldelli 2022}Baldelli, L., Brizi, Y., Filippucci R. On symmetric solutions for $(p, q)$-Laplacian equations in $\R^N$ with critical terms. The Journal of Geometric Analysis, 32(4)(2022): 120.

\bibitem{Cherfils 2005}Cherfils, L., Il¡¯Yasov, Y. On the stationary solutions of generalized reaction diffusion equations with $p$\&$q$-Laplacian. Communications on Pure and Applied Analysis, 4(1)(2005): 9-22.

\bibitem{Chang logarithmic 2023}Chang, X., Wang R., and Yan D. Ground states for logarithmic Schr\"{o}dinger equations on locally finite graphs. The Journal of Geometric Analysis, 33(7) (2023): 211.

\bibitem{Corra 2023}Corr\^{e}a, F.J.S., dos Santos, G.C., Tavares, L.S., Muhassua, S.S. Existence of solution for a singular elliptic system with convection terms. Nonlinear Analysis: Real World Applications, 66 (2022): 103549.

\bibitem{Chung 2005}Chung, S.Y., Berenstein, C.A. Berenstein. $\omega$-harmonic functions and inverse conductivity problems on networks. SIAM Journal on Applied Mathematics, 65(4)(2005): 1200-1226.

\bibitem{Elmoataz A2015}Elmoataz, A., Toutain, M., Tenbrinck, D. On the $p$-Laplacian and $\infty$-Laplacian on Graphs with Applications in Image and Data Processing. Siam Journal on Imaging Sciences, 8(4)(2015): 2412-2451.

\bibitem{Figueiredo 2011}Figueiredo, G.M. Existence of positive solutions for a class of $p\&q$ elliptic problems with critical growth on $\R^N$. Journal of mathematical analysis and applications, 378(2)(2011): 507-518.

\bibitem{Fife 2013}Fife, P.C. Mathematical aspects of reacting and diffusing systems. Springer-Verlag, Berlin-Heidelberg-New York, (1979).

\bibitem{Fulks 1960}Fulks, W. and Maybee, J.S. A singular non-linear equation. Osaka Journal of Mathematics, 12(1)(1960): 1-19.

\bibitem{Grigor'yan 2017}Grigor'yan, A., Lin, Y. Existence of positive solutions to some nonlinear equations on locally finite graphs. Science China Mathematics, 60(7)(2017): 1311-1324.

\bibitem{Grigor'yan A 2016 yamabe}Grigor'yan A., Lin Y., Yang Y. Yamabe type equations on graphs. Journal of Differential Equations, 261(9)(2016): 4924-4943.

\bibitem{Grigor'yan A 2018 book}Grigor¡¯yan A. Introduction to analysis on graphs. American Mathematical Society, Providence, Rhode Island, (2018).

\bibitem{k-w 2016}Grigor'yan, Alexander, Lin Y., Yang Y. Kazdan-Warner equation on graph. Calculus of Variations and Partial Differential Equations, 55(4)(2016): 1-13.

\bibitem{Goncalves 2007}Goncalves, J.V., Melo, A.L. On existence of $L^\infty$-ground states for singular elliptic equations in the presence of a strongly nonlinear term. Advanced Nonlinear Studies, 7(3)(2007): 475-490.

\bibitem{Jose Valdo 2007}Goncalves, J.V., Santos, C.A. Singular elliptic problems: Existence, non-existence and boundary behavior. Nonlinear Analysis: Theory, Methods \& Applications, 66(9)(2007): 2078-2090.

\bibitem{Han 2020}Han, X., Shao M. Existence and convergence of solutions for nonlinear biharmonic equations on graphs. Journal of Differential Equations, 268(7)(2020): 3936-3961.

\bibitem{Han X L 2021}Han, X., Shao M. $p$-Laplacian Equations on Locally Finite Graphs. Acta Mathematica Sinica, English Series, 37(11)(2021): 1645-1678.

\bibitem{Herrera 1992}Herrera, J.. Julio E. Reaction-diffusion equations in one dimension: particular solutions and relaxation. Physica D: Nonlinear Phenomena, 57(3-4)(1992): 249-266.

\bibitem{Hu 2024}Hu, Y., Wang, M. Life span of solutions to a semilinear parabolic equation on locally finite graphs. arxiv: 2405.18173.


\bibitem{Li 2009}Li G, L X. The existence of nontrivial solutions to nonlinear elliptic equation of $(p,q)$-Laplacian type on $\R^N$. Nonlinear Analysis: Theory, Methods \& Applications 71(5-6) (2009): 2316-2334.

\bibitem{Liu 2009}Liu, X., Guo, Y. Solutions for singular $p$-Laplacian equation in $\R^N$. Journal of Systems Science and Complexity, 22(4)(2009), 597-613.

\bibitem{Ou 2024}Ou, X., Zhang, X. Least energy sign-changing solutions for Kirchhoff-type equations with logarithmic nonlinearity on locally finite graphs. TWMS Journal of Pure and Applied Mathematics, 15(2024): 286-317.

\bibitem{Pan 2023}Pan, G., Ji, C. Existence and convergence of the least energy sign-changing solutions for nonlinear Kirchhoff equations on locally finite graphs. Asymptotic Analysis, 133(4)(2023): 463-482.

\bibitem{Pang 2024}Pang, Y., Xie, J., Zhang, X. Infinitely many solutions for three quasilinear Laplacian systems on weighted graphs. Boundary Value Problems, 2024(1)(2024): 45.

\bibitem{Peral I 1997}Peral, I. Multiplicity of solutions for the $p$-Laplacian. International Center for Theoretical Physics Lecture Notes, Trieste, (1997).

\bibitem{Sun Y 2001}Sun, Y., Wu, S. Combined effects of singular and superlinear nonlinearities in some singular boundary value problems. Journal of Differential Equations 176.2 (2001): 511-531.

\bibitem{Schneider 1989}Schneider, K., Murray, J.D. Mathematical Biology. Biomathematics, Springer-Verlag, Heidelberg, (1990).

\bibitem{Shaomeng 2023}Shao, M., Zhao, L., Yang, Y. Sobolev spaces on locally finite graphs. Proceedings of the American Mathematical Society, 153(02)(2025), 693-708.

\bibitem{Smoller 2012}Smoller, J. Shock waves and reaction¡ªdiffusion equations. Springer Verlag, New York, (1983).

\bibitem{Shao 2024}Shao, M., Yang, Y. Multiplicity and limit of solutions for logarithmic Schr\"{o}dinger equations on graphs. Journal of Mathematical Physics, 65(4)(2024): 17.

\bibitem{Sun 2008}Sun, Y., Li, S. Some remarks on a superlinear-singular problem: Estimates of $\lambda^*$. Nonlinear Analysis, 69(8)(2008): 2636-2650.

\bibitem{Ta V T 2010}Ta, V.T., Elmoataz, A., L\'{e}zoray, O. Nonlocal PDEs-based morphology on weighted graphs for image and data processing. IEEE transactions on Image Processing, 20(6)(2010) : 1504-1516.

\bibitem{Ta V T 2008}Ta, V.T., Elmoataz, A, L\'{e}zoray. O. Partial difference equations on graphs for mathematical morphology operators over images and manifolds. 2008 15th IEEE International Conference on Image Processing, 2008: 801-804.

\bibitem{Xie J 2018}Xie, J., Zhang, X. Infinitely many solutions for a class of fractional impulsive coupled systems with $(p, q) $-Laplacian. Discrete Dynamics in Nature and Society, 2018(1)(2018):9256192.

\bibitem{Yang P 2023}Yang, P., Zhang, X. Existence and multiplicity of nontrivial solutions for a $(p, q) $-Laplacian system on locally finite graphs. Taiwanese Journal of Mathematics, 28(3)(2024): 551-588.

\bibitem{Zhang  2024}Zhang, M., Lin, Y., Yang, Y. Fractional Laplace operator on finite graphs. arxiv: 2403.19987.


















\end{thebibliography}
\end{document}